\documentclass[10pt]{article}
\usepackage{amsmath}
\usepackage{amsfonts}
\usepackage{amssymb}
\usepackage{amsthm}

\begin{document}

\newcommand{\C}{{\mathbb{C}}}
\newcommand{\R}{{\mathbb{R}}}
\newcommand{\Z}{{\mathbb{Z}}}
\newcommand{\N}{{\mathbb{N}}}
\renewcommand{\l}{\left}
\renewcommand{\r}{\right}
\newcommand{\rl}[1]{{\rho_{{#1}}^L}}
\newcommand{\rr}[1]{{\rho_{{#1}}^R}}
\newcommand{\rlp}[3]{\rl{#1}\l( #2,#3 \r)}
\newcommand{\rrp}[3]{\rr{#1}\l( #2,#3 \r)}
\newcommand{\Vl}[3]{{V_{{#1}}^L\l({#2},{#3} \r)}}
\newcommand{\Bl}[3]{B_{{#1}}^L \l({#2},{#3} \r)}
\newcommand{\Vr}[3]{{V_{{#1}}^R\l({#2},{#3} \r)}}
\newcommand{\Br}[3]{B_{{#1}}^R \l({#2},{#3} \r)}
\newcommand{\grad}{\bigtriangledown}
\newcommand{\gl}{\grad_L}
\newcommand{\glx}{\grad_{L,x}}
\newcommand{\gly}{\grad_{L,y}}
\newcommand{\glz}{\grad_{L,z}}
\newcommand{\gr}{\grad_R}
\newcommand{\grx}{\grad_{R,x}}
\newcommand{\gry}{\grad_{R,y}}
\newcommand{\grz}{\grad_{R,z}}
\newcommand{\Xlp}[1]{\dil{X}{L}_{#1}}
\newcommand{\Xrp}[1]{\dil{X}{R}_{#1}}
\newcommand{\Ylp}[1]{\dil{Y}{L}_{#1}}
\newcommand{\Yrp}[1]{\dil{Y}{R}_{#1}}
\newcommand{\sS}{\mathcal{S}}
\newcommand{\sSz}{\sS_0}
\newcommand{\Opl}[1]{{\rm{Op}_L}\l( #1 \r)}
\newcommand{\Opr}[1]{{\rm{Op}_R}\l( #1 \r)}
\newcommand{\Opt}[1]{{\rm{Op}_T}\l( #1 \r)}
\newcommand{\dil}[2]{#1^{\l( #2 \r)}}
\renewcommand{\ker}[1]{{\rm Ker}\l( #1 \r)}
\newcommand{\ip}[2]{\left< #1 , #2 \right>_{L^2}}
\newcommand{\lieg}{\mathfrak{g}}
\newcommand{\lieh}{\mathfrak{h}}
\newcommand{\sA}{\mathcal{A}}
\newcommand{\Ltip}[2]{\l< #1,#2 \r>_{L^2}}
\newcommand{\GtG}{G\times G}
\newcommand{\SztSz}{\sSz\widehat{\otimes}\sSz}
\newcommand{\CM}{C_M^\infty\l(B\r)}
\newcommand{\CMtCM}{C_M^\infty\l( B\r) \widehat{\otimes} C_M^\infty\l( B\r)}
\newcommand{\Kt}{\widetilde{K}}
\newcommand{\Bt}{\widetilde{B}}
\newcommand{\Mt}{\widetilde{M}}
\newcommand{\Et}{\widetilde{E}}
\newcommand{\Zt}{\widetilde{Z}}
\newcommand{\phit}{\widetilde{\phi}}
\newcommand{\chit}{\widetilde{\chi}}
\newcommand{\chitt}{\widetilde{\chit}}
\newcommand{\psit}{\widetilde{\psi}}
\newcommand{\LpOpNorm}[3]{\l\| #1 \r\|_{L^{#2}\l( #3 \r)\circlearrowleft}}
\newcommand{\LpGOpNorm}[2]{\LpOpNorm{#1}{#2}{G}}
\newcommand{\LpGtGOpNorm}[2]{\LpOpNorm{#1}{#2}{\GtG}}
\newcommand{\LpNorm}[3]{\l\| #1 \r\|_{L^{#2}\l( #3 \r)}}
\newcommand{\LpNormG}[2]{\LpNorm{#1}{#2}{G}}
\newcommand{\LpNormGtG}[2]{\LpNorm{#1}{#2}{\GtG}}
\newcommand{\sT}{\mathcal{T}}
\newcommand{\Texp}[1]{{\text{T-exp}}\l( #1 \r)}
\newcommand{\sLp}[1]{\mathcal{L}\l( #1\r)}
\newcommand{\sJ}{\mathcal{J}}
\newcommand{\sM}{\mathcal{M}}
\newcommand{\M}{\sM}
\newcommand{\ML}{\M_L}
\newcommand{\MR}{\M_R}
\newcommand{\Ho}{\mathbb{H}^1}

\newtheorem{thm}{Theorem}[section]
\newtheorem{cor}[thm]{Corollary}
\newtheorem{prop}[thm]{Proposition}
\newtheorem{lemma}[thm]{Lemma}
\newtheorem{conj}[thm]{Conjecture}

\theoremstyle{definition}
\newtheorem{defn}[thm]{Definition}

\theoremstyle{remark}
\newtheorem{rmk}[thm]{Remark}

\title{An Algebra Containing the Two-Sided Convolution Operators}
\author{Brian Street}
\date{}

\maketitle

\begin{abstract}
We present an intrinsically defined algebra of operators containing
the right and left invariant Calder\'on-Zygmund operators
on a stratified group.  The operators in our algebra are pseudolocal
and bounded on $L^p$ ($1<p<\infty$).  This algebra provides an example
of an algebra of singular integrals that falls outside of the
classical Calder\'on-Zygmund theory.

\end{abstract}

\section{Introduction}
Let $G$ be a stratified Lie group.  That is, $G$ is connected, simply connected,
and its Lie Algebra $\lieg$ may be decomposed $\lieg =V_1 \bigoplus \cdots \bigoplus V_m$,
where $\l[ V_1 , V_k \r] = V_{k+1}$ for $1\leq k<m$ and $\l[V_1,V_m \r]=0$.
The Calder\'on-Zygmund theory for left (or right) invariant convolution
operators on $G$ is well-known (see \cite{SteinHarmonicAnalysis}, and 
Section \ref{SectionCZOps} for a review).  Given a distribution
kernel $K$ as in Definition \ref{DefnCZKer} one obtains two ``Calder\'on-Zygmund
singular integral operators'':
$$\Opl{K}f  := f * K$$
$$\Opr{K}f := K*f$$
The operators of the form $\Opl{K}$ form an algebra ($\Opl{K_1} \Opl{K_2} = \Opl{K_1 * K_2}$), are bounded on $L^p$ ($1<p<\infty$), and are pseudolocal.
The same is true for operators of the form $\Opr{K}$.  Also, if we consider:
$$\Opl{K_1} \Opr{K_2} f = \l( K_2 * f\r) * K_1 = K_2 * \l( f *K_1 \r) = \Opr{K_2} \Opl{K_1}f$$
we see that $\Opl{K_1}$ and $\Opr{K_2}$ commute.

Hence, it follows that:
\begin{equation*}
\Opl{K_1} \Opr{K_2} \Opl{K_3} \Opr{K_4} = \Opl{K_1 *K_3} \Opr{K_4 * K_2}
\end{equation*}
and so operators of the form $\Opl{K_1} \Opr{K_2}$ are closed under 
composition.  It is also evident that they are bounded on $L^p$ ($1<p<\infty$)
and are pseudolocal.

The main goal of this paper is to present an algebra of operators, which
contains operators of the form $\Opl{K_1}\Opr{K_2}$, and such that 
the operators in this algebra are
bounded on $L^p$ ($1<p<\infty$), and are pseudolocal.  Moreover,
the algebra will contain the so-called two-sided convolution operators, of
the form:
\begin{equation}\label{EqnIntroTwoSided}
T f (x) = \int K\l( y, z\r) f\l(z^{-1} x y^{-1}\r) \: dy dz 
\end{equation}
where $K$ is a product kernel (see Section \ref{SectionTwoSidedConvolutionOps}).
This algebra provides a naturally occurring example that falls outside
of the classical Calder\'on-Zygmund paradigm.

Operators that fall outside of the classical Calder\'on-Zygmund paradigm
often arise in the construction of parametricies of hypoelliptic
operators which are not maximally subelliptic.  In fact, one of the
original motivations for the present paper was the form of the parametrix
constructed in \cite{StreetAParametrixForKohnsOperator} for Kohn's
example of a sum of squares of complex vector fields, whose
commutators span the tangent space at each point, and such that
the sum of squares is hypoelliptic but not subelliptic (\cite{KohnHypoellipticityAndLossOfDerivatives}).  The parametrix is constructed from compositions of left and right
convolution operators on the three dimensional Heisenberg group, and
is therefore closely related to the algebra discussed in this paper.
It is our hope that the work in this paper will help to motivate
the proper algebras to use in other problems, where the
Calder\'on-Zygmund theory is no longer applicable (for instance,
as in \cite{NagelSteinTheDbarComplexOnDecoupledBoundaries}).

\par\noindent \textbf{Acknowledgements}
\par \noindent This project began as a collaboration with Eli Stein.  During
its early stages, we shared many very interesting conversations on this
subject.  Even as the project progressed, he continued to provide me
with suggestions and encouragement.  Finally,
virtually everything I know about Calder\'on-Zygmund theory, I learned
studying under him during my years as a graduate student.  On all these
counts, I am indebted to him.  I would also like to thank Alex Nagel,
with whom I shared some interesting conversations in the early
stages of this project.

\section{Statement of Results}\label{SectionStatmentOfResults}
Recall, $G$ is a stratified group (for a background on such groups
see \cite{FollandSubellipticEstimatesAndFunctionsSpacesOnNilpotentLieGroups}),
and as such, the lie algebra $\lieg = V_1 \bigoplus \cdots \bigoplus V_m$, 
with the $V_j$ satisfying
the relations
in the introduction.  

Fix a basis $\dil{X}{1}, \ldots, \dil{X}{n}$ for $V_1$, thought of as
elements of the tangent space to the identity in $G$.  Then we can think of each
$\dil{X}{j}$ as either a right invariant or a left invariant vector field,
call them $\dil{X_R}{j}$ and $\dil{X_L}{j}$ respectively.  From here,
we get the left and right gradients:
$$ \gl = \l( \dil{X_L}{1},\ldots,\dil{X_L}{n} \r)$$
$$ \gr = \l( \dil{X_R}{1},\ldots,\dil{X_R}{n} \r)$$
Our definitions will be in terms of $\gl$ and $\gr$, but will not
depend in any essential way on the specific choice of basis of $V_1$.  Throughout
this paper, we will use ordered multi-index notation.  Thus, for a 
finite sequence $s$ of numbers $1,\ldots, n$, we define $\l| s\r|$
to be the length of the sequence, and $\gl^s$, $\gr^s$ in the
usual way.  So that, for instance:
$$\gl^{\l( i,j,k \r)} = \dil{X_L}{i} \dil{X_L}{j} \dil{X_L}{k}$$
and $\l| \l( i,j,k \r) \r|=3$.

For $\epsilon\geq 0$, let $\rl{\epsilon}$ denote the Carnot-Carath\'eodory
distance on $G$ associated to the vector fields $\l\{ \gl, \epsilon \gr\r\}$
and $\rr{\epsilon}$ the one associated to the vector fields $\l\{ \gr, \epsilon \gl\r\}$ (see Section \ref{SectionCCBackground} and references therein for background on such metrics).
Let $\Bl{\epsilon}{x}{\delta}$ denote the ball centered at $x$ of radius
$\delta$ in the $\rl{\epsilon}$ metric, and $\Vl{\epsilon}{x}{\delta}$ its
volume.  Similarly, we define $B^R_\epsilon$ and $V^R_\epsilon$ in terms of $\rr{\epsilon}$.

\begin{defn}\label{DefnBumpFunc}
For $r_R\geq r_L>0$, we say $\phi\in C^\infty\l( G\r)$ is a normalized $r_L,r_R$ bump
function of order $M$ centered at $x\in G$ if $\phi$ is
supported in $\Bl{\frac{r_L}{r_R}}{x}{\frac{1}{r_L}}$, and 
$\forall \l| \alpha \r| + \l|\beta\r|\leq M$,
\begin{equation*}
\l| \gl^\alpha \gr^\beta \phi \r| \leq \frac{r_L^{\l| \alpha \r|}r_R^{\l| \beta\r|}}{\Vl{\frac{r_L}{r_R}}{x}{\frac{1}{r_L}}   }
\end{equation*}
When $r_L\geq r_R$, we replace $\Vl{\frac{r_L}{r_R}}{x}{\frac{1}{r_L}}$ with $\Vr{\frac{r_R}{r_L}}{x}{\frac{1}{r_R}}$.
\end{defn}

We define (for $0<r_L\leq r_R$)
\begin{equation*}
\begin{split}
B& \l( r_L, r_R, N_L, N_R, m, x, y \r) \\
 &= r_L^{N_L} r_R^{N_R} \l( 1+ r_L \rlp{\frac{r_L}{r_R}}{x}{y} \r)^{-m} \frac{1}{\Vl{\frac{r_L}{r_R}}{x}{\frac{1}{r_L}+\rlp{\frac{r_L}{r_R}}{x}{y}}}
\end{split}
\end{equation*}
and when $r_R\geq r_L$, we reverse the roles of $r_L$ and $r_R$ and of the
left and right vector fields.  Before
we define our algebra rigorously, let us write the definition
while being a little loose with quantifiers.
We say that $T\in \sA$ 
if for every $m\geq 0$, and for every $\phi_{\dil{r_L}{1},\dil{r_R}{1}}^x$ normalized $\dil{r_L}{1},\dil{r_R}{1}$ bump functions centered at $x$
and every $\phi_{\dil{r_L}{2},\dil{r_R}{2}}^y$
 normalized $\dil{r_L}{2},\dil{r_R}{2}$
bump functions centered at $y$ (we suppress the order for the moment),
we have:
\begin{equation*}
\begin{split}
&\l|\ip{\phi_{\dil{r_L}{1},\dil{r_R}{1}}^x}{\gl^{\alpha_1}\gr^{\beta_1}T\gl^{\alpha_2}\gr^{\beta_2} \phi_{\dil{r_L}{2},\dil{r_R}{2}}^y } \r|\\
&\quad\quad\leq C B\l( \dil{r_L}{1}\wedge \dil{r_L}{2} ,\dil{r_R}{1}\wedge\dil{r_R}{2}, \l|\alpha_1 \r| + \l| \alpha_2 \r|, \l|\beta_1\r|+\l|\beta_2\r|, m,x,y \r)
\end{split}
\end{equation*}
where $\l|\alpha_1 \r| + \l| \alpha_2 \r|$ and $\l|\beta_1\r|+\l|\beta_2\r|$
must be sufficiently large depending on $m$, and $C$ is uniform in the
choice of normalized bump function, $\dil{r_L}{1},\dil{r_R}{1},\dil{r_L}{2},\dil{r_R}{2}$, and in $x,y$.  Here, and in the rest of the paper, $a\wedge b$ denotes the minimum of $a$ and $b$, while $a\vee b$ denotes the maximum.  Rigorously, our definition is:

\begin{defn}\label{DefnsA}
We define $\sA$ 
to be the set of those operators $T:\sS\l( G\r) \rightarrow \sS\l(G\r)'$, such that 
for all $m\geq 0$,
there exists $N_0$, such that for all $N_L, N_R\geq N_0$ 
there exists $C>0$, $M\in \N$ such that for all
$x,y\in G$, all $\dil{r_L}{1},\dil{r_L}{2}, \dil{r_R}{1},\dil{r_R}{2} >0$, all
$\phi_{\dil{r_L}{1},\dil{r_R}{1}}^x$ normalized $\dil{r_L}{1},\dil{r_R}{1}$ bump functions centered at $x$ of order $M$, all $\phi_{\dil{r_L}{2},\dil{r_R}{2}}^y$ normalized
$\dil{r_L}{2},\dil{r_R}{2}$ bump functions centered at $y$ of order $M$, 
and all $\l| \alpha_1\r|+ \l| \alpha_2\r| = N_L$, 
$\l|\beta_1 \r| +\l|\beta_2\r|=N_R$, we have:
\begin{equation}\label{EqnBoundsADefn}
\begin{split}
&\l|\ip{\phi_{\dil{r_L}{1},\dil{r_R}{1}}^x}{\gl^{\alpha_1}\gr^{\beta_1}T\gl^{\alpha_2}\gr^{\beta_2} \phi_{\dil{r_L}{2},\dil{r_R}{2}}^y } \r|\\
&\quad\quad\leq C B\l( \dil{r_L}{1}\wedge \dil{r_L}{2} ,\dil{r_R}{1}\wedge\dil{r_R}{2}, N_L, N_R, m,x,y \r)
\end{split}
\end{equation}
\end{defn}

\begin{rmk}
We will see \it a posteriori \rm that $N_0=Q+m+1$ will work.  See
Remark \ref{RmkSmallN}.
\end{rmk}

\par\noindent We will show:
\begin{enumerate}
\item The operators in $\sA$ are the same as those in $\sA'$ (defined below; see Section \ref{SectionsAEquivsAp}).
\item The operators in $\sA'$ extend uniquely to bounded operators on
$L^p$, $1<p<\infty$ (Section \ref{SectionLpBounded}).
\item If $T\in \sA$, then $T^{*}\in \sA$, where $T^{*}$ denotes the 
$L^2$ adjoint of $T$ (Remark \ref{RmkDenfsAreSymmetric}).
\item The operators in $\sA'$ form an algebra (Remark \ref{RmksApAlgebra}).
\item The operators in $\sA'$ are pseudolocal (Section \ref{SectionPseudolocal}).
\item Two-sided convolution operators (and therefore the right and left
Calder\'on-Zygmund operators) are in $\sA'$ (Corollary \ref{CorTwosidedAreProduct}).
\end{enumerate}

Our main technical result is that the operators
in $\sA$ are the same as those in $\sA'$.
To define $\sA'$, we need a preliminary definition.

\begin{defn}\label{DefnElemKernels}
We say that $\phi\l( x,z\r)\in C^\infty\l( G\times G\r)$
is an $r_L,r_R$ elementary kernel if, for every $m$ and every
$\alpha_1, \beta_1, \alpha_2,\beta_2$, there exists a $C=C\l( m, \alpha_1, \alpha_2,\beta_1,\beta_2\r)$ such that
\begin{equation}\label{EqnEleZero}
\begin{split}
&\l| \glx^{\alpha_1} \glz^{\alpha_2} \grx^{\beta_1} \grz^{\beta_2} \phi \l(x,z\r)\r| 
\\&\quad \leq C B\l( r_L,r_R, \l|\alpha_1\r|+\l|\alpha_2\r|, \l|\beta_1\r|+\l|\beta_2\r|, m, x,z \r)
\end{split}
\end{equation}
and, for every $N_1,N_2,N_3,N_4\in \N$, and every $\l| \alpha_1\r| =N_1$,
$\l|\alpha_2\r| = N_2$, $\l| \beta_1\r| =N_3$, $\l| \beta_2\r| =N_4$,
there exist functions $\psi_{\alpha_1,\alpha_2,\beta_1,\beta_2}\in C^\infty\l( G\times G\r)$ such that
\begin{equation*}
\phi = r_L^{-N_1-N_2} r_R^{-N_3-N_4} \sum_{\alpha_1,\alpha_2,\beta_1,\beta_2} \glx^{\alpha_1} \glz^{\alpha_2} \grx^{\beta_1} \grz^{\beta_2} \psi_{\alpha_1,\alpha_2,\beta_1,\beta_2}
\end{equation*}
and the $\psi$ satisfy (\ref{EqnEleZero}) with different constants.
Finally, we say $E$ is an $r_L,r_R$ elementary operator if the Schwartz kernel of $E$ 
is an elementary kernel.
\end{defn}

For each $r_L,r_R$, Definition \ref{DefnElemKernels} implicitly defines
a family of seminorms of the elementary kernels (ie, the least possible
$C$ in (\ref{EqnEleZero}), and the least possible $C$ obtained from all choices
of $\psi$, etc.).  If $L:C^\infty\l( G\times G\r)\rightarrow C^\infty\l( G\times G\r)$ is a linear map
that takes $r_L,r_R$ elementary kernels to $r_L, r_R$ elementary kernels,
continuously, it makes sense to ask if it does so \it uniformly \rm in
$r_L,r_R$, since we may order the semi-norms consistently as $r_L$
and $r_R$ vary.

\begin{defn}\label{DefnsAp}
We define $\sA'$ to be those operators $T:\sSz\l(G\r)\rightarrow \sSz\l( G\r)$ 
such that for each $r_L, r_R$, and every $E$ an
$r_L,r_R$ elementary operator, $T E$ is an
$r_L, r_R$ elementary operator, and this map is uniformly continuous in $r_L, r_R$.
Here $\sSz$ is the set of Schwartz functions, all of whose moments vanish.
\end{defn}
\begin{rmk}
The operators in $\sA'$ are \it a priori \rm defined only on $\sSz$.
To see that they are the same as those in $\sA$, we first extend
them as bounded operators on $L^2$, and then prove that the
extended operator is in $\sA$.
\end{rmk}

\begin{rmk}\label{RmksApAlgebra}
It is evident that if $T_1,T_2\in \sA'$  
then $T_1T_2\in \sA'$.  We will show that the operators
in $\sA'$ are the same as those in $\sA$, and therefore $\sA$ forms
an algebra.
\end{rmk}

The operators in $\sA$ may be thought of as
``smoothing of order $0$.''  In Section \ref{SectionClosingRemarks}
we define the analogous concept of operators which are smoothing
of other orders.  In Section \ref{SectionClosingRemarks} we
also discuss an alternative to Definition \ref{DefnsA}, and
why a definition like Definition \ref{DefnsA} seems to be necessary.

\begin{rmk}\label{RmkDenfsAreSymmetric}
Definitions \ref{DefnsA} and \ref{DefnElemKernels} may not seem to be
symmetric (eg, if $T\in \sA$ is $T^{*}\in \sA$? and if $E$ is an
$r_L,r_R$ elementary kernel, is $E^{*}$?), however they are.  Indeed,
despite the fact that $B\l( \cdot, \cdot, \cdot,\cdot,\cdot, x,y \r)$
is not symmetric in $x$ and $y$, it follows from the results in Section \ref{SectionCCBackground} that there exists a $C>0$ such that:
\begin{equation*}
\frac{1}{C} B\l( \cdot,\cdot,\cdot,\cdot,\cdot, y,x\r) \leq B\l(\cdot,\cdot,\cdot,\cdot,\cdot, x,y\r)\leq C B\l( \cdot,\cdot,\cdot,\cdot,\cdot, y,x\r) 
\end{equation*}
\end{rmk}

Some words on notation.  When we refer to the ``unit ball'', we are
always referring to the set $\l\{ x : \l\| x \r\|<1 \r\}$, where
$\l\| \cdot \r\|$ is defined in Section \ref{SectionCZOps}.  $A\lesssim B$
will always mean $A\leq CB$, where $C$ is some constant, independent
of any relevant parameters, and $A\approx B$ means $A\lesssim B$ and
$B\lesssim A$.  Sometimes we will have a sum of positive numbers
of the form
$$\sum_{n\geq 0} a_n$$
and we will have
$$\sum_{n\geq 0} a_n\lesssim \sum_{n\geq 0} r^na_0$$
for some $r$, $0<r<1$.  In this case we will say the series $\sum_{n\geq 0} a_n$ falls
off geometrically or even ``is geometric,'' and we will use the
fact that in this case $\sum_{n\geq 0} a_n\approx a_0$.

\section{Calder\'on-Zygmund Operators}\label{SectionCZOps}
In this section, we will remind the reader of the standard theory
of Calder\'on-Zygmund convolution operators on $G$.  Our goal is three-fold:  first
to fix notation, second to present these concepts in a few different ways,
each of which will be useful in understanding our more complicated algebra,
and finally we will need these characterizations to show that these
Calder\'on-Zygmund operators are in our algebra.

Recall, $G$ is a stratified group, and so, as in the introduction,
the Lie algebra $\lieg=V_1\bigoplus \cdots \bigoplus V_n$, where
the $V_j$ satisfy the relations in the introduction.  The exponential
map $\exp:\lieg \rightarrow G$ is a diffeomorphism.  We define
dilations of $\lieg$, which for $r>0$ are given by $r\cdot X = r^j X$ for
$X\in V_j$.  These dilations induce automorphisms of $G$ by
$re^X = e^{r\cdot X}$.  If we identify $G$ with $\lieg$ via the
exponential map, Lebesgue measure becomes Haar measure for $G$, and
$d \l( rx\r) = r^Q dx$ for some $Q\in \N$.  We call $Q$ the ``homogeneous
dimension'' of $G$.
For a function $\phi:G\rightarrow \C$ and $r>0$, we define
$\dil{\phi}{r}\l(x\r) = r^Q \phi\l(rx\r)$. 
Let $\l\| \cdot \r\|: G\rightarrow \R^{+}$ be a smooth homogeneous norm.
See \cite{FollandSubellipticEstimatesAndFunctionsSpacesOnNilpotentLieGroups}
for a more in depth discussion.

For a background on the material presented here, see \cite{SteinHarmonicAnalysis}
and \cite{NagelRicciSteinSingularIntegralsWithFlagKernels}.  Indeed, we will
be following the presentation of ``product kernels'' from
\cite{NagelRicciSteinSingularIntegralsWithFlagKernels} later in
this section.

\begin{defn}\label{KNormalizedBumpFunc}
A $k$-normalized bump function on $G$ is a $C^k$ function supported
on the unit ball with $C^k$ norm bounded by $1$.  The definitions that
follow turn out to not depend in any essential way on $k$, and
so we shall speak of \it normalized bump functions\rm, thereby suppressing
the dependence on $k$.
\end{defn}

\begin{defn}\label{DefnCZKer}
A Calder\'on-Zygmund kernel on $G$, is a distribution $K$ on $G$, which
coincides with a $C^\infty$ function away from $0$, and satisfies:
\begin{enumerate}
\item (Differential inequalities) For each ordered multi-index $\alpha$, there
is a constant $C_\alpha$ so that
\begin{equation*}
\l| \gl^\alpha K \l( x\r) \r|\leq C_\alpha \l\| x \r\|^{-Q-\l|\alpha\r|} 
\end{equation*}
one may, equivalently, use $\gr$ in place of $\gl$.
\item (Cancellation conditions) Given any normalized bump function $\phi$,
and any $r>0$,
\footnote{We will abuse notation and write the pairing between distributions and test functions as an integral.}
$$\int K\l( x\r) \phi\l(rx\r) dx$$
is bounded independent of $\phi$ and $r$.
\end{enumerate}
\end{defn}

\begin{prop}\label{PropDecompOfCZ}
Let $K$ be a distribution on $G$.  Then, $K$ is a Calder\'on-Zygmund
kernel if and only if there exists a sequence $\l\{ \phi_j \r\}_{j\in \Z}\subset \sSz$,
forming a bounded subset of $\sSz$, such that
$$K= \sum_{j\in \Z} \dil{\phi_j}{2^j}$$
where this sum is taken in distribution (any such sum converges in distribution).  In this case,
\begin{equation}\label{EqnOpDecompCZ}
\Opl{K} = \sum_{j\in \Z} \Opl{\dil{\phi_j}{2^j}}
\end{equation}
where this sum is taken in the strong operator topology as bounded operators on $L^p$ ($1<p<\infty$).
In particular, $\Opl{K}$ is a bounded operator on $L^p$ ($1<p<\infty$).
In addition, (\ref{EqnOpDecompCZ}) converges in the topology of
bounded convergence as operators $\sSz\rightarrow\sSz$.
All of the above can be done uniformly over a bounded subset of 
Calder\'on-Zygmund kernels.  All of the above holds for $\Opr{K}$ as well.
\end{prop}
\begin{proof}
This result is essentially contained in the proofs of Theorem 2.2.1, Theorem 2.6.1, and Proposition 2.7.1
of \cite{NagelRicciSteinSingularIntegralsWithFlagKernels}.  
The only part not appearing in that paper is the convergence in the
topology of bounded convergence.  This follows in a manner completely
analogous to Theorem \ref{ThmDecompOfsAp}.
We
leave the details to the interested reader.
\end{proof}

\begin{thm}\label{ThmCZLikesAp}
Let $T:\sSz\l( G\r) \rightarrow C_0^\infty \l( G\r)'$.  Then,
$T=\Opl{K}$ (when restricted to $\sSz$), where $K$ is a Calder\'on-Zygmund kernel, if and only if 
for every $\phi\in \sSz$ and every $r>0$,
\begin{equation*} 
\Opl{K} \Opl{\dil{\phi}{r}} = \Opl{\dil{\psi_r}{r}}
\end{equation*}
where $\psi_r\in \sSz$, and as $\phi$ ranges over a bounded set in $\sSz$,
and $r$ ranges over $r>0$, we have that $\psi_r$ ranges over a bounded set in $\sSz$.
\end{thm}

We defer the proof to Section \ref{ProofOfCZLikesAP}. Theorem \ref{ThmCZLikesAp} should be interpreted in the following way:
we think of operators of the form $\Opl{\dil{\phi}{r}}$, with $\phi\in \sSz$, as our ``$r$ elementary
operators'' in analogy with Definition \ref{DefnElemKernels}.  
Theorem
\ref{ThmCZLikesAp} simply says that $T$ is a Calder\'on-Zygmund operator
if and only if
composition with $T$ takes
$r$ elementary operators to $r$ elementary operators uniformly,
in analogy with Definition \ref{DefnsAp}.

We now turn to an equivalent way of considering Calder\'on-Zygmund
operators that is analogous to Definition \ref{DefnsA}.  Let
$K$ be a Calder\'on-Zygmund kernel, let $T=\Opl{K}$, and let $\phi, \psi$
be normalized bump functions.  Define:
$$\phi^x_r\l( y \r) = r^Q \phi \l( r\l( y^{-1}x\r)\r)$$
and similarly for $\psi^x_r$.  The cancellation condition of
Definition \ref{DefnCZKer} shows that $\l| T \phi^{x_0}_r \r|\lesssim r^Q$,
on $\l\| x^{-1} x_0\r\|\leq 2r^{-1}$.  Combining this with the
growth condition, one sees:
\begin{equation}\label{EqnSupCancelZ}
\l| T \phi_r^{x_0} \l(x\r) \r| \lesssim
\begin{cases}
r^Q & \text{if $\l\| x^{-1}x_0 \r\|\leq r^{-1}$,}\\
\l\| x^{-1}x_0 \r\|^{-Q} & \text{if $\l\| x^{-1}x_0 \r\|\geq r^{-1}$
.}
\end{cases}
\end{equation}
Conversely, it is clear that Equation (\ref{EqnSupCancelZ}) implies
the cancellation condition of Definition \ref{DefnCZKer}.  To see that
it also implies the growth condition (where there are no derivatives involved), merely choose $\phi$ so that
as $r\rightarrow \infty$, $\phi_r^{x_0}\rightarrow \delta_{x_0}$.

Now suppose $s>r$, we see from (\ref{EqnSupCancelZ}),
\begin{equation}\label{EqnSupCancelZt}
\begin{split}
\l|\Ltip{\psi^{x_1}_s}{T\phi^{x_0}_r}\r| &\lesssim
\begin{cases}
r^Q & \text{if $\l\| x_1^{-1}x_0 \r\| \leq r^{-1}$,}\\
\l\| x_1^{-1} x_0 \r\|^{-Q} & \text{if $\l\| x_1^{-1}x_0 \r\|\geq r^{-1}$}
\end{cases}\\
&\lesssim \frac{1}{\l( \frac{1}{r}+\l\| x_1^{-1}x_0 \r\| \r)^Q}\\
&\lesssim \frac{1}{V\l( x_0, \frac{1}{r}+ \l\|x_1^{-1}x_0 \r\| \r)}
\end{split}
\end{equation}
where $V\l( x, \delta\r)$ denotes the volume of $\l\{ y: \l\| y^{-1}x \r\|<\delta \r\}$.  So we see that (\ref{EqnSupCancelZt}) is implied by 
(\ref{EqnSupCancelZ}).  The converse is true as well, as can be seen
by taking $\psi$ such that $\psi^{x_1}_s \rightarrow \delta_{x_1}$
as $s\rightarrow \infty$.  From these considerations, the following theorem follows
easily:

\begin{thm}
Suppose $T:\sS \rightarrow C_0^\infty\l( G\r)'$, and is
left invariant.  Then, as operators on $\sSz$, $T=\Opl{K}$ where $K$ is a Calder\'on-Zygmund
kernel if and only if for all $m$, and all $\phi$ and $\psi$ normalized
bump functions, and all $\alpha,\beta$ ordered multi-indicies such that
$\l| \alpha \r| +\l|\beta\r|$ is sufficiently large depending on $m$,
we have:
\begin{equation}\label{EqnCZBothCancel}
\l| \Ltip{\psi_s^{x_1}}{\gl^\alpha T\gl^\beta \phi_r^{x_0}} \r| \lesssim
\l( 1+ t\l\| x_1^{-1}x_0 \r\| \r)^{-m} \frac{t^{\l|\alpha\r|+\l|\beta\r|}}{V\l( x_0, \frac{1}{t}+\l\| x_1^{-1} x_0 \r\| \r)}
\end{equation}
where $t=s\wedge r$, and the bound is uniform in $s$, $r$, $x_0$, $x_1$,
and choices of normalized bump functions.  This is analogous to
Definition \ref{DefnsA}.
\end{thm}

\begin{rmk}
We will see later that using the cancellation condition on both
sides simultaneously as in (\ref{EqnCZBothCancel}) seems to be
necessary in our situation.  See Section \ref{SectionClosingRemarks}.
\end{rmk} 

\subsection{Proof of Theorem \ref{ThmCZLikesAp}}\label{ProofOfCZLikesAP}

\begin{lemma}\label{LemmaDecompForsSz}
Given $\phi\in \sSz$, there exists $\psi\in \sSz^n$ such that
\begin{equation}\label{EqnDefnsSz}
\phi = \gl \cdot \psi
\end{equation}
Moreover, for each continuous semi-norm $\l| \cdot \r|$ on $\sSz$,
the infimum over all such $\psi$ of $\sum \l| \psi_j \r|$ is a
continuous semi-norm on $\sSz$.

For $\phi \in \sS$, we say $\phi\in \sS_{\l(m\r)}$ if
$\phi$ can be written as in (\ref{EqnDefnsSz}) with $\psi\in {\sS}_{\l(m-1\r)}^n$;
where $\sS_{\l(0\r)}=\sS$.  Then, $\sSz = \cap_{m} \sS_{\l(m\r)}$.
\end{lemma}
\begin{proof}
This lemma is well known.
\end{proof}

\begin{lemma}\label{LemmaCZCancel}
Suppose $\phi_1, \phi_2\in \sSz$, then 
$$\dil{\phi_1}{2^j}*\dil{\phi_2}{2^k} = 2^{-\l| j-k\r|} \dil{\psi}{2^{l}}$$
where $l$ can be either $j\wedge k$ or $j\vee k$, $\psi\in \sSz$, and when $\phi_1, \phi_2$ vary over a bounded set of $\sSz$, and $j,k$ vary over $\Z$, $\psi$
varies over a bounded set of $\sSz$.
\end{lemma}
\begin{proof}
We prove the result first in the case $l= j\wedge k$ and $j\leq k$, the case $k<j$ is similar.  Note
that:
\begin{equation*}
\dil{\phi_1}{2^j} *\dil{\phi_2}{2^k} (x) = 2^{\l(j+k\r)Q} \int \phi_1 \l( 2^j \l( xy^{-1} \r) \r) \phi_2 \l(2^k y\r) dy
\end{equation*}
and so replacing $x$ with $2^{-j} x$ and doing a change of variables $u = 2^{j} y^{-1}$,
we see that we can just prove the lemma for $j=0$, $k\geq 0$.

Writing $\phi_2 = \gr \cdot \phi_2'$ as in Lemma \ref{LemmaDecompForsSz}
(using right invariant vector fields, instead of left),
we see that
\begin{equation*}
\begin{split}
\phi_1 * \dil{\phi_2}{2^k} &= 2^{kQ-k} \int \phi_1 \l( xy^{-1} \r) \gr\cdot \phi_2'(2^k y) dy \\
&= 2^{kQ-k} \int \l(\gr \phi_1\r) \l( xy^{-1} \r) \cdot \phi_2' (2^ky) dy 
\end{split}
\end{equation*}
Repeating this process $Q$ more times, we see:
\begin{equation}\label{EqnAfterByParts}
\phi_1 * \dil{\phi_2}{2^k} = 2^{-k} \sum_m \int \psi_{1,m} \l( xy^{-1} \r)  \psi_{2,m}(2^k y) dy
\end{equation}
where $m$ ranges over a finite set, and the $\psi_{1,m}, \psi_{2,m}\in \sSz$
(and range over a bounded set as $\phi_1, \phi_2$ do).

To see that $\psi_1 * \l( 2^{-kQ}\dil{\psi_2}{2^k}\r)$ is rapidly decreasing (independent
of $k$), we now need to merely apply the fact that if $f_1$ and
$f_2$ are two bounded rapidly decreasing functions, then
$f_1 * f_2$ is rapidly decreasing (note that $\psi_2 \l( 2^k y \r)$
decreases faster than $\psi_2\l(y\r)$).

Since, $\gr \l( \phi_1 * \dil{\phi_2}{2^k}\r) = \l( \gr \phi_1 \r) * \phi_2$,
we see that $\phi_1 * \dil{\phi_2}{2^k}\in \sS$.  To see it is really in
$\sSz$, we use the fact that $\phi_1 = \gr \cdot \phi_1'$, and therefore,
\begin{equation*}
\phi_1 * \dil{\phi_2}{2^k} = \l( \gr\cdot \phi_1' \r) * \dil{\phi_2}{2^k} 
=\gr \cdot \l( \phi_1' * \dil{\phi_2}{2^k} \r)
\end{equation*}
repeating this process and applying Lemma \ref{LemmaDecompForsSz},
completes the proof of the case $l=j\wedge k$.

Turning to the case when $l=j\vee k$, we again assume $j=0$ and 
now assume $k> 0$,
the other cases being similar.  A computation similar to the one leading
up to (\ref{EqnAfterByParts}) shows that
\begin{equation}\label{EqnAfterByParts2}
\phi_1 * \dil{\phi_2}{2^k} = 2^{-Nk} \sum_m \int \psi_{1,m} \l( xy^{-1} \r)  \psi_{2,m}(2^k y) dy
\end{equation}
for any fixed $N$, where the $\psi_{1,m},\psi_{2,m}$ are as above.  Thus,
if one wishes to show that
\begin{equation}\label{EqnBoundConvo1}
\l\| 2^k x \r\|^N \phi_1 * \dil{\phi_2}{2^k}
\end{equation}
is bounded by $C2^{\l( Q-1\r)k}$, one merely needs to apply (\ref{EqnAfterByParts2}) and use the fact that 
$$\int \psi_{1,m}\l(xy^{-1}\r) \psi_{2,m}\l( 2^k y \r)$$ 
is rapidly decreasing, as shown above.  In fact, we get the stronger result
that (\ref{EqnBoundConvo1}) is bounded independent of $k$.  Derivatives work as before,
yielding the result.

\end{proof}

\begin{proof}[Proof of Theorem \ref{ThmCZLikesAp}]
First, suppose that $T=\Opl{K}$ where $K$ is a Calder\'on-Zygmund
kernel.  We prove the result for $r=2^l$ for some $l\in \Z$.  The
more general result follows from this by moving $r$ to the
closest such $2^l$, via replacing $\phi$ by $\phi\l( \frac{r}{2^l} x\r)$.

Applying Proposition \ref{PropDecompOfCZ}, we may write:
\begin{equation*}
\Opl{K} = \sum_{k\in \Z} \Opl{\dil{\psi_k}{2^k}}
\end{equation*}
with the $\psi_k\in \sSz$ uniformly in $k$ (even as $j$ varies).
We now apply Lemma \ref{LemmaCZCancel} to see:
\begin{equation*}
\begin{split}
\Opl{K} \Opl{\dil{\phi}{2^j}} &= \sum_{k\in \Z} \Opl{\dil{\psi_k}{2^k}} \Opl{\dil{\phi}{2^j}}\\
&= \sum_{k\in \Z} \Opl{\dil{\phi}{2^j}*\dil{\psi_k}{2^k}} \\
&= \sum_{k\in \Z} 2^{-\l| k-j \r|}\Opl{\dil{\phi_k}{2^j}} \\
&=: \Opl{\dil{\phi_0}{2^j}}
\end{split}
\end{equation*}
where $\phi_0\in \sSz$ ranges over a bounded set as the relevant
variables change.

Conversely, suppose $T$, satisfies the conditions on the theorem. 
We wish to show that $T=\Opl{K}$, where $K$ is a Calder\'on-Zygmund
operator.  We know that $I$ (the identity) is a
Calder\'on-Zygmund operator, and therefore, by Proposition \ref{PropDecompOfCZ}:
\begin{equation*}
I = \sum_{k\in Z} \Opl{\dil{\phi_k}{2^k}}
\end{equation*}
with the convergence in the topology of bounded convergence as operators
$\sSz\rightarrow\sSz$.  Hence,
\begin{equation*}
\begin{split}
T &= T I\\
&= \sum_{k\in \Z} T \Opl{\dil{\phi_k}{2^k}}\\
&= \sum_{k\in \Z} \Opl{\dil{\psi_k}{2^k}}
\end{split}
\end{equation*}
with the $\psi_k$ forming a bounded subset of $\sSz$.  Thus,
by Proposition \ref{PropDecompOfCZ}, $T$ is a Calder\'on-Zygmund
operator.
\end{proof}

\subsection{Two-Sided Convolution Operators}\label{SectionTwoSidedConvolutionOps}
In this section, we turn to the definition of so-called
``product kernels,'' and use that to define the relevant
``two-sided convolution operators'' that we will be studying
(see (\ref{EqnIntroTwoSided})).  Our main reference for
product kernels is \cite{NagelRicciSteinSingularIntegralsWithFlagKernels},
and we refer the reader there for any further reading.

\begin{defn}
A product kernel on $\GtG$ is a distribution $K\l(x,y\r)$ on $\GtG$, which
coincides with a $C^\infty$ function away from $\l\{ x=0\r\} \cup \l\{ y=0\r\}$
and which satisfies:
\begin{enumerate}
\item (Differential inequalities)  For each ordered multi-indicies $\alpha_1,\alpha_2$,
there is a constant $C=C\l( \alpha_1,\alpha_2\r)$ such that:
\begin{equation*}
\l| \glx^{\alpha_1} \gly^{\alpha_2} K\l( x,y\r) \r|\leq C \l| x\r|^{-Q-\l|\alpha_1\r|} \l| y \r|^{-Q-\l| \alpha_2 \r|}
\end{equation*}
the definition remains unchanged if we replace $\gl$ by $\gr$.
\item (Cancellation conditions) Given any normalized bump function $\phi$
on $G$, and any $R>0$, the distributions:
\footnote{$\int K\l( x,y\r) \phi\l(x \r) dx$ denotes that distribution which,
when paired with the test function $\psi\l( y \r)$, equals $\int K\l( x,y\r) \phi\l(x \r) \psi\l( y\r)dx\: dy$}
\begin{equation*}
K_{\phi,R}\l( x\r) = \int K\l( x,y\r) \phi\l( Ry\r) dy
\end{equation*}
\begin{equation*}
K_{\phi,R}\l( y\r) = \int K\l( x,y\r) \phi\l( Rx\r)dx
\end{equation*}
are Calder\'on-Zygmund kernels, uniformly in $\phi$ and $R$.
\end{enumerate}
For such a kernel we define $\Opl{K}$ and $\Opr{K}$, acting on functions
in $C_0^\infty \l( \GtG\r)$ to be convolution with $K$ on the right and
left over the group $\GtG$, respectively.
\end{defn}

For a function $f:G\rightarrow \C$ and $r_1,r_2>0$, we define:
\begin{equation*}
\dil{f}{r_1,r_2}\l( x,y\r) = r_1^Q r_2^Q f\l( r_1 x, r_2 y\r)
\end{equation*}

\begin{defn}
Let $\SztSz$ denote the set of those functions $f\in \sS\l( \GtG\r)$
such that for every multi-index $\alpha$,
\begin{equation*}
\int x^\alpha f\l( x,y\r) dx =0
\end{equation*}
\begin{equation*}
\int y^\alpha f\l( x,y\r) dy =0
\end{equation*}
\end{defn}

\begin{rmk}\label{RmkNuclear}
We note that $\SztSz$ is nothing more than the tensor product of the
nuclear space $\sSz$ with itself.  This
explains our notation.  See \cite{TrevesTopologicalVectorSpaces} for
a background on tensor products.  We will not use any deep results
about tensor products, however they will provide us with one
small convenience.  Since the above tensor product agrees with the
projective tensor product, when we wish to prove a result about
$f\in \SztSz$, it will often suffice to prove the result for 
$f\l(x,y\r)=\phi_1\l(x\r) \phi_2\l(y\r)$, $\phi_1,\phi_2\in \sSz$.  
We will use this fact freely in the sequel.
In particular, each $f\l( x,y\r)\in \SztSz$ can be decomposed:
\begin{equation}\label{EqnDecompSztSz}
f = \glx \cdot \gry \cdot g
\end{equation}
where $g\in \l(\SztSz\r)^{n^2}$; this can be easily seen
for elementary tensor products by Lemma \ref{LemmaDecompForsSz},
and therefore holds for all elements of $\SztSz$.
\end{rmk}

\begin{prop}
Let $K$ be a distribution on $\GtG$.  Then, $K$ is a product kernel if
and only if there exists a sequence $\l\{ \phi_{j,k} \r\}_{j,k\in \Z} \subset \SztSz$,
forming a bounded subset of $\SztSz$, such that:
\begin{equation*}
K = \sum_{j,k\in \Z} \dil{\phi_{j,k}}{2^j,2^k}
\end{equation*} 
where this sum is taken in distribution (any such sum converges in distribution).  In this case,
\begin{equation*}
\Opl{K} = \sum_{j,k\in \Z} \Opl{\dil{\phi_{j,k}}{2^j,2^k}}
\end{equation*}
where this sum converges in the strong operator topology as maps
$L^p\l( \GtG\r)\rightarrow L^p\l( \GtG\r)$, for ($1<p<\infty$).
In particular, $\Opl{K}$ extends to a bounded operator on $L^p\l( \GtG\r)$, ($1<p<\infty$).  All of the above can be done uniformly for kernels forming
a bounded subset of the product kernels.
A similar result holds for $\Opr{K}$.  
\end{prop}
\begin{proof}
This is essentially contained in the proofs of Theorem 2.2.1, 
Theorem 2.6.1, and Proposition 2.7.1 of \cite{NagelRicciSteinSingularIntegralsWithFlagKernels}.  We leave the details to the interested reader.
\end{proof}

Given a product kernel $K$ on $\GtG$, we may define the
two-sided convolution operator $\Opt{K}$ as:
\begin{equation*}
\Opt{K}f \l( x \r) = \int K\l( y,z \r) f\l( z^{-1} x y^{-1}  \r) dy\: dz
\end{equation*}
To see that this makes sense, for $x$ fixed, after the $z$ integration,
the integrand left over is $O\l( \l\| y\r\|^{-2Q}\r)$ and so converges
absolutely for $y$ large.  For $y$ small, this makes sense using the
usual pairing of distributions and test functions.
Note that,
\begin{equation*}
\Opt{K_1\l( x \r) K_2\l( y\r)} = \Opl{K_1} \Opr{K_2}
\end{equation*}
and so two-sided convolution operators contain the right and left invariant Calder\'on-Zygmund
operators as special cases.

Since our algebra $\sA$ will contain the two-sided convolution operators
corresponding to product kernels, the $L^p$ boundedness ($1<p<\infty$)
of two-sided convolution operators will follow from the $L^p$
boundedness of operators in $\sA$.  However, we will actually
need the $L^p$ boundedness of two-sided convolution operators to
\it prove \rm the $L^p$ boundedness of operators in $\sA$,
and therefore we now turn to proving the $L^p$ boundedness of
two-sided convolution operators. 

We now introduce a formal trick that allows us to consider
two-sided convolution operators on $G$ as convolution operators on 
$\GtG$.  This may be found in \cite{KisilConnectionBetweenTwoSidedAndOneSided}.
For a function $f:G\rightarrow \C$, we may define a new function $Ef:\GtG\rightarrow \C$ by
\begin{equation*}
\l( Ef\r) \l( x,y \r) = f\l( y^{-1} x\r)
\end{equation*}
Then, a simple computation yields that, whenever it makes sense,
\begin{equation*}
E \Opt{K} f = \Opl{\Kt} Ef
\end{equation*}
where $\Kt \l( x, y\r) = K\l( x, y^{-1}\r)$.  Note that
$\Kt$ is a product kernel if and only if $K$ is.
Note also that,
\begin{equation}\label{EqnTwoSidedAlg}
\begin{split}
E \Opt{K_1} \Opt{K_2} &= \Opl{\Kt_1}\Opl{\Kt_2} E \\&= \Opl{\Kt_1 *\Kt_2} E \\&= E \Opt{\widetilde{\Kt_1 * \Kt_2}}
\end{split}
\end{equation}
and so two-sided convolution operators form an algebra (since product kernels
form an algebra under convolution, see \cite{NagelRicciSteinSingularIntegralsWithFlagKernels}).

\begin{lemma}\label{LemmaTwoSidedBounded}
Suppose $K$ is a product kernel and $K$ has compact support.
Then, $\Opt{K}$ extends to a bounded operator on $L^p\l( G\r)$,
and moreover, $$\LpGOpNorm{\Opt{K}}{p} \leq \LpGtGOpNorm{\Opl{\Kt}}{p}$$
\end{lemma}
\begin{proof}
The proof uses transference.  We follow the outline of a similar
argument on page 483 of \cite{SteinHarmonicAnalysis}, where the proof
falls under the heading ``method of descent.''
We suppose that $K\l( x,y\r)$ is supported on $\l\| x \r\|, \l\| y\r\| \leq M$.
Consider, for $R>0$,
\begin{equation}\label{EqnTrans1}
\begin{split}
\LpNormG{\Opt{K}f}{p}^p &= \LpNormG{\l( E\Opt{K}f\r) \l( \cdot,x_2 \r)}{p}^p,\:\:\forall x_2\in G \\
&= \LpNormG{ \l( \Opl{\Kt} Ef\r) \l(\cdot,x_2\r) }{p}^p, \:\:\forall x_2\in G\\
&= \frac{1}{R^Q} \int_{\l\| x_2 \r\|\leq R} \LpNormG{\l(\Opl{\Kt}Ef\r) \l( \cdot, x_2 \r)}{p}^p dx_2
\end{split}
\end{equation}
But,
\begin{equation*}
\begin{split}
\l(\Opl{\Kt}Ef\r) \l( x_1, x_2\r) &= \int \Kt\l( y^{-1}x_1, z^{-1}x_2 \r) Ef\l(y,z\r) dy\:dz
\end{split}
\end{equation*}
and since in (\ref{EqnTrans1}), we are only considering $\l\| x_2 \r\|\leq R$,
and by the support of $\Kt$ only considering $\l\| z^{-1}x_2\r\|\leq M$, we
see that we are only integrating over $\l\| z \r\|\leq M+R$.  Hence,
in (\ref{EqnTrans1}) we may replace $Ef\l( y,z\r)$ with
$$F\l(y,z\r) := \l(Ef\l( y,z \r)\r) \chi_{\l\{\l\| z\r\|\leq M+R\r\}}$$
Thus,
\begin{equation*}
\begin{split}
&\LpNormG{\Opt{K}f}{p}^p = \frac{1}{R^Q} \int_{\l\| x_2\r\|\leq R} 
\LpNormG{\l(\Opl{\Kt}F\r)\l( \cdot, x_2 \r)}{p}^p dx_2 \\
&\quad\leq \frac{1}{R^Q} \LpNormGtG{\Opl{\Kt}F}{p}^p \\
&\quad\leq \frac{1}{R^Q} \LpGtGOpNorm{\Opl{\Kt}}{p}^p \LpNormGtG{F}{p}^p\\
&\quad\leq \frac{1}{R^Q} \LpGtGOpNorm{\Opl{\Kt}}{p}^p \int \l| \l( Ef \r)\l( y,z \r)  \chi_{\l\{ \l\| z \r\|\leq M+R \r\}}\r|^p dy\:dz\\
& \quad= \frac{1}{R^Q} \LpGtGOpNorm{\Opl{\Kt}}{p}^p \l( M+R \r)^Q \LpNormG{f}{p}^p\\
&\quad\xrightarrow[R\rightarrow\infty]{} \LpGtGOpNorm{\Opl{\Kt}}{p}^p \LpNormG{f}{p}^p
\end{split}
\end{equation*}
Completing the proof.
\end{proof}
\begin{cor}\label{CorTwoSidedBounded}
Let $K$ be a product kernel, then $\Opt{K}$ extends uniquely to a bounded operator
on $L^p\l( G \r)$, $1<p<\infty$, and
\begin{equation*}
\LpGOpNorm{\Opt{K}}{p}\leq \LpGtGOpNorm{\Opl{\Kt}}{p}
\end{equation*}
\end{cor}
\begin{proof}
It suffices to show that for $f\in C_0^\infty \l( G\r)$,
$$\LpNormG{\Opt{K} f}{p} \leq \LpGtGOpNorm{\Opl{\Kt}}{p} \LpNormG{f}{p}$$
and this follows from Lemma \ref{LemmaTwoSidedBounded} and a
simple limiting argument.
\end{proof}

\begin{prop}\label{PropDecompOfTwosided}
Suppose $K$ is a product kernel, and suppose that 
$\l\{ \phi_{j,k} \r\}\subset \SztSz$ is a bounded subset
such that 
$$K = \sum_{j,k\in \Z} \dil{\phi_{j,k}}{2^j,2^k}$$
where this sum is taken in distribution.  Then,
$$\Opt{K} = \sum_{j,k\in \Z} \Opt{\dil{\phi_{j,k}}{2^j,2^k}}$$
where this sum converges in the strong operator topology
on $L^p\l( G\r)$, $1<p<\infty$.
\end{prop}
\begin{proof}
Corollary \ref{CorTwoSidedBounded} tells us that the operators:
$$\sum_{|j|,|k|\leq N} \Opt{\dil{\phi_{j,k}}{2^j,2^k}} = \Opt{ \sum_{|j|,|k| \leq N} \dil{\phi_{j,k}}{2^j,2^k} }$$
are uniformly bounded on $L^p$.  It is easy to see that for
$f\in C_0^\infty$,
$$\sum_{j,k\in \Z} \Opt{\dil{\phi_{j,k}}{2^j,2^k}}f = \Opt{K} f$$
where this sum is taken in distribution.  
Putting these two facts together, and using the fact that
functions of the form $f= X_L X_R g$, $g\in C_0^\infty$, $X_L=\gl^{\alpha}$,
$\l| \alpha \r| =1$, $X_R=\gr^\beta$, $\l|\beta\r| = 1$, span
a dense subset of $L^p$, it suffices to show that:
\begin{equation}\label{EqnStrong1}
\sum_{j,k\in \Z} \Opt{\dil{\phi_{j,k}}{2^j,2^k}}f
\end{equation}
converges in $L^p$, for all such $f$.
We separate
(\ref{EqnStrong1}) into four sums, and we first consider:
\begin{equation}\label{EqnStrong2}
\begin{split}
\sum_{j\geq 0, k\geq 0} \Opt{\dil{\phi_{j,k}}{2^j,2^k}} f &= \sum_{\l| \alpha\r|=1,\l| \beta\r|=1}\sum_{j\geq 0, k\geq 0} \Opt{\dil{ \l( \glx^\alpha \grx^\beta \psi_{j,k}\r)}{2^j,2^k}} f
\\&= \sum_{\alpha,\beta} \sum_{j\geq 0, k\geq 0} 2^{-j-k} \Opt{\dil{\psi_{j,k}}{2^j,2^k}} \gl^\alpha \gr^\beta f
\end{split}
\end{equation}
here we have applied
(\ref{EqnDecompSztSz}) and $\psi_{j,k}$ also depends on $\alpha$, $\beta$ and ranges over a bounded
subset of $\SztSz$ as $\alpha$, $\beta$, $j$, and $k$ vary.  Thus, Corollary \ref{CorTwoSidedBounded}
tells us that (\ref{EqnStrong2}) converges in $L^p$.  We now consider:
\begin{equation}\label{EqnStrong3}
\begin{split}
\sum_{j< 0, k< 0} \Opt{\dil{\phi_{j,k}}{2^j,2^k}} f &= \sum_{j< 0, k< 0} \Opt{\dil{\phi_{j,k}}{2^j,2^k}} X_L X_R g\\
&= \sum_{j< 0, k< 0} \Opt{X_L X_R\dil{\phi_{j,k}}{2^j,2^k}} g\\
&= \sum_{j< 0, k< 0} 2^{j+k} \Opt{\dil{\l( X_L X_R\phi_{j,k}\r)}{2^j,2^k}} g
\end{split}
\end{equation}
and since $X_L X_R\phi_{j,k}$ ranges over a bounded subset of
$\SztSz$, we again have by Corollary \ref{CorTwoSidedBounded}
that (\ref{EqnStrong3}) converges in $L^p$.
Finally, the sums where $j<0, k\geq0$ and $k<0, j\geq0$ follow
from a combination of the two methods above.  We leave the details to
the reader.
\end{proof}

We state, without proof, a result similar to Theorem \ref{ThmCZLikesAp},
that gives a characterization of two-sided convolution operators.  We
leave out to proof for two reasons.  Firstly, it it entirely analogous
to the proof of Theorem \ref{ThmCZLikesAp}.  Secondly, we will not
use this characterization for anything other than motivation for
Definition \ref{DefnsAp}.

\begin{thm}
Suppose $T:\sSz \rightarrow C_0^\infty\l( G\r)'$.  Then,
$T=\Opt{K}$ (when restricted to $\sSz$), where $K$ is a product kernel, if and only if 
for every $\phi\in \SztSz$, and every $r_L,r_R>0$,
\begin{equation*}
T \Opt{\dil{\phi}{r_L,r_R}} = \Opt{\dil{\psi_{r_L,r_R}}{r_L,r_R}}
\end{equation*}
where $\psi_{r_L,r_R}\in \SztSz$, and ranges over a bounded subset of
$\SztSz$ as $\phi$ ranges over a bounded subset and $r_L$ and $r_R$ vary.
\end{thm}

\section{Carnot-Carath\'eodory Distances}\label{SectionCCBackground}
In this section, we review the metrics defined naturally in terms
of a given family of vector fields (often called Carnot-Carath\'eodory
metrics, or sub-Riemannian metrics).  Our main references for this section
are \cite{NagelSteinWaingerBallsAndMetricsDefinedByVectorFields,NagelSteinDifferentiableControlMetricsAndScaledBumpFunctions}
however we will need to restate many of the results from those
papers in a slightly stronger way; though no new proofs will
be required.  The expert in these topics may skip this section (except,
perhaps, for Section \ref{SectionNewDist}), given
the understanding that all the facts we will use about such distances
are true \it uniformly \rm for $\rl{\epsilon},\rr{\epsilon}$ for
$\epsilon\in \l[0,1\r]$, where $\rl{\epsilon},\rr{\epsilon}$ were
defined in the introduction.

Let $\Omega \subset \R^N$ be a connected open set, and let $Y_1, Y_2, \cdots Y_q$
be a list (possibly with repetitions) of real $C^\infty$ vector fields
on $\Omega$.  Associate to each $Y_j$ an integer $d_j\geq 1$, called
the formal degree of $Y_j$.  Following \cite{NagelSteinWaingerBallsAndMetricsDefinedByVectorFields,NagelSteinDifferentiableControlMetricsAndScaledBumpFunctions},
we define:

\begin{defn}[Definition 2.1.1 in \cite{NagelSteinDifferentiableControlMetricsAndScaledBumpFunctions}]\label{DefnHomogType}
The list of vector fields and associated formal degrees $\l\{ \l( Y_j, d_j\r)\r\}$ is said to be of finite
homogeneous type on $\Omega$ if:
\begin{enumerate}
\item For all $1\leq j,k \leq q$,
\begin{equation*}
\l[ Y_j, Y_k \r] = \sum_{d_l\leq d_j+d_k} c_{j,k}^l \l( x\r) Y_l
\end{equation*}
where $c_{j,k}^l \in C^\infty \l( \Omega \r)$.
\item At each point $x\in \Omega$, $\l\{ Y_1\l(x\r), \ldots, Y_q\l(x\r) \r\}$
spans the tangent space at $x$.
\end{enumerate}
\end{defn}

A fundamental example of Definition \ref{DefnHomogType} (and the
only one we will use) is given by a set of vector fields $X_1,\ldots, X_n$
on $\Omega$ such that all the iterated commutators of length at most
$m$ span the tangent space at each point.  We take $Y_1, \ldots, Y_q$
to be a list of all these commutators, with the degree of $Y_j$
being the length of the commutator from which it arises.

\begin{defn}[Definition 2.1.2 of \cite{NagelSteinDifferentiableControlMetricsAndScaledBumpFunctions}]
Let $Y=\l\{ \l( Y_1, d_1\r) ,\ldots, \l( Y_q, d_q\r)\r\}$ be a list of vector fields and formal degrees which are
of finite homogeneous type on $\Omega$.  For each $\delta>0$ let
$C\l( \delta ,Y \r)$ denote the set of absolutely continuous curves
$\phi:\l[ 0,1 \r]\rightarrow \Omega$ which satisfy:
\begin{equation*}
\phi'\l( t\r) = \sum_{j=1}^q a_j\l( t\r) Y_j\l( \phi\l( t\r)\r)\quad \text{with} \quad \l| a_j\l(t\r) \r|\leq \delta^{d_j}
\end{equation*}
for almost all $t\in \l[0,1\r]$.  For $x,y\in \Omega$, set
\begin{equation*}
\rho_Y\l( x,y\r) = \inf \l\{ \delta>0 | \bigg(\exists \phi\in C\l(\delta\r) \bigg) \bigg( \phi\l(0\r)=x, \phi\l(1\r)=y\bigg) \r\}
\end{equation*}
The function $\rho_Y$ is called the control metric on $\Omega$, generated
by $Y$.
\end{defn}

\begin{rmk}\label{RmkLeftInvDist}
If we take $Y$ to be that list of vector fields generated by
the left invariant vector fields of order $1$ (ie, we take
$\gl^\alpha$, $\l|\alpha\r|=1$ to be the vector fields
whose iterated commutators up to some order span the tangent space, and
use these to generate $Y$ as discussed above),
then the induced metric $\rho_Y \l( x,y\r)$ is equivalent to
$\l\| y^{-1} x\r\|$.  Indeed, it is easy to see that they are both
left invariant, and both homogeneous of order $1$ with respect
to the dilations on the group, and the equivalence then follows
from a simple compactness argument.
\end{rmk}

If $Y=\l\{ \l(Y_1,d_1\r),\cdots, \l( Y_q,d_q \r) \r\}$ is of finite
homogeneous type, and if $I=\l( i_1,\ldots, i_N \r)$ is an ordered
$N$-tuple of integers, with each $i_j\leq q$, we define:
\begin{equation*}
\lambda_I^Y \l( x \r) = {\rm det} \l( Y_{i_1}, \ldots, Y_{i_q} \r) \l(x\r)
\end{equation*}
where we regard each $Y_i$ as an $N$-tuple of smooth functions, and
$\lambda_I^Y$ is then the determinant of the corresponding $N\times N$
matrix.  We also set:
$$d\l( I \r) = d_{i_1} + \cdots + d_{i_N}$$
and define:
$$\Lambda_Y \l( x,\delta \r) = \sum_{I} \l| \lambda_I^Y \l( x\r)\r|\delta^{d\l(I\r)}$$

\begin{defn}
Let $S$ be a set of lists of vector fields $Y=\l\{ \l( Y_1,d_1 \r),\ldots, \l(Y_q,d_q\r) \r\}$ of homogeneous type, where $q$ may vary with $Y$.
We say $S$ is bounded if there is a uniform bound for $q$ for all
$Y\in S$, and the following hold:
\begin{enumerate}
\item There is an $M$ such that $d_j\leq M$ for all formal degrees
associated to some $Y\in S$.
\item We insist that the set of vector fields listed in some $Y\in S$,
thought of as sections of $T\Omega$, form a bounded set in the usual
topology of smooth sections of $T\Omega$.
\item The $c_{j,k}^l$ from Definition \ref{DefnHomogType} may be chosen
from a bounded subset of $C^\infty$ uniformly for $Y\in S$.
\item For every compact set $K\subset \Omega$, there exists a $c>0$ such
that for each $Y\in S$ we have
$\Lambda_Y \l(x,1\r)\geq c$
for $x\in K$.
\end{enumerate}
\end{defn}

The relevance of such bounded sets $S$
is that many of the results in 
\cite{NagelSteinWaingerBallsAndMetricsDefinedByVectorFields,NagelSteinDifferentiableControlMetricsAndScaledBumpFunctions}
hold uniformly for $Y\in S$, with no changes to the proof.  We shall need
some of these results and state them below.  We remark that these
bounded sets are the precompact sets in a natural topology on
the set of families of vector fields of homogeneous type; though we
do not expound on this further, as it will be of no use to us in the sequel.
Fix, for the remainder of this section, such a bounded set $S$. 
We will remind the reader of the results from
\cite{NagelSteinWaingerBallsAndMetricsDefinedByVectorFields,NagelSteinDifferentiableControlMetricsAndScaledBumpFunctions} that we will use,
and make explicit their uniformity in $S$.  All of the results in
this section follow by merely keeping track of the constants'
dependence on $Y$ in \cite{NagelSteinWaingerBallsAndMetricsDefinedByVectorFields,NagelSteinDifferentiableControlMetricsAndScaledBumpFunctions}. 

\begin{rmk}
The reader wishing to prove the results in this section may find it useful
to recall that the inverse function theorem remains true uniformly for
compact subsets of $C^{\infty}$.  Ie, if $R\subset C^{\infty}$ is a compact set, and
if $x$ is a point such that for all $f\in R$, the Jacobian determinant of $f$
at $x$ is non-zero (and hence, has absolute value bounded below, independent of $f$), then there exists an open neighborhood $U$ (independent
of $f$)
containing $x$ such that for all $f\in R$, $f:U\rightarrow f\l( U\r)$
is a diffeomorphism.  The essential point here is actually that $R$ is
a compact subset of $C^1$.
\end{rmk}

For a list of vector fields and formal degrees
$Y=\l\{ \l(Y_1,d_1\r),\ldots,\l( Y_q,d_q \r) \r\}$ of homogeneous type,
we define:
$$B_Y \l( x,\delta \r) = \l\{ y: \rho_Y\l( x,y \r)<\delta \r\}$$

\begin{defn}\label{DefnLocEquiv}
We say that two functions $\rho_1,\rho_2:\Omega\times \Omega\rightarrow \l[0,\infty\r]$
are locally equivalent if for every $x_0\in \Omega$ there exists an open
set $U$ containing $x_0$ such that for every compact set $K\subset\subset U$
there is a constant $C$ such that if $x_1,x_2\in K$,
$$\frac{1}{C} \rho_1\l( x_1,x_2 \r)\leq \rho_2\l( x_1,x_2\r) \leq C \rho_1\l( x_1,x_2\r)$$
\end{defn}

\cite{NagelSteinWaingerBallsAndMetricsDefinedByVectorFields} defines 
other pseudo-distances that are locally equivalent to $\rho_Y$, but can be
easier to work with.  We remind the reader of two of them that we shall use.
The definition of the first is similar to that of $\rho_Y$, but only allows constant
linear combinations of the vectors $Y_1,\ldots, Y_q$.  For $\delta>0$
let $C_2\l( \delta,Y\r)$ denote the class of smooth curves $\phi:\l[ 0,1\r]\rightarrow \Omega$
such that:
$$\phi'\l( t\r) = \sum_{j=1}^q a_j Y_j\l( \phi\l( t\r)\r)$$
with $\l| a_j \r|<\delta^{d_j}$.  Define:
$$\rho_{Y,2} \l(x,y\r) = \inf \l\{ \delta>0 : \bigg( \exists \phi\in C_2\l( \delta, Y \r)\bigg) \bigg( \phi\l(0\r)=x, \phi\l(1\r)=y \bigg) \r\}$$

\begin{thm}[Theorem 2 from \cite{NagelSteinWaingerBallsAndMetricsDefinedByVectorFields}]\label{ThmLocalEquiv}
$\rho_{2,Y}$ is locally equivalent to $\rho_Y$, with constants that
can be chosen uniformly for $Y\in S$.
\end{thm}

The definition of the second locally equivalent metric allows us
to single out $N$ of the vector fields $Y_{i_1},\ldots, Y_{i_N}$.  For
each $N$-tuple $I=\l( i_1,\ldots, i_N\r)$, let $C_3\l( \delta,I,Y\r)$
denote the class of smooth curves $\phi:\l[ 0,1\r] \rightarrow \Omega$
such that:
\begin{equation*}
\phi'\l( t\r) = \sum_{j=1}^N a_j Y_{i_j}\l(\phi\l(t\r)\r)
\end{equation*}
with $\l|a_j\r| < \delta^{d\l(Y_{i_j}\r)}$.  We define
\begin{equation*}
\rho_{Y,3}\l( x,y\r) = \inf\l\{ \delta>0: \bigg( \exists I\exists \phi\in C_3\l(\delta,I,Y\r)  \bigg)\bigg( \phi(0)=x, \phi(1)=y \bigg) \r\}
\end{equation*}

\begin{thm}[Theorem 3 from \cite{NagelSteinWaingerBallsAndMetricsDefinedByVectorFields}]\label{ThmCanUseMaximalDet}
$\rho_{3,Y}$ is locally equivalent to $\rho_Y$, with constants that can be chosen uniformly for $Y\in S$.  Moreover, for every $x_0\in \Omega$,
there exists an open set $U$ containing $x_0$ such that for every compact
set $K\subset\subset U$ we have the following for all $Y\in S$:
if for a fixed $x\in K$ we have:
\begin{equation*}
\delta^{d\l( I\r)} \l|\lambda_I^Y\l( x\r)\r|\geq \epsilon \Lambda_Y\l( x,\delta\r) 
\end{equation*}
then, there exists a $C$ depending on $\epsilon$ and $K$, but not on $Y$, such that for
every $y\in K$,
\begin{equation*}
\rho_{3,Y}\l( x,y\r) \geq C \inf \l\{ \delta>0 : \bigg( \exists\phi\in C_3\l( \delta,I,Y\r) \bigg)\bigg(\phi(0)=x, \phi(1)=y\bigg)  \r\}
\end{equation*}
\end{thm}

\begin{thm}[Theorem 1 from \cite{NagelSteinWaingerBallsAndMetricsDefinedByVectorFields}]\label{ThmEstBalls}
For every compact set $K\subset\subset \Omega$, there are constants $C_1, C_2$
such that for all $x\in K$ and all $Y\in S$,
$$0<C_1\leq \frac{\l|B_Y\l( x,\delta\r)\r|}{\Lambda_Y\l( x,\delta \r)}\leq C_2$$
where here, and in the rest of the paper, $\l| E\r|$ denotes
the Lebesgue measure of $E$.
\end{thm}

\begin{cor}[Corollary of Theorem \ref{ThmEstBalls}]\label{CorEstVolBall}
For every compact set $K\subset\subset \Omega$ there is a constant $C$
such that for all $Y\in S$ and all $x\in K$,
$$\l| B\l( x,2\delta \r) \r|\leq C \l|B\l( x,\delta \r)\r|$$
\end{cor}

\begin{thm}[Lemma 3.1.1 from \cite{NagelSteinDifferentiableControlMetricsAndScaledBumpFunctions}]\label{ThmBumpFuncs}
Let $E\subset\subset \Omega$ be compact.
There exist constants $\delta_0,\epsilon_0,\sigma_0>0$ such that
for each $x\in E$, each $0<\delta< \delta_0$, and each $Y\in S$, there is a function
$\phi=\phi_{x,\delta,Y}\in C^\infty\l(\Omega\r)$ such that:
\begin{enumerate}
\item For all $y\in \Omega$, $0\leq \phi\l(y\r)\leq 1$.
\item $\phi\l(y\r) =0$ when $\rho\l( x,y\r)>\sigma_0\delta$, $\phi\l( y\r)=1$
when $\rho\l( x,y\r)<\epsilon_0\delta$.
\item For every ordered multi-index $\alpha$,
$$\sup_{y\in \Omega}\l| \grad_Y^\alpha \phi\l(y\r) \r|\leq C_\alpha \delta^{-d\l( \alpha \r)}$$
where $\grad_Y=\l( Y_1,\ldots, Y_q\r)$, and $d\l( \alpha\r)$ is the formal degree
of $\grad_Y^\alpha$ with each $Y_j$ having formal degree $d_j$.
\end{enumerate}
\end{thm}
\begin{rmk}
Actually, the stronger results given by Theorem 3.3.1 and Theorem 3.3.2
of \cite{NagelSteinDifferentiableControlMetricsAndScaledBumpFunctions}
are true uniformly for $Y\in S$, but we will not need these
results.
\end{rmk}

	\subsection{A New Distance}\label{SectionNewDist}
		Given two metrics, $\rho_1$ and $\rho_2$, one may define a third function
$\rho_2\circ \rho_1:\Omega\times\Omega \rightarrow \l[ 0,\infty\r]$ defined
by:
$$\rho_2\circ \rho_1\l( x,y\r) = \inf\l\{ \delta>0 : \exists z\in \Omega, \rho_1\l(x,z\r)\leq \delta, \rho_2\l(z,y\r)\leq \delta \r\}$$
Suppose $X=\l\{ \l(X_1,c_1 \r), \ldots, \l( X_p,c_p \r)\r\}$
and $Y=\l\{ \l( Y_1,d_1\r),\ldots, \l( Y_q,d_q \r) \r\}$ are two
list of vector fields that are of finite homogeneous type.  Suppose
also that $\l[ X_j,Y_k \r]=0$ for every $1\leq j\leq p$, $1\leq k\leq q$.

We will see (under an assumption) that
$$\rho_Y \circ \rho_X = \rho_{X\cup Y}$$
In fact, even without our assumption, our proof works locally.  However,
the condition that the $X_j$ commute with the $Y_k$ is so
restrictive this is a moot point (see Remark \ref{RmkRemarkOnOurAssumptions}).

Before we speak about our assumption, a word of notation.  If
$\l( Z,d\r)$ appears in both $X$ and $Y$ we count it as appearing
twice in $X\cup Y$, equivalently, we replace $\l( Z,d\r)$ with
$\l( 2Z,d\r)$ in $X\cup Y$.

Our assumption is as follows: for every $\delta>0$ and every
$\l|a_j\l( t\r)\r|\leq \delta$, $1\leq j\leq p$, $a_j$ measurable,
and every $x\in \Omega$, there exists a unique solution in $\Omega$ to:
$$\phi\l( 0\r) = x$$
$$\phi'\l(t\r) = \sum_{j=1}^p a_j\l( t\r) X_j\l(\phi\l(t\r) \r)$$
and similarly for the $Y_k$.  We denote this solution by the
time-ordered exponential (also known as the product integral):
\begin{equation*}
\phi\l( t\r) = \Texp{ \int_0^t \sum_{j=1}^p a_j\l( s\r) X_j ds}x
\end{equation*}
See \cite{DollardFriedmanProductIntegrationWithApplicationsToDifferentialEquations,GamkrelidzeAgrachevOrdinaryDifferentialEquationsOnVectorBundles} for a background on product integration.

If $Z\l( s\r)$ is a family of vector fields (and is locally integrable),
then:
\begin{equation*}
\Texp{\int_0^t Z\l( s\r) ds}x = \lim_{N\rightarrow \infty} \exp\l( \int_{\frac{N-1}{N}t}^t Z\l(s\r) ds\r)\cdots \exp\l( \int_{0}^{\frac{t}{N}} Z\l( s\r) ds \r)x
\end{equation*}
From here we see that if $Z_1\l( s\r)$ and $Z_2\l( r\r)$ commute for
every $s$ and $r$, then:
\begin{equation}\label{EqnTimeCommute}
\Texp{\int_0^t Z_1\l(s\r) +Z_2\l( s\r) ds}x =\Texp{\int_0^t Z_1\l(s\r) ds}\Texp{\int_0^t Z_2\l( s\r) ds}x 
\end{equation}

\begin{thm}\label{ThmNewDist}
Under the setup above, we have:
$$\rho_Y \circ \rho_X = \rho_{X\cup Y}$$
\end{thm}
\begin{proof}
Suppose that $\rho_Y \circ \rho_X \l( x,y\r) <\delta$, so that there
exists a $z$ with $\rho_X\l( x, z\r) <\delta$ and $\rho_Y\l( z, y\r)<\delta$.
Let $\phi_X,\phi_Y:\l[ 0,1\r]\rightarrow \Omega$ be absolutely continuous
curves such that $\phi_X\l( 0\r) =x$, $\phi_X\l( 1\r) =z$, $\phi_Y\l( 0\r) =z$,
and $\phi_Y \l( 1\r) =y$, and such that
$$\phi_X\l( t\r) = \Texp{\int_0^t \sum_{j=1}^p a_j\l( s\r)X_j ds}x$$
$$\phi_Y\l( t\r) = \Texp{\int_0^t \sum_{k=1}^q b_k\l( s\r)Y_k ds}z$$
with $\l| a_j\r|< \delta^{c_j}$, $\l| b_j\r|< \delta^{d_j}$.
But then,
\begin{equation*}
\gamma\l( t\r) = \Texp{\int_0^t \sum_{j=1}^p a_j\l( s\r)X_j + \sum_{k=1}^q b_k\l( s\r) Y_k ds}x
\end{equation*}
is a path from $x$ to $z$.  Indeed,
\begin{equation*}
\begin{split}
\gamma\l( 1\r) &= \Texp{\int_0^1 \sum_{j=1}^p a_j\l( s\r)X_j + \sum_{k=1}^q b_k\l( s\r) Y_k ds}x\\
&=\Texp{\int_0^1 \sum_{k=1}^q b_k\l( s\r) Y_k ds}  \Texp{\int_0^1 \sum_{j=1}^p a_j\l( s\r)X_j ds}x\\
&= \Texp{\int_0^1 \sum_{k=1}^q b_k\l( s\r) Y_k ds}z\\
&= y
\end{split}
\end{equation*}
But, we also have that:
\begin{equation*}
\gamma'\l( t\r) = \sum_{j=1}^p a_j\l( t\r) X_j\l( \gamma\l(t\r)\r) + \sum_{k=1}^q b_k\l( t\r) Y_k\l( \gamma\l( t\r)\r)
\end{equation*}
and since $\l|a_j\r| < \delta^{c_j}$ and $\l| b_j \r|< \delta^{d_j}$
we see that $\rho_{X\cup Y}\l( x, y\r) < \delta$.

Conversely, suppose $\rho_{X\cup Y}\l( x,y\r) <\delta$.  Then there is a
path of the form:
\begin{equation*}
\phi\l( t\r) = \Texp{\int_0^t \sum_{j=1}^p a_j\l( s\r) X_j + \sum_{k=1}^q b_k\l( s \r) Y_k ds}x
\end{equation*}
with $\phi\l( 1\r) = y$, $\l|a_j\r|<\delta^{c_j}$, $\l|b_k\r|<\delta^{d_k}$.  Define
\begin{equation*}
\phi_X\l( t\r) = \Texp{\int_0^t \sum_{j=1}^p a_j\l( s\r) X_j ds}x
\end{equation*}
and let $z=\phi_X\l( 1\r)$.  Define:
\begin{equation*}
\phi_Y\l( t\r) = \Texp{\int_0^t \sum_{k=1}^q b_k\l( s \r) Y_k ds} z
\end{equation*}
Note that $\phi_Y\l( 1\r) = \phi\l( 1\r) = y$ by (\ref{EqnTimeCommute}).
It is easy to see that $\phi_X\in C\l( \delta, X\r)$, $\phi_Y\in C\l( \delta, Y\r)$, and it then follows that $\rho_X\l( x, z\r) < \delta$ and $\rho_Y\l( z, y\r)<\delta$, showing that $\rho_Y\circ\rho_X\l( x,y\r)<\delta$ and completing the proof.
\end{proof}

\begin{rmk}\label{RmkNoTimeOrdered}
If one wished to show only that $\rho_Y\circ\rho_X$ is locally equivalent
to $\rho_{X\cup Y}$ in Theorem \ref{ThmNewDist} (which would be sufficient
for our purposes), then the proof is a bit easier.  Indeed, Theorem
\ref{ThmLocalEquiv} would allow us to replace the exponentials with
variable coefficients with ones with constant coefficients.  Then
the same proof yields the result, without any need for time ordered
exponentials, nor the need for our assumption.  
In spite of this, we believe that the
proof of Theorem \ref{ThmNewDist} helps to elucidate the situation.
\end{rmk}

\begin{rmk}\label{RmkRemarkOnOurAssumptions}
One example of such $X_j$ and $Y_k$ is as follows:  let $X_j$ be a spanning
set of the right invariant vector fields on some Lie group, and let
$Y_k$ be a spanning set of the left invariant vector fields (and
we may even restrict them to a small connected open set).  It is
not hard to see that this is the only example.
\end{rmk}


\section{The Distances $\rl{\epsilon}$ and $\rr{\epsilon}$}
Given a finite set $F$ of vector fields such that $F$ along with 
the commutators of all orders of elements of $F$ up to some fixed
order $m$:
$$\l[ X_1,\l[ X_2,\ldots,\l[X_{n-1},X_{n}\r]\ldots\r]\r], \quad n\leq m, \quad X_j\in F$$
span the tangent space at each point (it is often said that such a
set $F$ satisfies H\"ormander's condition), we associate a list
of vector fields of finite homogeneous type as in Section \ref{SectionCCBackground},
by taking the list of all commutators up to order $m$ and associating
to a commutator of length $k$ degree $k$.  That is to say the elements
of $F$ are given degree $1$, elements of the form $\l[ X_1, X_2\r]$
are given degree $2$ and so forth.  
Call this list of vector fields $\sLp{F}$.
From this list of vector fields
of finite homogeneous type, we get a metric $\rho_{\sLp{F}}$.

We define, as in the introduction for $\epsilon\in \l[0,1\r]$,
$$\rl{\epsilon} = \rho_{\sLp{\l\{\gl, \epsilon \gr\r\}}}$$
$$\rr{\epsilon} = \rho_{\sLp{\l\{\gr, \epsilon \gl\r\}}}$$
Here, our set $\Omega$ from Section \ref{SectionCCBackground} is the
entire group $G$.  Now it is easy to see that 
$$\l\{ \sLp{\l\{\gl, \epsilon \gr\r\} }, \sLp{\l\{\gr, \epsilon \gl\r\}} | \epsilon\in [0,1]\r\}$$
is a bounded set as in Section \ref{SectionCCBackground}.
Thus, all of the theorems from that section hold uniformly for
$\epsilon\in \l[0,1\r]$.

We define $\Bl{\epsilon}{x}{\delta}$ to be the ball of radius $\delta$
centered at $x$ in the $\rl{\epsilon}$ metric, and we define
$\Vl{\epsilon}{x}{\delta}$ to be its volume.  Similarly,
we define $\Br{\epsilon}{x}{\delta}$ and $\Vr{\epsilon}{x}{\delta}$.
Note that all of the relevant quantities from Section \ref{SectionCCBackground}
are homogeneous of an appropriate degree.  For instance,
$$\rlp{\epsilon}{rx}{ry} = r \rlp{\epsilon}{x}{y}$$
$$\Vl{\epsilon}{rx}{\delta} = r^Q \Vl{\epsilon}{x}{\frac{\delta}{r}}$$
From such considerations, it is easy to see that all of the results
from Section \ref{SectionCCBackground} hold globally, instead of locally.
That is, many of the results are true on any fixed compact set $E$.
Take that compact set $E$ to be the closed unit ball.  Then to
see that the result holds globally, merely scale everything down
until it fits into the unit ball, and apply the result on the unit ball.
As all the quantities are homogeneous of the proper degrees, this extends
the results.  In the same manner, we may even take $\delta_0=\infty$ in
Theorem \ref{ThmBumpFuncs}.

\begin{rmk}\label{RmkScalingPropsOfMetrics}
We have the following scaling properties of the distances $\rl{\epsilon}$
and $\rr{\epsilon}$:
\begin{equation*}
r\rlp{\epsilon}{x}{y} = \rlp{\epsilon}{rx}{ry} = \rho_{\l\{ \frac{1}{r}\gl, \frac{1}{r}\epsilon\gr \r\}}\l(x,y\r)
\end{equation*}
and similarly for $\rr{\epsilon}$.  The first equality just follows
by homogeneity of the vector fields and was discussed above.  To see
that the first term equals the last term, note that:
$$C\l( \delta, \sLp{\gl, \epsilon \gr} \r) = C\l( r\delta, \sLp{\frac{1}{r} \gl, \frac{1}{r}\epsilon \gr}\r)$$
as can be seen directly from the definition.
\end{rmk}

	\subsection{Relationship to Convolution Operators}\label{SectionRelToConvo}
		In this section, we investigate the relationship between
$\rl{\epsilon}$, $\rr{\epsilon}$ and two-sided
convolution operators on $G$.  We will use one simplifying piece of
notation.  For an operator $T$, we write $\ker{T}\l( x,y \r)$ for
the Schwartz kernel of $T$ when mapping from the $y$ variable
to the $x$ variable.

Define $B:=\l\{ x : \l\| x\r\| < 1\r\}$.  With $\chi_B$ denoting
the characteristic function of $B$, set (for $r_L,r_R>0$):
\begin{equation}\label{EqnIntersectBalls}
\begin{split}
K_{r_L,r_R}\l( x,y\r) &= \ker{\Opl{\dil{\chi_B}{r_L}}\Opr{\dil{\chi_B}{r_R}}}\l( x,y\r)\\
&= r_L^Q r_R^Q \int \chi_{\frac{1}{r_L}B}\l( z^{-1} x \r) \chi_{\frac{1}{r_R}B}\l( z y^{-1}\r) dz\\
&= \l| \l\{ z : \l\| z^{-1}x \r\|<\frac{1}{r_L}, \l\| zy^{-1} \r\|<\frac{1}{r_R} \r\}  \r|
\end{split}
\end{equation}

Recall that $\l\| z^{-1}x\r\| \approx \rlp{0}{x}{z}$ (see Remark \ref{RmkLeftInvDist}) and, similarly,
$\l\| z y^{-1}\r\| \approx \rrp{0}{y}{z}$.  Using that
$r_L\rlp{0}{x}{y} = \rho_{\sLp{\frac{1}{r_L}\gl}}$ (and similarly for 
$r_R\rr{0}$) (see Remark \ref{RmkScalingPropsOfMetrics}), and using
(\ref{EqnIntersectBalls}), we see that
there exists a constant $C$ (independent of $r_L,r_R$) such
that:
\begin{equation}\label{EqnCompBalls}
\begin{split}
& \l\{ \l(x,y\r) : \rho_{\sLp{\frac{1}{r_R}\gr}}\circ \rho_{\sLp{\frac{1}{r_L}\gl}}\l( x,y \r)\leq C^{-1} \r\} \subseteq \l\{ \l(x,y\r) : K_{r_L,r_R}\l(x,y\r) \ne 0 \r\}\\
&\quad\quad\quad\quad\subseteq \l\{ \l( x,y\r) : \rho_{\sLp{\frac{1}{r_R}\gr}}\circ \rho_{\sLp{\frac{1}{r_L}\gl}}\l( x,y \r)\leq C  \r\}
\end{split}
\end{equation}
Section \ref{SectionNewDist} tells us that:
$$\rho_{\sLp{\frac{1}{r_R}\gr}}\circ \rho_{\sLp{\frac{1}{r_L}\gl}} = \rho_{\sLp{\frac{1}{r_R}\gr}\cup\sLp{\frac{1}{r_L}\gl}}=\rho_{\sLp{\frac{1}{r_R}\gr,\frac{1}{r_L}\gl}}$$
while Remark \ref{RmkScalingPropsOfMetrics} tells us that (when $r_L\leq r_R$):
$$\rho_{\sLp{\frac{1}{r_R}\gr,\frac{1}{r_L}\gl}} = r_L\rho_{\sLp{\frac{r_L}{r_R}\gr, \gl}} = r_L \rl{\frac{r_L}{r_R}}$$
with a similar result with $r_R\leq r_L$.  For the remainder
of this section, we restrict our attention to the case $r_L\leq r_R$,
with the understanding that the case $r_R\leq r_L$ follows in the
same way with completely symmetric arguments.

Putting all of this together, we see that:
\begin{equation}\label{EqnSuppOfConvo}
\Bl{\frac{r_L}{r_R}}{x}{\frac{1}{Cr_L}} \subseteq \l\{ y: K_{r_L,r_R}\l( x,y\r)\ne 0 \r\} \subseteq \Bl{\frac{r_L}{r_R}}{x}{\frac{C}{r_L}}
\end{equation}
Applying Corollary \ref{CorEstVolBall}, we see that:
\begin{equation*}
\l| \l\{ y: K_{r_L,r_R}\l( x,y\r)\ne 0 \r\}\r|  \approx \Vl{\frac{r_L}{r_R}}{x}{\frac{1}{r_L}} 
\end{equation*}
Define:
$$M\l( x, r_L,r_R \r) = \sup_{y} K_{r_L,r_R} \l( x,y\r)$$
The main result of this section is the following theorem
and its corollary:
\begin{thm}\label{ThmBoundIntersection}
$M\l( x, r_L, r_R\r) \approx \frac{1}{\Vl{\frac{r_L}{r_R}}{x}{\frac{1}{r_L}}}$.
Moreover, there is a $\delta_0>0$ (independent of $r_L,r_R$) such that
for all $x$, there exists a $y_0$ with:
$$K_{r_L,r_R}\l( x, y\r) \approx \frac{1}{\Vl{\frac{r_L}{r_R}}{x}{\frac{1}{r_L}}}$$
for all $y\in \Bl{\frac{r_L}{r_R}}{y_0}{\frac{\delta_0}{r_L}}$.
\end{thm}
\begin{cor}\label{CorBoundIntersection}
We may take $y_0=x$ in Theorem \ref{ThmBoundIntersection}.
\end{cor}

\begin{rmk}\label{RmkConvoBoundBump}
One of our main uses for Corollary \ref{CorBoundIntersection} is as follows.
Take $\delta_0$ as in the corollary, and set $\chi = \chi_{\frac{1}{\delta_0}B}$.
Then, it is easy to see that
$$\phi_0\l( y\r):=\ker{\Opl{\dil{\chi}{r_L}}\Opr{\dil{\chi}{r_R}}}\l( x,y\r) \gtrsim \frac{1}{\Vl{\frac{r_L}{r_R}}{x}{\frac{1}{r_L}}}$$
for $y\in \Bl{\frac{r_L}{r_R}}{x}{\frac{1}{r_L}}$.
Thus, when we wish bound a function $\phi$ supported on $\Bl{\frac{r_L}{r_R}}{x}{\frac{1}{r_L}}$ and which is $\lesssim \frac{1}{\Vl{\frac{r_L}{r_R}}{x}{\frac{1}{r_L}}}$, it suffices to instead bound $\phi_0$. 
\end{rmk}


\begin{lemma}\label{EstMAbove}
$M\l( x, r_L, r_R\r) \approx M\l( x, \frac{r_L}{2}, \frac{r_R}{2}\r)$
\end{lemma}
\begin{proof}
It is clear that
\begin{equation*}
M\l( x, r_L, r_R\r) \leq 2^{2Q}M\l( x, \frac{r_L}{2}, \frac{r_R}{2}\r)
\end{equation*}
and so we focus only on the reverse inequality.

For the proof of this lemma, alone, we drop the assumption that
$r_L\leq r_R$.  Then, it suffices to show that:
$$M\l( x, r_L, \frac{r_R}{2}\r) \lesssim M\l( x, r_L, r_R\r)$$
and the remainder of the result will follow by symmetry.

Define $g = \chi_B * \chi_{2B}$.  Note that there exists a
$c>0$ such that $g(x)>c$ for $x\in 2B$.  Hence, we have:
\begin{equation*}
\begin{split}
K_{r_L,\frac{r_R}{2}}&\l(x,y\r) = 2^{-Q}\ker{\Opl{\dil{\chi_B}{r_L}}\Opr{\dil{\chi_{2B}}{r_R}}}\l( x,y\r) \\
&\lesssim \ker{\Opl{\dil{\chi_B}{r_L}} \Opr{\dil{g}{r_R}}}\l( x,y\r)\\
&= \ker{\Opl{\dil{\chi_B}{r_L}} \Opr{\dil{\chi_B}{r_R} * \dil{\chi_{2B}}{r_R}}} \l( x,y\r)\\
&= \ker{\Opl{\dil{\chi_B}{r_L}} \Opr{\dil{\chi_B}{r_R}} \Opr{\dil{\chi_{2B}}{r_R}}}\l( x,y\r)\\
&= \int K_{r_L,r_R}\l( x, z\r) \dil{\chi_{2B}}{r_R}\l( y^{-1} z \r) dz\\
&\leq \int M\l( x, r_L, r_R\r) \dil{\chi_{2B}}{r_R}\l( y^{-1} z \r) dz \\
&= M\l( x, r_L, r_R\r) \int \chi_{2B}\l( z\r) dz\\
&\approx M\l( x, r_L, r_R\r)
\end{split}
\end{equation*}
Completing the proof.
\end{proof}

Let $\phi\in C^\infty\l( G\r)$ be $\geq 0$, $\leq 1$, supported on $2B$, and
$=1$ on $B$.
Define:
\begin{equation*}
\Kt_{r_L,r_R} = \ker{\Opl{\dil{\phi}{r_L}}\Opr{\dil{\phi}{r_R}}}
\end{equation*}
\begin{equation*}
\Mt\l(x, r_L, r_R\r) = \sup_y \Kt_{r_L,r_R}\l( x, y\r)
\end{equation*}
So that the definition of $K$ shows:
\begin{equation*}
K_{r_L,r_R} \l( x, y\r) \leq \Kt_{r_L,r_R} \l( x,y\r) \leq 2^{2Q} K_{\frac{r_L}{2},\frac{r_R}{2}}\l(x,y\r)
\end{equation*}
and Lemma \ref{EstMAbove} then tells us:
\begin{equation*}
M\l( x, r_L, r_R\r) \approx \Mt\l( x, r_L, r_R\r)
\end{equation*}

\begin{lemma}\label{LemmaChangeKt}
There exists $C_0>0$ (independent of $r_L,r_R$) such that for any $1\geq \delta>0$, and any $\gamma\in C\l( \frac{\delta}{r_L}, \sLp{\gl, \frac{r_L}{r_R}\gr}\r)$, we have:
\begin{equation*}
\l|\Kt_{r_L, r_R}\l( x, \gamma(1)\r)-\Kt_{r_L, r_R}\l( x, \gamma(0)\r)\r|\leq C_0\delta \Mt\l( x, r_L,r_R\r)
\end{equation*}
\end{lemma}

Before we prove Lemma \ref{LemmaChangeKt}, let's see how it finishes
the proof of Theorem \ref{ThmBoundIntersection}.  Fix $x$ and let
$y_0$ be such that:
\begin{equation*}
\Kt_{r_L,r_R}\l( x,y_0\r) = \Mt\l( x, r_L,r_R\r)
\end{equation*}
Then if we take $\delta_0=\delta=\min\l\{\frac{1}{2C_0},1\r\}$ in Lemma \ref{LemmaChangeKt},
we see that for $y\in \Bl{\frac{r_L}{r_R}}{y_0}{\delta_0}$,
\begin{equation*}
\Kt_{r_L,r_R}\l( x,y\r) \geq \frac{1}{2} \Mt\l( x, r_L,r_R\r)
\end{equation*}

Suppose, for contradiction that $\Mt\l( x, r_L, r_R\r) >> \frac{1}{\Vl{\frac{r_L}{r_R}}{x}{\frac{1}{r_L}}}$ (here $>>$ just means not $\lesssim$).  Then,
\begin{equation*}
\begin{split}
\int \Kt_{r_L,r_R}\l( x,y\r) dy &\geq \int_{y\in \Bl{\frac{r_L}{r_R}}{y_0}{\frac{\delta_0}{r_L}}} \Kt_{r_L,r_R}\l( x, y\r) dy\\
&\geq \frac{1}{2} \int_{y\in \Bl{\frac{r_L}{r_R}}{y_0}{\frac{\delta_0}{r_L}}} \Mt\l( x, r_L,r_R\r) dy\\
&=\frac{1}{2} \Mt\l( x, r_L,r_R\r) \Vl{\frac{r_L}{r_R}}{y_0}{\frac{\delta_0}{r_L}}\\
&>>\frac{1}{\Vl{\frac{r_L}{r_R}}{x}{\frac{1}{r_L}}} \Vl{\frac{r_L}{r_R}}{x}{\frac{1}{r_L}}\\
&=1
\end{split}
\end{equation*}
Here we used that
\begin{equation}\label{EqnUpperBall}
\Vl{\frac{r_L}{r_R}}{y_0}{\frac{\delta_0}{r_L}}\approx \Vl{\frac{r_L}{r_R}}{x}{\frac{1}{r_L}}
\end{equation}
as can be seen by the fact that since $y_0$ is in the support of $\Kt_{r_L,r_R}\l( x, y\r)$ we must have $\rlp{\frac{r_L}{r_R}}{x}{y_0}\lesssim \frac{1}{r_L}$, and thus,
$$\Bl{\frac{r_L}{r_R}}{x}{\frac{1}{r_L}}\subset \Bl{\frac{r_L}{r_R}}{y_0}{C_1\frac{\delta_0}{r_L}}$$
$$\Bl{\frac{r_L}{r_R}}{y_0}{\frac{\delta_0}{r_L}}\subset\Bl{\frac{r_L}{r_R}}{x}{\frac{C_1}{r_L}}$$
for some $C_1$.

But we also have,
\begin{equation*}
\int \Kt_{r_L,r_R}\l(x,y\r) dy = \int \phi\l( x\r) dx \int \phi\l( y \r) dy \approx 1
\end{equation*}
achieving the contradiction.

Hence we see that 
$\Mt\l( x, r_L, r_R\r)\lesssim \frac{1}{\Vl{\frac{r_L}{r_R}}{x}{\frac{1}{r_L}}}$.
But, $\Mt\l( x, r_L, r_R\r) \approx M\l( x, r_L, r_R\r)$, and applying
(\ref{EqnCompBalls}) we see:
\begin{equation*}
\begin{split}
M\l( x, r_L, r_R\r) &\geq \frac{1}{\Vl{\frac{r_L}{r_r}}{x}{\frac{C}{r_L}}}\int_{y\in \Bl{\frac{r_L}{r_R}}{x}{\frac{C}{r_L}}} K_{r_L, r_R}\l( x, y\r) dy\\
&= \frac{1}{\Vl{\frac{r_L}{r_R}}{x}{\frac{C}{r_L}}}\int K_{r_L, r_R}\l( x, y\r) dy\\
&=\frac{1}{\Vl{\frac{r_L}{r_R}}{x}{\frac{C}{r_L}}}\\
&\approx \frac{1}{\Vl{\frac{r_L}{r_R}}{x}{\frac{1}{r_L}}}
\end{split}
\end{equation*}
And so we see that:
\begin{equation*}
M\l( x, r_L, r_R\r) \approx \Mt\l( x, r_L, r_R\r) \approx \frac{1}{\Vl{\frac{r_L}{r_R}}{x}{\frac{1}{r_L}}}
\end{equation*}
proving the first part of Theorem \ref{ThmBoundIntersection}.
Moreover, we therefore have for all $y\in \Bl{\frac{r_L}{r_R}}{y_0}{\frac{\delta_0}{r_L}}$,
$$\Kt_{r_L,r_R}\l( x,y\r)\approx \frac{1}{\Vl{\frac{r_L}{r_R}}{x}{\frac{1}{r_L}}}$$
and since:
$$\Kt_{r_L,r_R}\l( x, y\r) \lesssim K_{\frac{r_L}{2},\frac{r_R}{2}}\l( x,y\r)$$
we have for all $y\in \Bl{\frac{r_L}{r_R}}{y_0}{\frac{\delta_0}{r_L}}$,
$$K_{\frac{r_L}{2},\frac{r_R}{2}}\l( x,y\r) \approx \frac{1}{\Vl{\frac{r_L}{r_R}}{x}{\frac{1}{r_L}}} \approx \frac{1}{\Vl{\frac{r_L}{r_R}}{x}{\frac{2}{r_L}}}$$

Dividing $\delta_0$ by $2$, this proves the second part of
Theorem \ref{ThmBoundIntersection} for 
$K_{\frac{r_L}{2},\frac{r_R}{2}}$.
Now merely replace $r_L$ and $r_R$ with $2r_L$ and $2r_R$ to complete
the proof.

\begin{proof}[Proof of Lemma \ref{LemmaChangeKt}]
We again use time-ordered exponentials, as in Section \ref{SectionNewDist};
and we again remark that their use is unnecessary, given Theorem \ref{ThmLocalEquiv} (see Remark \ref{RmkNoTimeOrdered}),
however we shall use them as we believe it adds to the clarity of the
exposition.

Let $\l\{\l( \Ylp{j}, d_j\r)\r\}=\sLp{\gl}$ and let
$\l\{\l( \Yrp{j}, d_j\r)\r\}=\sLp{\gr}$ so that
$$\l\{ \l(\Ylp{j}, d_j\r), \l( \epsilon^{d_j} \Yrp{j}, d_j\r) \r\}=\sLp{\gl, \epsilon \gr}$$
We remark that $\Ylp{j} f\l( rx\r) = r^{d_j} \l( \Ylp{j}f\r) \l( rx\r)$,
and similarly for $\Yrp{j}$.  We also remark that, if we do our enumeration
consistently between right and left, we have $\Yrp{j} f\l( x^{-1}\r) = \l(-1\r)^{d_j}\l( \Ylp{j}f\r) \l( x^{-1}\r)$.
Suppose $\gamma\in C\l( \frac{\delta}{r_L}, \sLp{\gl, \frac{r_L}{r_R} \gr}\r)$
with $\gamma\l( 0\r) = z_0$, $\gamma\l( 1\r) =z_1$;
so that
\begin{equation*}
\gamma\l( t\r) = \Texp{ \int_0^t \sum_{j} a_j\l( s\r) \Ylp{j} + \l( \frac{r_L}{r_R} \r)^{d_j} b_j\l( s\r) \Yrp{j}  ds} z_0
\end{equation*}
with $\l| b_j\r|,\l| a_j\r| <\l(\frac{\delta}{r_L}\r)^{d_j}$.

As in the proof of Theorem \ref{ThmNewDist}, we define:
$$\gamma_L\l( t\r) = \Texp{\int_0^t \sum_{j} a_j\l( s\r) \Ylp{j} ds}z_0$$
$$\gamma_R\l( t\r) = \Texp{\int_0^t \sum_{j}\l( \frac{r_L}{r_R} \r)^{d_j} b_j\l( s\r) \Yrp{j}  ds} \gamma_L\l( 1\r)$$
so that $\gamma_R\l( 1\r) = z_1$.
Consider,
\begin{equation*}
\begin{split}
&\frac{d}{dt} \Kt_{r_L,r_R}\l( x, \gamma_R\l( t\r) \r) = r_L^Q r_R^Q \frac{d}{dt} \int \phi\l( r_L\l( y^{-1}x \r) \r) \phi\l( r_R\l( y\gamma_R\l(t\r)^{-1}  \r) \r) dy\\
&\quad = r_L^Q r_R^Q \int \phi\l( r_L\l( y^{-1}x \r) \r) \l( \sum_j \l(-\frac{r_L}{r_R}\r)^{d_j} r_R^{d_j} b_j\l( t\r)\l( \Ylp{j} \phi \r)\l(r_R\l( y\gamma_R\l(t\r)^{-1}  \r)\r) \r) dy
\end{split}
\end{equation*}
Hence, using that $\l| b_j\r|<\l(\frac{\delta}{r_L}\r)^{d_j}$, and using that
$\phi$ is a fixed $C^\infty$ function supported on $2B$, we see that:
\begin{equation*}
\begin{split}
&\l|\frac{d}{dt} \Kt_{r_L,r_R}\l( x, \gamma_R\l( t\r) \r) \r| \lesssim r_L^Q r_R^Q \int \chi_{2B}\l( r_L\l( y^{-1} \r) \r) \delta\chi_{2B}\l( r_R\l( y\gamma_R\l(t\r)^{-1}  \r)\r)  dy\\
&\quad\quad\lesssim \delta M\l( x,\frac{r_L}{2},\frac{r_R}{2}\r)\\
&\quad\quad\approx \delta M\l( x, r_L,r_R\r)
\end{split}
\end{equation*}

Now, consider:
\begin{equation*}
\begin{split}
&\Kt_{r_L,r_R}\l( x, \gamma_L\l( t\r) \r) = \int \dil{\phi}{r_L}\l( y^{-1}x\r) \dil{\phi}{r_R}\l( y\gamma_L\l(t\r)^{-1}  \r) dy\\
&\quad\quad = \int \dil{\phi}{r_L}\l( \gamma_L\l(t\r)^{-1} y \r) \dil{\phi}{r_R}\l( xy^{-1} \r)dy
\end{split}
\end{equation*}
and then a similar proof to the one above shows that:
\begin{equation*}
\l|\frac{d}{dt} \Kt_{r_L,r_R}\l( x, \gamma_L\l( t\r) \r) \r| \lesssim \delta M\l( x, r_L,r_R\r)
\end{equation*}
Completing the proof.
\end{proof}

\begin{proof}[Proof of Corollary \ref{CorBoundIntersection}]
It is easy to see that
$$\Opl{\dil{\chi_B}{r_L}}\Opr{\dil{\chi_B}{r_R}}$$
is a self-adjoint operator, and thus,
$$K_{r_L,r_R}\l( x,y\r) = K_{r_L, r_R}\l( y,x\r)$$
If we take $y_0$ as in Theorem \ref{ThmBoundIntersection}, then
we see:
\begin{equation*}
\begin{split}
&\ker{\Opl{\dil{\l(\chi_B*\chi_B\r)}{r_L}} \Opr{\dil{\l(\chi_B*\chi_B\r)}{r_R}}}\l(x,x\r) \\&\quad\quad=  \ker{\l[\Opl{\dil{\chi_B}{r_L}}\Opr{\dil{\chi_B}{r_R}}\r]^2}\l(x,x\r)\\
&\quad\quad =\int K_{r_L, r_R}\l( x,y \r) K_{r_L,r_R}\l( x,y\r) dy\\
&\quad\quad \gtrsim \int_{y\in \Bl{\frac{r_L}{r_R}}{y_0}{\frac{\delta_0}{r_L}}} \l(\frac{1}{\Vl{\frac{r_L}{r_R}}{x}{\frac{1}{r_L}}}\r)^2 dy\\
&\quad \quad \approx \int_{y\in \Bl{\frac{r_L}{r_R}}{y_0}{\frac{\delta_0}{r_L}}} \frac{1}{\Vl{\frac{r_L}{r_R}}{y_0}{\frac{\delta_0}{r_L}}}\frac{1}{\Vl{\frac{r_L}{r_R}}{x}{\frac{1}{r_L}}} dy\\
&\quad\quad = \frac{1}{\Vl{\frac{r_L}{r_R}}{x}{\frac{1}{r_L}}} 
\end{split}
\end{equation*}
where we have applied (\ref{EqnUpperBall}) to get the second to last line.

Since $\chi_B*\chi_B\leq \chi_{4B}$, we have shown:
\begin{equation*}
\begin{split}
K_{\frac{r_L}{4}, \frac{r_R}{4}}\l( x, x\r) &\gtrsim \frac{1}{\Vl{\frac{r_L}{r_R}}{x}{\frac{1}{r_L}}}\\
&\approx \frac{1}{\Vl{\frac{r_L}{r_R}}{x}{\frac{4}{r_L}}}\\
&\approx M\l( x, \frac{r_L}{4}, \frac{r_R}{4}\r)
\end{split}
\end{equation*}
Replacing $r_L, r_R$ with $4r_L, 4r_R$, we see that:
\begin{equation*}
\Kt_{r_L,r_R}\l( x,x\r) \geq K_{r_L,r_R}\l( x, x\r) \approx M\l( x, r_L, r_R\r)\approx \Mt\l( x, r_L,r_R\r)
\end{equation*}
Lemma \ref{LemmaChangeKt} then tells us that there is a $\delta_0>0$ such
that for $y\in \Bl{\frac{r_L}{r_R}}{x}{\delta_0}$,
\begin{equation*}
\begin{split}
K_{\frac{r_L}{2},\frac{r_R}{2}}\l( x, y\r) &\geq \Kt_{r_L,r_R}\l( x,y\r)\\
&\gtrsim \Mt\l( x, r_L,r_R\r) \\
&\approx M\l( x, \frac{r_L}{2},\frac{r_R}{2}\r)
\end{split}
\end{equation*}
proving the result for $\frac{r_L}{2},\frac{r_R}{2}$ in place of $r_L, r_R$.
\end{proof}

We close this section with some simple inequalities that will be of use
in the sequel.
\begin{lemma}\label{LemmaChangeScale}
For $r_L,r_R, r_L^1,r_R^1 > 0$,
\begin{equation*}
K_{r_L^1,r_R^1}\l( x,z\r) \leq
\begin{cases}
\l(\frac{r_Lr_R^1}{r_L^1r_R}\r)^Q K_{r_L^1,\frac{r_L^1r_R}{r_L}}\l( x,z\r) & \text{ if $\frac{r_L}{r_R}\geq \frac{r_L^1}{r_R^1}$}\\
\l(\frac{r_Rr_L^1}{r_R^1r_L}\r)^Q K_{\frac{r_R^1r_L}{r_R},r_R^1}\l( x,z\r) & \text{ if $\frac{r_L}{r_R}\leq \frac{r_L^1}{r_R^1}$}
\end{cases}
\end{equation*}
\end{lemma}
\begin{proof}
This follows directly from the definition.
\end{proof}

\begin{cor}\label{CorChangeVolScale}
Suppose $r_L\leq r_R$ and $r_L^1\leq r_R^1$.  Then, we have:
\begin{equation*}
\frac{1}{\Vl{x}{\frac{r_L^1}{r_R^1}}{\delta}}\lesssim
\begin{cases}
\l(\frac{r_Lr_R^1}{r_L^1r_R}\r)^Q \frac{1}{\Vl{\frac{r_L}{r_R}}{x}{\delta}} & \text{ if $\frac{r_L}{r_R}\geq \frac{r_L^1}{r_R^1}$}\\
\l(\frac{r_L^1r_R}{r_Lr_R^1}\r)^Q \frac{1}{\Vl{\frac{r_L}{r_R}}{x}{\frac{r_Rr_L^1}{r_Lr_R^1}\delta}} & \text{ if $\frac{r_L}{r_R}\leq \frac{r_L^1}{r_R^1}$}
\end{cases}
\end{equation*}
In the case when $r_L\leq r_R$ but $r_L^1\geq r_R^1$, we have
\begin{equation}\label{EqnChangeVolScale2}
\frac{1}{\Vr{x}{\frac{r_R^1}{r_L^1}}{\delta}}\lesssim
\begin{cases}
\l( \frac{r_Lr_R^1}{r_L^1r_R} \r)^Q \frac{1}{\Vl{\frac{r_L}{r_R}}{x}{\frac{r_R^1}{r_L^1}\delta}} & \text{ if $\frac{r_L}{r_R}\geq\frac{r_L^1}{r_R^1}$}\\
\l( \frac{r_L^1r_R}{r_Lr_R^1} \r)^Q \frac{1}{\Vl{\frac{r_L}{r_R}}{x}{\frac{r_R}{r_L}\delta}} & \text{ if $\frac{r_L}{r_R}\leq\frac{r_L^1}{r_R^1}$}
\end{cases}
\end{equation}
\end{cor}
\begin{proof}
This follows by combining Corollary \ref{CorEstVolBall} with Lemma \ref{LemmaChangeScale}.
\end{proof}

\section{Bump Functions and Elementary Operators}\label{SectionBumpFuncsAndElem}
	We are now in a position to better understand the normalized bump
functions from Definition \ref{DefnBumpFunc} and the elementary
kernels from Definition \ref{DefnElemKernels}.  
The intuition for the normalized bump functions is easy to understand.
Indeed, if $\phi$ and $\psi$ are two $k$-normalized bump functions
in the sense of Definition \ref{KNormalizedBumpFunc}, then
$$\ker{\Opl{\dil{\phi}{r_L}}\Opr{\dil{\psi}{r_R}}}\l( x, \cdot\r)$$
is essentially an $r_L,r_R$ normalized bump function centered at $x$ of some
order, dependent on $k$.  This follows from (\ref{EqnSuppOfConvo}) and
Theorem \ref{ThmBoundIntersection}.  Following this analogy, we have:

\begin{lemma}\label{LemmaConvoOnBump}
Suppose $\phi_{r_L,r_R}^x$ is an $r_L,r_R$ normalized bump function centered at $x$,
and $\psi$ is a $C^\infty$ function supported in the unit ball $B$.  Then,
$$\Opl{\dil{\phi}{r_L'}}\phi_{r_L,r_R}^x$$
is a constant times a $r_L\wedge r_L',r_R$ normalized bump function, except
perhaps with support on a ball with a constant times the radius of the
support of a normalized bump function.  The order of $\Opl{\dil{\phi}{r_L'}}\phi_{r_L,r_R}^x$ will depend on the order of $\phi_{r_L,r_R}^x$ in a way
implicit in the proof.
\end{lemma}
\begin{proof}
The support and bounds of $\Opl{\dil{\phi}{r_L'}}\phi_{r_L,r_R}^x$ are easy to
see.  Indeed, fixing $\chi$ as in Remark \ref{RmkConvoBoundBump}, we see
that 
$$\l|\phi_{r_L,r_R}^x\l( z\r)\r| \lesssim \ker{\Opl{\dil{\chi}{r_L}}\Opr{\dil{\chi}{r_R}}}\l( z,x\r)$$
Thus,
\begin{equation*}
\begin{split}
&\l|\Opl{\dil{\phi}{r_L'}}\phi_{r_L,r_R}^x\l(z\r)\r|\lesssim\ker{\Opl{\dil{\phi}{r_L'}}\Opl{\dil{\chi}{r_L}}\Opr{\dil{\chi}{r_R}}}\l(x,z\r) 
\\&\quad\lesssim \ker{\Opl{\dil{\chit}{r_L\wedge r_L'}}\Opr{\dil{\chi}{r_L}}}\l( x,z\r)
\end{split}
\end{equation*}
where $\chit$ has some fixed bound and is supported in some fixed ball.
Thus the bounds for $\Opl{\dil{\phi}{r_L'}}\phi_{r_L,r_R}^x\l(z\r)$
follow from (\ref{EqnSuppOfConvo}) and
Theorem \ref{ThmBoundIntersection}.  It only remains to bound the derivatives
of $\Opl{\dil{\phi}{r_L'}}\phi_{r_L,r_R}^x\l(z\r)$.

For $\gr$ derivatives, this is easy.  Indeed,
$$\gr^{\alpha} \Opl{\dil{\phi}{r_L'}}\phi_{r_L,r_R}^x\l(z\r) = \Opl{\dil{\phi}{r_L'}} \gr^{\alpha}\phi_{r_L,r_R}^x\l(z\r)$$
and then the result follows from the definition of a normalized bump function
and our previous bounds.  Similarly, if $r_L'\leq r_L$, we have:
$$\gl^{\alpha} \Opl{\dil{\phi}{r_L'}}\phi_{r_L,r_R}^x\l(z\r) = r_L'^{\l|\alpha\r|}\Opl{\dil{\gl^{\alpha}\phi}{r_L'}} \phi_{r_L,r_R}^x\l(z\r)$$
The only problem that remains is when $r_L\leq r_L'$.  In that case we use
the following result:
\begin{equation}\label{EqnCommuteXl}
\gl^\alpha \Opl{\dil{\phi}{r_L'}} = \sum_{\l|\beta\r|=\l|\alpha\r|}\Opl{\dil{\phi_\beta}{r_L}}\gl^\beta
\end{equation}
where $\phi_\beta$ is of the same form as $\phi$.
From (\ref{EqnCommuteXl}), $\gl$ derivatives follow much in the same way
as $\gr$ derivatives.  We leave the details to the reader.
When one takes $\gl$ and $\gr$ derivatives simultaneously, the result
follows from a combination of the above two methods.

(\ref{EqnCommuteXl}) is well known, and so to save space and not introduce
too much extra notation, we prove it only in the case of the Heisenberg group (see Section \ref{SectionHeisenbergGroup} for the notation used here).
Indeed, suppose $\gl=\l( X_L,Y_L\r)$ and $\gr=\l( X_R,Y_R\r)$, with 
$\l[ X_L,Y_L\r] = 4\partial_t = -\l[ X_R,Y_R\r]$, and $X_L = X_R+4y\partial_t$.  Then,
\begin{equation*}
\begin{split}
&X_L \Opl{\dil{\phi}{r_L'}} = \l( X_R + 4yT\r) \Opl{\dil{\phi}{r_L'}}\\
&\quad= \Opl{\dil{\phi}{r_L'}} X_L + 4\frac{1}{r_L'}\Opl{\dil{y\phi}{r_L'}}\partial_t\\
&\quad= \Opl{\dil{\phi}{r_L'}} X_L + \frac{1}{r_L'}\Opl{\dil{y\phi}{r_L'}}\l[ X_L, Y_L\r]\\
&\quad= \Opl{\dil{\phi}{r_L'}} X_L + \Opl{\dil{X_Ly\phi}{r_L'}}Y_L - \Opl{\dil{Y_Ly\phi}{r_L'}}X_L
\end{split}
\end{equation*}
The proof of the more general result follows in a similar fashion.
\end{proof}

\begin{thm}\label{ThmGensisElemKer}
Suppose $\phi\in \SztSz$, then $\Opt{\dil{\phi}{r_L,r_R}}$ is an $r_L,r_R$
elementary operator.  Indeed, this is true uniformly as $\phi$
varies over a bounded set, and $r_L,r_R>0$ vary.
\end{thm}
\begin{proof}
By Remark \ref{RmkNuclear} it suffices to prove the result for
elementary tensor products.  Ie, we replace $\phi\l( x,y\r)$ with
$\phi\l(x\r) \psi\l(y\r)$, where $\phi,\psi\in \sSz$.  
Thus we are concerned with showing:
$$\Opl{\dil{\phi}{r_L}} \Opr{\dil{\psi}{r_R}}$$
is an $r_L,r_R$ elementary operator.

Since $\phi,\psi\in \sSz$, we may apply Lemma \ref{LemmaDecompForsSz},
for every $N_1,N_2,N_3,N_4\in \N$, we may write,
$$\phi = \sum_{\l|\alpha_1\r|=N_1}\sum_{\l|\alpha_2\r|=N_2} \gl^{\alpha_1}\gr^{\alpha_2} \phi_{\alpha_1,\alpha_2}$$
$$\psi = \sum_{\l|\beta_1\r|=N_3}\sum_{\l|\beta_2\r|=N_4} \gr^{\beta_1}\gl^{\beta_2} \psi_{\beta_1,\beta_2}$$
with $\phi_{\alpha_1,\alpha_2},\psi_{\beta_1,\beta_2}\in \sSz$.
Therefore,
\begin{equation*}
\begin{split}
&\Opl{\dil{\phi}{r_L}}\Opr{\dil{\phi}{r_R}} = r_L^{-N_1-N_2}r_R^{-N_3-N_4}
\\&\quad\times\sum_{\substack{\l|\alpha_1\r|=N_1, \l|\alpha_2\r|=N_2\\ \l|\beta_1\r| =N_3 , \l| \beta_2\r|=N_4}} \l(-1\r)^{\l|\alpha_2\r|+\l|\beta_2\r|} \gl^{\alpha_2}\gr^{\beta_1}\Opl{\dil{\phi_{\alpha_1,\alpha_2}}{r_L}} \Opl{\dil{\psi_{\beta_2,\beta_2}}{r_R}} \gl^{\alpha_2}\gr^{\beta_2}
\end{split}
\end{equation*}
and so to show that $\Opl{\dil{\phi}{r_L}}\Opr{\dil{\phi}{r_R}}$ is
an elementary operator, it suffices to show that:
$$\ker{\Opl{\dil{\phi}{r_L}}\Opr{\dil{\phi}{r_R}}}$$
satisfies the estimates of (\ref{EqnEleZero}).  For this purpose,
it will suffice to just assume $\phi,\psi\in \sS$.  Moreover, it is easy
to reduce the problem to the case when $\l| \alpha_1\r|=\l|\alpha_2\r|=\l|\beta_1\r|=\l|\beta_2\r|=0$, and so we prove it only in this case, leaving the
details to the reader.  Henceforth, we will only need that $\phi$ and $\psi$
are rapidly decreasing.

It is easy to see that it suffices to prove (\ref{EqnEleZero}) for
$r_L=2^{j_0}$, $r_R=2^{k_0}$ for $j_0,k_0\in \Z$ (a completely unnecessary reduction,
but it makes notation a little easier).  We also assume $k_0\geq j_0$, the
other situation being similar.
Let $\chi_0\l(x\r)=\chi_B\l(2x\r)$ ($B$ is $\l\|x\r\|<1$),
and $\chi_1\l(x\r) = \chi_{0}\l( x\r)-\chi_{0}\l( 2x\r)$, so that:
$$\chi_0\l(  x\r) + \sum_{j=1}^\infty \chi_1\l( 2^{-j} x  \r) = 1$$
Define $\phi_j$ by the equation:
\begin{equation*}
2^{\l(-j\r)Q} \phi_j\l( 2^{-j} x \r) = 
\begin{cases}
\chi_1\l( 2^{-j} x\r) \phi\l( x \r) & \text{if $j>0$,}\\
\chi_0\l( x\r)  \phi\l(  x\r) & \text{if $j=0$}
\end{cases} 
\end{equation*}
and define $\psi_k$ in a similar manner; so that:
\begin{equation*}
\sum_{j=0}^\infty \dil{\phi_j}{2^{-j}} =\phi
\end{equation*}
\begin{equation*}
\sum_{k=0}^\infty \dil{\psi_k}{2^{-k}} =\psi
\end{equation*}
and $\phi_j\l( x\r), \psi_k\l(x\r)$ are supported where $\l\|x \r\|\leq 1$.
Using that $\phi$ is rapidly decreasing, we see that for any $N\in \N$:
\begin{equation*}
\begin{split}
\l|\dil{\phi_j}{2^{j_0-j}}\l( x\r)\r| &\lesssim 2^{j_0 Q} \l( 1+2^{j} \r)^{-N}\chi_{\l\{\l\| x \r\|\leq 2^{j-j_0}\r\}}\l(x\r)\\
&\lesssim \l( 1+2^j\r)^{-N+Q} 2^{\l(j_0-j\r)Q}\chi_B\l(2^{j_0-j}x  \r)\\
&= \l( 1+2^j\r)^{-N'} \dil{\chi_B}{2^{j_0-j}}\l(x\r)
\end{split}
\end{equation*}
and similarly,
\begin{equation*}
\l| \dil{\psi_k}{2^{k_0-k}}\l(x\r) \r| \lesssim \l( 1+2^k\r)^{-N'}\dil{\chi_B}{2^{k_0-k}}
\end{equation*}
As in Section \ref{SectionRelToConvo}, we will use the notation:
$$K_{r_L,r_R}=\ker{\Opl{\dil{\chi_B}{r_L}}\Opr{\dil{\chi_B}{r_R}}}$$
We are ready to compute our main bound (we use (\ref{EqnSuppOfConvo}), Lemma \ref{LemmaChangeScale}, and Theorem \ref{ThmBoundIntersection} freely below):
\begin{equation*}
\begin{split}
&\ker{\Opl{\dil{\phi}{2^{j_0}}}\Opr{\dil{\psi}{2^{k_0}}}}\l( x,z\r) \\
&=
\sum_{\substack{j\geq 0\\ k\geq 0}} \ker{\Opl{\dil{\phi_j}{2^{j_0-k}}}\Opr{\dil{\psi_k}{2^{k_0-k}}}}\l( x,z\r)\\
&\lesssim \sum_{\substack{j\geq 0\\k\geq 0}} \l(1+2^j \r)^{-N_1} \l(1+2^{k}\r)^{-N_2} \ker{\Opl{\dil{\chi_B}{2^{j_0-j}}}\Opr{\dil{\chi_B}{2^{k_0-k}}}}\l( x,z\r)\\
&\lesssim \sum_{\substack{j\geq 0\\k\geq 0}} \l(1+2^j \r)^{-N_1} \l(1+2^{k}\r)^{-N_2}
\begin{cases}
2^{\l( j-k\r)Q} K_{2^{j_0-j},2^{k_0-j}}\l( x,z\r) & \text{ if $j\geq k$}\\
2^{\l( k-j\r)Q} K_{2^{j_0-k}, 2^{k_0-k}}\l(x,z\r) & \text{ if $k\geq j$}
\end{cases}\\
&\lesssim \sum_{\substack{j\geq 0\\k\geq 0}} \l(1+2^j \r)^{-N_1} \l(1+2^{k}\r)^{-N_2}
\begin{cases}
2^{\l( j-k\r)Q} \frac{\chi_{\l\{\rlp{2^{j_0-k_0}}{x}{z}\lesssim 2^{j-j_0}\r\}}}{\Vl{2^{j_0-k_0}}{x}{2^{j-j_0}}} & \text{ if $j\geq k$}\\
2^{\l( k-j\r)Q} \frac{\chi_{\l\{\rlp{2^{j_0-k_0}}{x}{z}\lesssim 2^{k-j_0}\r\}}}{\Vl{2^{j_0-k_0}}{x}{2^{k-j_0}}} & \text{ if $k\geq j$}\\
\end{cases}
\end{split}
\end{equation*}
Both the sum when $j\geq k$ and the sum when $k\geq j$ fall off faster than
a geometric series (for $N_1$ and $N_2$ chosen sufficiently large), and therefore are bounded by their first term.
For the sum when $j\geq k$, the first term is when $k$ is zero and
when $2^{j-j_0}\approx \rlp{2^{j_0-k_0}}{x}{z}$ (or when $j=0$ if such a
$j$ is less than $0$), with a similar result when $k\geq j$.  Hence, we see that:
\begin{equation*}
\begin{split}
&\ker{\Opl{\dil{\phi}{2^{j_0}}}\Opr{\dil{\psi}{2^{k_0}}}}\l( x,z\r) \\&\quad\quad\lesssim
\frac{1}{\l(1+2^{j_0}\rlp{2^{j_0-k_0}}{x}{z} \r)^N \Vl{2^{j_0-k_0}}{x}{2^{-j_0}+\rlp{2^{j_0-k_0}}{x}{z}}}
\end{split}
\end{equation*}
for any $N\geq 0$, completing the proof.
\end{proof}

Theorem \ref{ThmGensisElemKer} along with Proposition \ref{PropDecompOfTwosided}
show that every two-sided convolution operator can be decomposed as a
sum of elementary kernels.  In fact, this will be true for every
operator in $\sA'$; moreover, this will characterize $\sA'$.  We devote the rest of this section to proving these
facts.  Our first step is an analog of Lemma \ref{LemmaCZCancel}.

\begin{lemma}\label{LemmaCompElemKer}
Suppose $E_{2^{j_1},2^{k_1}}$ is a $2^{j_1},2^{k_1}$
elementary operator, and $E_{2^{j_2},2^{k_2}}$ is a
$2^{j_2},2^{k_2}$ elementary operator.  Then,
$$E_{2^{j_1},2^{k_1}}E_{2^{j_2},2^{k_2}}=2^{-\l|j_1-j_2 \r|-\l|k_1-k_2\r|}E_{2^{j_3},2^{k_3}}$$
where 
$j_3$ can be either $j_1$ or $j_2$ and $k_3$ can
be either $k_1$ or $k_2$ and
$E_{2^{j_3},2^{k_3}}$ is a $2^{j_3},2^{k_3}$ elementary
operator uniformly 
as $j_1,j_2,k_1,k_2$ vary over $\Z$ with constants only depending
on the constants for $E_{2^{j_1},2^{k_1}}$ and $E_{2^{j_2},2^{k_2}}$.
\end{lemma}
\begin{proof}
Let $\phi_{j_i,k_i}=\ker{E_{2^{j_i},2^{k_i}}}$, $i=1,2$.  Thus,
we are interested in the function
$$\phi_{j_3,k_3}\l( x,z\r) = \int \phi_{j_1,k_1}\l( x,y\r) \phi_{j_2,k_2}\l( y,z\r)dy$$
Suppose, for a moment, that $j_2\geq j_1$.  Then we see, from
Definition \ref{DefnElemKernels}, that
\begin{equation*}
\begin{split}
\int \phi_{j_3,k_3}\l( x,z\r) &=\sum_{\l|\alpha\r|=N} 2^{-Nj_2}\int \phi_{j_1,j_2}\l( x,y\r) \gly^{\alpha} \psi_{j_2,k_2,\alpha}\l( y,z\r) dy \\ 
&= \sum_{\l| \alpha\r| =N} 2^{N\l( j_1-j_2 \r)} \int \psi_{j_1,k_2,\alpha}\l( x,y\r) \phi_{j_2,k_2,\alpha}\l( y,z\r) dy
\end{split}
\end{equation*}
where the $\psi_{j_i,k_i,\alpha}$ are uniformly $2^{j_i},2^{k_1}$ elementary kernels.  Hence,
it suffices to consider only terms of the form:
\begin{equation*}
2^{N\l( j_1-j_2\r)} \int \psi_{j_1,k_1}\l( x,y\r) \psi_{j_2,k_2}\l( y,z\r) dy
\end{equation*}
where $N$ is any fixed large integer (which may depend on the semi-norm
we wish to estimate).  Doing the same argument for $k_1,k_2$, we see that it
suffices to consider only terms of the form:
\begin{equation*}
2^{-N\l| j_1-j_2\r|-N\l| k_1-k_2 \r|} \int \psi_{j_1,k_1}\l( x,y\r) \psi_{j_2,k_2}\l( y,z\r) dy
\end{equation*}

We proceed in the case when $k_1\geq j_1$ and $k_2\geq j_2$.  The three
other cases follow with only minor changes to the proof, and we 
leave those details to the interested reader.  
Define
\begin{equation*}
\chi_l\l( r\r)=
\begin{cases}
1 & \text{if $l>0$ and $2^{l-1}\leq r<2^l$,}\\
0 & \text{if $l>0$ and $r\geq 2^l$ or $r<2^{l-1}$,}\\
1 & \text{if $l=0$ and $r<1$,}\\
0 & \text{if $l=0$ and $r\geq 1$.}
\end{cases}
\end{equation*}
In the following, $N_1,N_2$ can be any two fixed large
integers we choose, and we let $\chi$ be as in Remark \ref{RmkConvoBoundBump}:
\begin{equation}\label{EqnCalcCompElem}
\begin{split}
&\l|\int \psi_{j_1,k_1}\l( x,y\r) \psi_{j_2,k_2}\l( y,z\r) dy\r| \\&\leq \sum_{\substack{l_1\geq 0\\l_2\geq 0}}
\l|\int \chi_{l_1}\l( 2^{j_1}\rlp{2^{j_1-k_1}}{x}{y} \r) \psi_{j_1,k_1}\l( x,y\r) \chi_{l_2}\l( \rlp{2^{j_2-k_2}}{y}{z}\r) \psi_{j_1,k_2}\l( y,z\r) dy\r|\\
&\lesssim \sum \int \frac{\chi_{\l\{\rlp{2^{j_1-k_1}}{x}{y}\leq 2^{l_1-j_1}\r\}}  \chi_{\l\{\rlp{2^{j_2-k_2}}{y}{z}\leq 2^{l_2-j_2}\r\}}}{ 2^{N_1l_1}  \Vl{2^{j_1-k_1}}{x}{2^{l_1-j_1}} 2^{N_2l_2}  \Vl{2^{j_2-k_2}}{x}{2^{l_2-j_2}} } dy \\
&\lesssim \sum 2^{-N_1l_1-N_2l_2} \int \ker{\Opl{\dil{\chi}{2^{j_1-l_1}}}\Opr{\dil{\chi}{2^{k_1-l_1}}}}\l( x,y\r)\\
&\quad\quad\quad\quad\times \ker{\Opl{\dil{\chi}{2^{j_2-l_2}}}\Opr{\dil{\chi}{2^{k_2-l_2}}}}\l( y,z\r) dy \\
&= \sum 2^{-N_1l_1-N_2l_2} \ker{\Opl{\dil{\chi}{2^{j_1-l_1}} * \dil{\chi}{2^{j_2-l_2}} }    \Opr{\dil{\chi}{2^{k_2-l_2}} * \dil{\chi}{2^{k_1-l_1}} }} \l( x,z\r)\\
&= \sum 2^{-N_1l_1-N_2l_2} \ker{\Opl{\dil{\chit}{2^{\l( j_1-l_1 \r)\wedge \l( j_2-l_2 \r)}}}\Opr{\dil{\chitt}{2^{\l( k_1-l_1 \r)\wedge \l( k_2-l_2 \r)}}}}\l(x,z\r)
\end{split}
\end{equation}
where $\chit$ and $\chitt$ are non-negative bounded functions with support in a fixed
bounded set, with these bounds independent of $j_1,j_2,k_1,k_2,l_1,l_2$.
Note that we could have achieved the same left hand side for (\ref{EqnCalcCompElem})
in the three other cases where we allow $k_1<j_1$ or $k_2<j_2$ or both. 
We will now drop our assumption $k_1\geq j_1$, $k_2\geq j_2$, though we
will return to it at the end.

Let $\Bt$ be a large fixed ball containing the support of $\chit$ and $\chitt$,
and define:
$$\Kt_{r_1, r_2}\l( x, z\r) = \ker{\Opl{\dil{\chi_{\Bt}}{r_1}} \Opr{\dil{\chi_{\Bt}}{r_2}}  } \l( x,z\r) $$
so that we have:
\begin{equation}\label{EqnCompElemBound}
\begin{split}
&2^{-N\l|j_1-j_2\r| -N\l|k_1-k_2\r|}\l|\int \psi_{j_1,k_1}\l( x,y\r) \psi_{j_2,k_2}\l( y,z\r) dy\r| \\&\lesssim 2^{-N\l|j_1-j_2\r| -N\l|k_1-k_2\r|}\sum_{\substack{l_1\geq 0 \\ l_2\geq 0}} 2^{-N_1l_1-N_2l_2}\Kt_{2^{\l( j_1-l_1 \r)\wedge \l( j_2-l_2 \r)}, 2^{\l( k_1-l_1 \r)\wedge \l( k_2-l_2 \r)}}\l( x, z\r)
\end{split}
\end{equation}

We now proceed in proving the lemma in the case when $j_3=j_1$ and
$k_3=k_1$.  The case when $j_3=j_2$ and $k_3=k_2$ is completely
symmetric.  The remaining two cases follow by similar arguments, 
and we leave those proofs to the reader.  We also work in the case when $k_1\geq j_1$,
the other case being symmetric.

We separate the RHS of (\ref{EqnCompElemBound}) into 4 sums:  depending
on whether $j_1-l_1\leq j_2-l_2$ and whether $k_1-l_1\leq k_2-l_2$.  The
first case we deal with is the sum over those $l_1$ and $l_2$ such that
$j_1-l_1\leq j_2-l_2$ and $k_1-l_1\leq k_2-l_2$.  In this case, we
need only take $N=1$.  In what follows, we will use (\ref{EqnSuppOfConvo}), Lemma \ref{LemmaChangeScale}, and Theorem \ref{ThmBoundIntersection}
freely (indeed, we will use their analogs for $\Kt$ which follow from
the methods in Section \ref{SectionRelToConvo}; note the $\Kt$ we use here differs slightly
from the one in Section \ref{SectionRelToConvo}).
We have:
\begin{equation*}
\begin{split}
&2^{-\l|j_1-j_2\r| -\l|k_1-k_2\r|}\sum_{l_1,l_2} 2^{-N_1l_1-N_2l_2}\Kt_{2^{\l( j_1-l_1 \r)
}, 2^{\l( k_1-l_1 \r)}}\l( x, z\r)\\
&\lesssim 2^{-\l|j_1-j_2\r| -\l|k_1-k_2\r|} \sum_{l_1\geq 0} 2^{-N_1l_1} \Kt_{2^{\l( j_1-l_1 \r)},2^{\l( k_1-l_1 \r)}}\l( x, z\r)\\
&\lesssim 2^{-\l|j_1-j_2\r| -\l|k_1-k_2\r|} \sum_{l_1\geq 0}2^{-N_1l_1} \frac{\chi_{\l\{ \rlp{2^{j_1-k_1}}{x}{z}\lesssim 2^{l_1-j_1} \r\}}}{\Vl{2^{j_1-k_1}}{x}{2^{l_1-j_1}}}
\end{split}
\end{equation*}
This sum falls off geometrically, and is therefore bounded
by a multiple of its first term, which occurs when
$\rlp{2^{j_1-k_1}}{x}{z}\approx 2^{l_1-j_1}$, or when $l_1=0$ (whichever
$l_1$ is greater).
Thus, we have that this sum is:
\begin{equation*}
\lesssim 2^{-\l| j_1-j_2\r| - \l| k_1-k_2\r|} \frac{1}{\l( 1+ 2^{j_1}\rlp{2^{j_1-k_1}}{x}{z} \r)^{N_1} \Vl{2^{k_1-j_1}}{x}{2^{-j_1}+\rlp{2^{j_1-k_1}}{x}{z}}}
\end{equation*}
which is $2^{-\l| j_1-j_2\r| - \l| k_1-k_2\r|}$ times the bound for a $2^{j_1},2^{k_1}$
elementary kernel.

We now turn to the case when $j_2-l_2\leq j_1-l_1$ and $k_2-l_2\leq k_1-l_1$.
As we estimate this case, we will use the fact that we may choose $N$ to
be large and this will allow us to absorb some terms by changing $N$.
When we do this, we will replace $N$ by $N'$ and then by $N''$, etc.
\begin{equation*}
\begin{split}
&\sum_{l_1,l_2} 2^{-N\l|j_1-j_2\r|-N\l|k_1-k_2\r|-N_1l_1-N_2l_2} \Kt_{2^{j_2-l_2},2^{k_2-l_2}}\l( x,z\r)\\
&\lesssim \sum_{l_2} 2^{-N\l|j_1-j_2\r|-N\l|k_1-k_2\r|-N_2l_2} \Kt_{2^{j_2-l_2},2^{k_2-l_2}}\l( x,z\r) \\
&\lesssim \sum_{l_2} 2^{-N\l|j_1-j_2\r|-N\l|k_1-k_2\r|-N_2l_2}
\\&\quad\quad\times \begin{cases}
2^{\l(j_1-k_1+k_2-j_2\r)Q} \Kt_{2^{j_2-l_2},2^{j_2-l_2+k_1-j_1}}\l(x,z\r) & \text{ if $j_1-k_1\geq j_2-k_2$}\\
2^{\l(j_2-k_2+k_1-j_1\r)Q} \Kt_{2^{k_2-l_2+j_1-k_1}, 2^{k_2-l_2}}\l( x,z\r) & \text{ if $j_1-k_1\leq j_2-k_2$}
\end{cases}\\
&\lesssim \sum_{l_2} 2^{-N'\l|j_1-j_2\r|-N'\l|k_1-k_2\r|-N_2l_2}\begin{cases}
\Kt_{2^{j_2-l_2},2^{j_2-l_2+k_1-j_1}}\l(x,z\r) \\
\Kt_{2^{k_2-l_2+j_1-k_1}, 2^{k_2-l_2}}\l( x,z\r) 
\end{cases}\\
&\lesssim \sum_{l_2} 2^{-N'\l|j_1-j_2\r|-N'\l|k_1-k_2\r|-N_2l_2}\begin{cases}
\frac{\chi_{\l\{ \rlp{2^{j_1-k_1}}{x}{z}\lesssim 2^{l_2-j_2}   \r\}}}{\Vl{2^{j_1-k_1}}{x}{2^{l_2-j_2}}}\\
\frac{\chi_{\l\{ \rlp{2^{j_1-k_1}}{x}{z}\lesssim 2^{k_1-j_1+l_2-k_2}   \r\}}}{\Vl{2^{j_1-k_1}}{x}{ 2^{k_1-j_1+l_2-k_2}  }}
\end{cases}
\end{split}
\end{equation*}
This sum falls off geometrically, and is therefore bounded by its first term.
Thus, we have:
\begin{equation*}
\begin{split}
&\sum_{l_1,l_2} 2^{-N\l|j_1-j_2\r|-N\l|k_1-k_2\r|-N_1l_1-N_2l_2} \Kt_{2^{j_2-l_2},2^{k_2-l_2}}\l( x,z\r)\\
&\lesssim  2^{-N'\l|j_1-j_2\r|-N'\l|k_1-k_2\r|} 
\\&\quad\quad\times\begin{cases}
\frac{1}{\l( 1+2^{j_2}\rlp{2^{j_1-k_1}}{x}{z} \r)^{N_2}  \Vl{2^{j_1-k_1}}{x}{2^{-j_2}+\rlp{2^{j_1-k_1}}{x}{z} }  }\\
\frac{1}{\l( 1+ 2^{k_2+j_1-k_1}\rlp{2^{j_1-k_1}}{x}{z} \r)^{N_2} \Vl{2^{j_1-k_1}}{x}{2^{k_1-j_1-k_2}+\rlp{2^{j_1-k_1}}{x}{z}}}
\end{cases}\\
&\lesssim \frac{2^{-N''\l|j_1-j_2\r|-N''\l|k_1-k_2\r|}}{\l( 1+2^{j_1}\rlp{2^{j_1-k_1}}{x}{z} \r)^{N_2}\Vl{2^{j_1-k_1}}{x}{2^{-j_1}+\rlp{2^{j_1-k_1}}{x}{z}}  }
\end{split}
\end{equation*}
thereby completing the bound in this case.

We now turn to the case when $j_1-l_1\leq j_2-l_2$ and $k_1-l_1\geq k_2-l_2$.
\begin{equation*}
\begin{split}
&\sum_{l_1,l_2} 2^{-N\l|j_1-j_2\r|-N\l|k_1-k_2\r|-N_1l_1-N_2l_2} \Kt_{2^{j_1-l_1 },2^{k_2-l_2}}\l( x,z\r)\\
&\lesssim \sum_{l_1,l_2} 2^{-N\l|j_1-j_2\r|-N\l|k_1-k_2\r|-N_1l_1-N_2l_2} 2^{\l( k_1-k_2+l_2-l_1\r)Q}\Kt_{2^{j_1-k_1+k_2-l_2}, 2^{k_2-l_2}}\l( x,z\r)
\end{split}
\end{equation*}
where, when we applied Lemma \ref{LemmaChangeScale}, only the latter case
($k_1-l_1\geq k_2-l_2$) applies.  Continuing our bound, we have:
\begin{equation*}
\begin{split}
&\lesssim \sum_{l_2} 2^{-N'\l|j_1-j_2\r|-N'\l|k_1-k_2\r|-N_2'l_2} \Kt_{2^{j_1-k_1+k_2-l_2}, 2^{k_2-l_2}}\l( x,z\r)
\end{split}
\end{equation*}
but this is just the lower case for our computation when $k_1-l_1\geq k_2-l_2$
and $j_1-l_1\geq j_2-l_2$.  Thus, we have:
\begin{equation*}
\begin{split}
&\sum_{l_1,l_2} 2^{-N\l|j_1-j_2\r|-N\l|k_1-k_2\r|-N_1l_1-N_2l_2} \Kt_{2^{j_1-l_1 },2^{k_2-l_2}}\l( x,z\r)\\
&\lesssim \frac{2^{-N''\l|j_1-j_2\r|-N''\l|k_1-k_2\r|}}{\l(
 1+2^{j_1}\rlp{2^{j_1-k_1}}{x}{z} \r)^{N_2'}\Vl{2^{j_1-k_1}}{x}{2^{-j_1}+\rlp{2^{
j_1-k_1}}{x}{z}}  }
\end{split}
\end{equation*}

Finally, when $j_1-l_1\geq j_2-l_2$ and $k_1-l_1\leq k_2-l_2$, the proof
proceeds as in the previous case, but now one ends up with the upper case
for our computation when $k_1-l_1\geq k_2-l_2$
and $j_1-l_1\geq j_2-l_2$.  Putting all of this together, we have:
\begin{equation*}
\begin{split}
&\l| \int \phi_{j_1,k_1}\l( x,y\r) \phi_{j_2,k_2}\l( y,z\r)dy\r| 
\\&\lesssim 
\frac{2^{-\l| j_1-j_2 \r|-\l| k_1-k_2 \r|}}{\l( 1+2^{j_1}\rlp{2^{j_1-k_1}}{x}{z} \r)^N  \Vl{2^{j_1-k_1}}{x}{2^{-j_1}+\rlp{2^{j_1-k_1}}{x}{z}} }
\end{split}
\end{equation*}
for any $N$ we choose.

Now let's turn to derivatives.  Fix ordered multi-indicies $\alpha_1,\beta_1,\alpha_2,\beta_2$, and consider:
\begin{equation*}
\begin{split}
&\l| \glx^{\alpha_1}\glz^{\alpha_2}\grx^{\beta_1}\grz^{\beta_2} \int \phi_{j_1,k_1}\l( x,y\r) \phi_{j_2,k_2}\l( y,z\r) dy \r|\\
&\leq\sum \l| \glx^{\alpha_1}\glz^{\alpha_2}\grx^{\beta_1}\grz^{\beta_2} 2^{-\l|\alpha_2\r|\l| k_1-k_2 \r|-\l| \beta_2 \r|\l|j_1-j_2\r|}\int \psi_{j_1,k_1}\l( x,y\r) \psi_{j_2,k_2}\l( y,z\r) dy \r|
\end{split}
\end{equation*}
where this is some finite sum and the $\psi_{j_i,k_i}$ are elementary
kernels, as we saw at the start of the proof.  Applying the definition of elementary kernels, we see:
\begin{equation*}
\begin{split}
&\l| \glx^{\alpha_1}\glz^{\alpha_2}\grx^{\beta_1}\grz^{\beta_2} 2^{-\l|\alpha_2\r|\l| k_1-k_2 \r|-\l| \beta_2 \r|\l|j_1-j_2\r|} \int \psi_{j_1,k_1}\l( x,y\r) \psi_{j_2,k_2}\l( y,z\r) dy \r|\\
&= 2^{j_1\l(\l| \alpha_1\r|+\l|\alpha_2\r|\r)+k_1\l( \l|\beta_1 \r| +\l|\beta_2 \r|\r)} \l| \int \psit_{j_1,k_1}\l( x,y\r) \psit_{j_2,k_2}\l( y,z\r) dy \r|
\end{split}
\end{equation*}
where the $\psit$ are also elementary kernels.  Hence, our bound for the
composition of two elementary kernels proved above, applied to:
\begin{equation*}
\l| \int \psit_{j_1,k_1}\l( x,y\r) \psit_{j_2,k_2}\l( y,z\r) dy\r|
\end{equation*}
gives the proper bound from the definition of $2^{j_1},2^{k_1}$ elementary
kernels for:
\begin{equation*}
\l| \glx^{\alpha_1}\glz^{\alpha_2}\grx^{\beta_1}\grz^{\beta_2} \int \phi_{j_1,k
_1}\l( x,y\r) \phi_{j_2,k_2}\l( y,z\r) dy \r|
\end{equation*}

Finally, we need to see that we may ``pull out'' derivatives as in 
Definition \ref{DefnElemKernels}.
Pulling out $x$ derivatives, works easily:
\begin{equation*}
\begin{split}
&\int \phi_{j_1,k_1}\l( x,y\r) \phi_{j_2,k_2}\l( y,z\r) dy\\
&= \sum_{\substack{\l|\alpha\r|=N_1\\\l| \beta \r|=N_3}} 2^{-j_1N_1-k_1N_3}\glx^{\alpha}\grx^{\beta}\int
\psi_{j_1,k_1,\alpha,\beta}\l( x,y\r) \phi_{j_2,k_2}\l( y,z\r) dy
\end{split}
\end{equation*}
and is therefore $2^{-j_1N_1-k_1N_3}$ times a sum of terms of the same form.

Pulling out $z$ derivatives takes one more step.  Indeed, fix $N_2$ and $N_4$
and suppose we wish to pull out $z$ left derivatives of order $N_2$ and
right derivatives of order $N_4$ (as in Definition \ref{DefnElemKernels}).
Then consider,
\begin{equation*}
\begin{split}
&\int \phi_{j_1,k_1}\l( x,y\r) \phi_{j_2,k_2}\l( y,z\r) dy\\
&=\sum 2^{-N_2\l| j_1-j_2 \r|-N_4\l| k_1-k_2\r|}\int \psi_{j_1,k_1}\l( x,y\r) \psi_{j_2,k_2}\l( y,z\r) dy\\
&= \sum \sum_{\substack{\l|\alpha\r|=N_2\\\l| \beta \r|=N_4}} \glz^{\alpha} \grz^{\beta} 2^{-N_2j_1-N_4k_1}\int \psi_{j_1,k_1}\l( x,y\r) \psi_{j_2,k_2,\alpha,\beta}\l( y,z\r) dy
\end{split}
\end{equation*}
where the sum above is, as usual, a finite sum over such terms, as we
saw in the beginning of the proof.  This completes the proof.
\end{proof}

\begin{lemma}\label{LemmaElemKerAreSchwarz}
Suppose $\phi$ is an $r_L,r_R$ elementary kernel, $r_L\leq r_R$.  Then if we define
$\dil{\psi}{r_L}\l( x\r) =\phi\l( x,0\r)$, we have that $\psi\in \sSz$,
uniformly for $\phi$ which are uniformly $r_L,r_R$ elementary kernels,
with constants independent of $r_L,r_R$.  When $r_R\leq r_L$, we have:
$\dil{\psi}{r_R}\l( x\r) =\phi\l( x,0\r)$
yields $\psi$ uniformly in $\sSz$.
\end{lemma}
\begin{proof}
This is a simple consequence of the definitions.
\end{proof}

\begin{cor}\label{CorElemKerOnsSz}
Suppose $\phi\in \sSz$, $l\in Z$, and $E_{2^{j},2^{k}}$ is a 
$2^{j},2^{k}$ elementary kernel.  Then, 
\begin{equation}\label{EqnToProveOnsSz}
E_{2^{j},2^{k}} \dil{\phi}{2^{l}} = 2^{-\l| l-j \r|-\l|l-k\r|}\dil{\psi}{2^{l}}
\end{equation}
where $\psi\in \sSz$.  As $\phi$ and $E_{2^{j},2^{k}}$ range 
over bounded sets, so
does $\psi$.  Moreover, this is true uniformly in $j,k,l$.
\end{cor}
\begin{proof}
Define $\dil{\psit}{2^l}=E_{2^{j},2^{k}} \dil{\phi}{2^{l}}$.  First, let us see
that it will suffice to show that $\psit$ is rapidly decreasing (uniformly,
in the relevant parameters).
Indeed, suppose we have that it is rapidly decreasing.  First, let
us see how to obtain the
factor
$2^{-\l| l-j \r|-\l|l-k\r|}$ in (\ref{EqnToProveOnsSz}).
Consider, in the case $l\leq j$,
\begin{equation*}
\begin{split}
E_{2^{j},2^{k}} \dil{\phi}{2^{l}} &= \sum_{\l|\alpha\r|=1}2^{-j}\Et_{2^{j},2^{k},\alpha} \gl^{\alpha} \dil{\phi}{2^{l}}\\
&= 2^{l-j} \sum_{\l|\alpha\r|=1} \Et_{2^{j},2^{k},\alpha} \dil{\phit_\alpha}{2^{l}}
\end{split}
\end{equation*}
and so is a finite sum of terms of the same form, but now with a factor of
$2^{l-j}$ out front.

On the other hand, if $j\leq l$, we may apply Lemma \ref{LemmaDecompForsSz}
to see:
\begin{equation*}
\begin{split}
E_{2^{j},2^{k}} \dil{\phi}{2^{l}} &= \sum_{\l| \alpha\r| =1} 2^{-l} E_{2^{j},2^{k}}\gl^{\alpha} \dil{\phit_{\alpha}}{2^{l}}\\
&= \sum_{\l| \alpha\r|=1} 2^{j-l} \Et_{2^{j},2^{k},\alpha} \dil{\phit_\alpha}{2^{l}}
\end{split}
\end{equation*}
and so it is a finite sum of terms of the same form, but now with a factor
of $2^{j-l}$ out front.  In a similar manner we may obtain a factor
of $2^{-\l| k-l\r|}$ out front.

Thus we have seen that, given that $\psit$ is rapidly decreasing, we
have that $\psi$ in the statement of the corollary, is rapidly decreasing,
uniformly in the relevant parameters.  Let us turn to derivatives 
of $\psi$.  It is easy to see from the Definition \ref{DefnElemKernels}
that if $\l|\alpha\r|=1$,
\begin{equation}\label{EqnElemCommutes}
\gl^\alpha E_{2^j,2^k} = \sum_{\l|\beta\r|=1} \Et_{2^j,2^k} \gl^\beta
\end{equation}
Hence,
\begin{equation*}
\begin{split}
\gl^\alpha E_{2^j,2^k} \dil{\phi}{2^{l}} &= \sum_{\l|\beta\r|=1} \Et_{2^j,2^k} \gl^\beta \dil{\phi}{2^{l}}\\
&= 2^l \sum_{\l|\beta\r|=1} \Et_{2^j,2^k} \dil{\phit_\beta}{2^l}
\end{split}
\end{equation*}
and so is $2^l$ times a finite sum of terms of the same form.  Thus,
$\psi$ behaves properly under derivatives, and we have shown that $\psi$
is uniformly in $\sS$.

To see $\psi$ is uniformly in $\sSz$, we need to ``pull out'' derivatives.
However, we merely use the other direction of (\ref{EqnElemCommutes}) to
see:
\begin{equation*}
\begin{split}
E_{2^j,2^k} \dil{\phi}{2^{l}} &= 2^{-l} \sum_{\l|\alpha\r|=1} E_{2^j,2^k} \gl^\alpha \dil{\phit_\alpha}{2^{l}}
\\&= \sum_{\l|\alpha\r|=1}\sum_{\l|\beta\r|=1} 2^{-l} \gl^{\beta} \Et_{2^j,2^k,\beta} \dil{\phit_\alpha}{2^{l}}
\end{split}
\end{equation*}
and so one can ``pull out'' derivatives.

Thus, we turn to proving that $\psit$ is rapidly decreasing.  In fact,
by the argument earlier in this proof, it suffices to show that
for each $M$, there exists an $N$ such that:
\begin{equation*}
2^{-\l| l-j\r|N - \l|l-k\r|N}\l|\dil{\psit}{2^l}\l( x\r)\r| \lesssim 2^{lQ} \l( 1+ 2^l \l|x\r|\r)^{-M}
\end{equation*}
And therefore, it suffices to show that if $l_0=\min\l\{l,j,k \r\}$, then
\begin{equation*}
\l| E_{2^j,2^k} \phi \l(x\r) \r|\lesssim 2^{l_0Q}\l( 1+2^{l_0}\l|x\r| \r)^{-M}
\end{equation*}
to do this, we will show that if we redefine $\psit$ to be:
\begin{equation*}
\dil{\psit}{2^{l_0}} =  E_{2^j,2^k} \dil{\phi}{2^l}
\end{equation*}
then we have that $\psit\in \sSz$ uniformly in the relevant parameters.
We proceed in the cases when $l_0=j$ or $l_0=l$.  The case when $l_0=k$
is similar to that when $l_0=j$.

We will next prove that
$$E_{2^{j},2^{k}} \Opl{\dil{\phi}{2^{l}}} = 2^{-\l| j-l \r|} E_{2^{l_0},2^{k}}$$
where $E_{2^{l_0},2^{k}}$ is a $2^{l_0},2^{k}$ elementary kernel.  Then,
the result will follow from the fact that:
$$E_{2^{j},2^{k}} \dil{\phi}{2^{l}} = \ker{E_{2^{j},2^{k}} \Opl{\dil{\phi}{2^{l}}}}\l( \cdot,0\r)$$
and applying Lemma \ref{LemmaElemKerAreSchwarz}.

We consider the identity operator $I$ as a right convolution operator.
Then, we may apply Proposition \ref{PropDecompOfCZ} to see that:
\begin{equation*}
I=\sum_{k_2} \Opr{\dil{\psi_{k_2}}{2^{k_2}}}
\end{equation*}
with this sum converging strongly in $L^2$.  It is easy to see
that everything we're dealing with in this proof is continuous on $L^2$,
hence,
\begin{equation*}
\begin{split}
E_{2^{j},2^{k}} & \Opl{\dil{\phi}{2^{l}}} = E_{2^{j},2^{k}} I \Opl{\dil{\phi}{2^{l}}}\\
&=\sum_{k_2} E_{2^{j},2^{k}} \Opr{\dil{\psi_{k_2}}{2^{k_2}}}\Opl{\dil{\phi}{2^{l}}}\\
&=\sum_{k_2} E_{2^{j},2^{k}} E_{2^{l},2^{k_2}}\\
&=\sum_{k_2} 2^{-\l| j-l \r|-\l|k-k_2\r|} E_{2^{l_0},2^{k}}\\
&=2^{-\l| j-l \r|} E_{2^{l_0},2^{k}}
\end{split}
\end{equation*}
where we have used Lemma \ref{LemmaCompElemKer} and Theorem \ref{ThmGensisElemKer}, completing the proof.
\end{proof}

\begin{thm}\label{ThmDecompOfsAp}
Suppose for each $j,k\in \Z$ we have $E_{2^j,2^k}$ a $2^j,2^k$ elementary
operator, uniformly in $j,k$.  Then,
\begin{equation}\label{EqnDecompsA}
T=\sum_{j,k} E_{2^j,2^k}
\end{equation}
converges in the topology of bounded convergence as operators 
$\sSz\rightarrow \sSz$, and also converges in the strong operator topology
as bounded operators $L^2\rightarrow L^2$.  Moreover, $T\in \sA'$.  Conversely,
every operator in $\sA'$ can be decomposed as in (\ref{EqnDecompsA}).
\end{thm}
\begin{proof}
The convergence of the sum
$$T=\sum_{j,k} E_{2^j,2^k}$$
in the topology of bounded convergence as operators $\sSz\rightarrow \sSz$
follows directly from Corollary \ref{CorElemKerOnsSz}, thinking of a fixed
element $\phi\in \sSz$ as $\phi=\dil{\phi}{2^0}$.  To see that the sum
(\ref{EqnDecompsA}) converges in the strong operator topology $L^2\rightarrow L^2$,
we apply the Cotlar-Stein lemma.  Indeed, the adjoint of a $2^j,2^k$
elementary operator is again a $2^j,2^k$ elementary operator (see Remark \ref{RmkDenfsAreSymmetric}), and
therefore,
\begin{equation*}
E_{2^{j_1},2^{k_1}}^{*}E_{2^{j_2},2^{k_2}} = 2^{-\l| j_1-j_2\r| -\l| k_1-k_2\r|} \Et_{2^{j_2},2^{k_2}}
\end{equation*}
by Lemma \ref{LemmaCompElemKer}.  Thus, to see that the sum converges
in the strong operator topology $L^2\rightarrow L^2$, it suffices
to show that the operators $\Et_{2^{j_2},2^{k_2}}$ are uniformly
bounded $L^2\rightarrow L^2$.  This follows easily from Lemma \ref{LemmaElemBoundedByMax}.  Alternatively, it is easy to see that
the $\Et_{2^{j_2},2^{k_2}}$ are uniformly NIS operators corresponding
to the metrics $\rl{2^{j_2-k_2}}$ (or $\rr{2^{k_2-j_2}}$ if $k_2\leq j_2$),
which in turn correspond to spaces which are uniformly spaces of homogeneous type in the sense of
\cite{DavidJourneSemmesOperateursDeCalderonZygmund}, by the
remarks in Section \ref{SectionCCBackground}.  Thus 
the uniform $L^2$ boundedness for $\Et_{2^{j_2},2^{k_2}}$ follows
by usual proofs that NIS operators are bounded on $L^2$ 
(see \cite{KoenigOnMaixmalSobolevAndHolderEstimates,NagelRicciSteinSingularIntegralsWithFlagKernels}).

To see $T\in \sA'$, note that:
\begin{equation*}
\begin{split}
T E_{2^{j_0},2^{k_0}} &= \sum_{j,k} E_{2^j,2^k} E_{2^{j_0},2^{k_0}}\\
&= \sum_{j,k} 2^{-\l| j-j_0 \r|-\l|k-k_0 \r|} E_{2^{j_0},2^{k_0}}\\
&= E_{2^{j_0},2^{k_0}}
\end{split}
\end{equation*}
where we have applied Lemma \ref{LemmaCompElemKer}, where $E_{2^{j_0},2^{k_0}}$
is a $2^{j_0},2^{k_0}$ elementary operator.  For more general $r_L,r_R$
elementary operators, merely think of an $r_L,r_R$ elementary operator
as a $2^{j_0},2^{k_0}$ elementary operator, where $j_0,k_0$ are chosen
to minimize $\l| 2^{j_0}-r_L \r|+\l|2^{k_0}-r_R\r|$.

For the converse, suppose $T\in \sA$.  Thinking of the identity
$I$ as a two-sided convolution operator, we may apply Proposition
\ref{PropDecompOfTwosided} to see that:
\begin{equation*}
I = \sum_{j,k\in \Z} \Opt{\dil{\phi_{j,k}}{2^j,2^k}}
\end{equation*}
where $\l\{ \phi_{j,k} \r\}\subset \SztSz$ is a bounded set.  Applying
Theorem \ref{ThmGensisElemKer}, we see:
\begin{equation*}
I = \sum_{j,k\in \Z} E_{2^j,2^k}
\end{equation*}
(where this sum converges strongly in $L^2$, and as we have seen earlier
in the proof, in the topology of bounded convergence $\sSz\rightarrow \sSz$).
Hence, we see:
\begin{equation*}
T= TI = \sum_{j,k\in \Z} TE_{2^j,2^k} = \sum_{j,k\in \Z} \widetilde{E}_{2^j,2^k}
\end{equation*}
completing the proof.
\end{proof}

\begin{cor}\label{CorTwosidedAreProduct}
Let $K$ be a product kernel.  Then, $\Opt{K}\in \sA'$.
\end{cor}
\begin{proof}
This is a combination of Proposition \ref{PropDecompOfTwosided}
and Theorems \ref{ThmGensisElemKer} and \ref{ThmDecompOfsAp}.
\end{proof}

\section{Equivalence of $\sA$ and $\sA'$}\label{SectionsAEquivsAp}
	In this section, we show that $\sA$ and $\sA'$ are the same spaces.
To begin with, we will need a better understanding of the function $B$
defined in Section \ref{SectionStatmentOfResults} by:
\begin{equation*}
\begin{split}
B& \l( r_L, r_R, N_L, N_R, m, x, y \r) \\
 &= r_L^{N_L} r_R^{N_R} \l( 1+ r_L \rlp{\frac{r_L}{r_R}}{x}{y} \r)^{-m} \frac{1}{\Vl{\frac{r_L}{r_R}}{x}{\frac{1}{r_L}+\rlp{\frac{r_L}{r_R}}{x}{y}}}
\end{split}
\end{equation*}

\begin{lemma}\label{LemmaSumB}
If $N_L,N_R\geq Q+m+1$,
we have,
\begin{equation}\label{EqnSumB}
\sum_{\substack{j\leq j_0\\k\leq k_0}} B\l( 2^j, 2^k, N_L,N_R, m, x,y \r) \approx B\l( 2^{j_0}, 2^{k_0}, N_L, N_R, m, x, y \r)
\end{equation}
and as a simple corollary:
\begin{equation*}
\sum_{j\leq j_0} B\l( 2^j, 2^{k_0}, N_L, N_R, m, x, y\r) \approx B\l( 2^{j_0},2^{k_0},N_L,N_R,m,x,y\r)
\end{equation*}
\end{lemma}
\begin{proof}
$\gtrsim$ is clear, and so we focus on $\lesssim$.  Without loss
of generality, we may assume $j_0\leq k_0$.  We separate the sum
(\ref{EqnSumB}) into the sum when $j\leq k$ and the sum when $k\leq j$.
We consider, first, the easier case when $j\leq k$:
\begin{equation*}
\begin{split}
&\sum_{\substack{j\leq k\\j\leq j_0\\k\leq k_0}} B\l( j,k,N_L,N_R,m,x,y\r) 
= \sum 2^{N_Lj +N_Rk} \frac{\l( 1+ 2^j\rlp{2^{j-k}}{x}{y}\r)^{-m}}{\Vl{2^{j-k}}{x}{2^{-j}+\rlp{2^{j-k}}{x}{y}}}\\
&\quad = \sum 2^{N_Lj+N_Rk}\frac{\l( 1+2^{j-j_0}2^{j_0} \rlp{2^{j_0-k_0}2^{k_0-k+j-j_0}}{x}{y} \r)^{-m}}{\Vl{2^{j-k}}{x}{2^{j_0-j}2^{-j_0} + \rlp{2^{j_0-k_0}2^{k_0-k+j-j_0}}{x}{y}}}
\end{split}
\end{equation*}
we now use the elementary fact that $\frac{1}{2} \rl{\epsilon} \leq \rl{2\epsilon}\leq \rl{\epsilon}$, to see:
\begin{equation*}
\begin{split}
&\lesssim \sum 2^{N_L j +N_R k} \frac{ \l( 1+ 2^{j-j_0+k-k_0}2^{j_0}\rlp{2^{j_0-k_0}}{x}{y} \r)^{-m}  }{\Vl{2^{j-k}}{x}{2^{k-k_0}\l( 2^{-j_0}+\rlp{2^{j_0-k_0}}{x}{y} \r)}}
\end{split}
\end{equation*}
Applying Corollary \ref{CorChangeVolScale}, we have:
\begin{equation*}
\lesssim 
\begin{cases}
\sum 2^{N_L j +N_R k+\l(j_0-j+k-k_0\r)Q}\frac{ \l( 1+ 2^{j-j_0+k-k_0}2^{j_0}\rlp{2^{j_0-k_0}}{x}{y} \r)^{-m}  }{\Vl{2^{j_0-k_0}}{x}{2^{k-k_0}\l( 2^{-j_0}+\rlp{2^{j_0-k_0}}{x}{y} \r)}} & \text{ if $j_0-k_0\geq j-k$}\\
\sum 2^{N_L j +N_R k+\l( j-j_0+k_0-k \r)Q} \frac{ \l( 1+ 2^{j-j_0+k-k_0}2^{j_0}\rlp{2^{j_0-k_0}}{x}{y} \r)^{-m}  }{\Vl{2^{j_0-k_0}}{x}{2^{j-j_0}\l( 2^{-j_0}+\rlp{2^{j_0-k_0}}{x}{y} \r)}} & \text{ if $j_0-k_0\leq j-k$}\\
\end{cases}
\end{equation*}
But for $\delta<1$, $\Vl{2^{j-k}}{x}{\delta r}\gtrsim \delta^Q \Vl{2^{j-k}}{x}{r}$ (this is a consequence of Theorem \ref{ThmEstBalls}), and so we have,
\begin{equation*}
\begin{split}
&\lesssim
\sum 2^{N_L j +N_R k+\l( j_0-j+k_0-k\r)Q}\frac{ \l( 1+ 2^{j-j_0+k-k_0}2^{j_0}\rlp{2^{j_0-k_0}}{x}{y} \r)^{-m}  }{\Vl{2^{j_0-k_0}}{x}{ 2^{-j_0}+\rlp{2^{j_0-k_0}}{x}{y} }} 
\\&\lesssim B\l( 2^{j_0},2^{k_0},N_L,N_R, m, x,y \r)
\end{split}
\end{equation*}
since this sum is geometric, provided $N_L,N_R\geq Q+m+1$.

We now turn to the sum when $k\leq j$:
\begin{equation}\label{EqnBBound1}
\begin{split}
&\sum_{\substack{k\leq j\\j\leq j_0\\k\leq k_0}} B\l( j,k,N_L,N_R, m, x,y\r) = \sum 2^{N_L j + N_R k} \frac{\l( 1+ 2^k\rrp{2^{k-j}}{x}{y}\r)^{-m}}{\Vr{2^{k-j}}{x}{2^{-k}+\rrp{2^{k-j}}{x}{y}}}
\end{split}
\end{equation}
Using the fact that $\rr{\epsilon} = \epsilon \rl{\frac{1}{\epsilon}}$ (here
we have removed the restriction $\epsilon\leq 1$) we see:
\begin{equation*}
\rr{2^{k-j}} = \rr{2^{k_0-j_0}2^{j_0-j+k-k_0}} \geq 2^{j_0-j} \rr{2^{k_0-j_0}} = 2^{k_0-j} \rl{2^{j_0-k_0}} 
\end{equation*}
Plugging this into (\ref{EqnBBound1}), we see that (\ref{EqnBBound1}) is
\begin{equation*}
\lesssim \sum 2^{N_Lj+N_Rk} \frac{\l( 1+ 2^{k_0-j_0+k-j} 2^{j_0}\rlp{2^{j_0-k_0}}{x}{y}\r)^{-m}}{\Vr{2^{k-j}}{x}{2^{-k}+ 2^{k_0-j}\rlp{2^{j_0-k_0}}{x}{y}}}
\end{equation*}
using that $j_0\leq k_0$,
\begin{equation*}
\lesssim \sum 2^{N_Lj+N_Rk} \frac{\l( 1+ 2^{j_0-j} 2^{j_0}\rlp{2^{j_0-k_0}}{x}{y}\r)^{-m}}{\Vr{2^{k-j}}{x}{2^{-k}+ 2^{k_0-j}\rlp{2^{j_0-k_0}}{x}{y}}}
\end{equation*}
Applying Corollary \ref{CorChangeVolScale} and using the fact that
the indicies we are summing over satisfy $k_0-k\geq j_0-j$ and so
we are in the lower case of (\ref{EqnChangeVolScale2}), and thus
\begin{equation*}
\lesssim \sum 2^{N_Lj+N_Rk+\l(k_0-k+j-j_0\r)Q} \frac{\l( 1+ 2^{j_0-j} 2^{j_0}\rlp{2^{j_0-k_0}}{x}{y}\r)^{-m}}{\Vl{2^{j_0-k_0}}{x}{2^{k_0-k}2^{-j_0}+ 2^{k_0-j+k_0-j_0}\rlp{2^{j_0-k_0}}{x}{y}}}
\end{equation*}
using that $j\leq j_0\leq k_0$, we see:
\begin{equation*}
\begin{split}
&\lesssim \sum 2^{N_Lj+N_Rk+\l(k_0-k\r)Q} \frac{\l( 1+ 2^{j_0-j} 2^{j_0}\rlp{2^{j_0-k_0}}{x}{y}\r)^{-m}}{\Vl{2^{j_0-k_0}}{x}{2^{-j_0}+ \rlp{2^{j_0-k_0}}{x}{y}}}\\
&\lesssim B\l( j_0,k_0,N_L,N_R,m,x,y\r)
\end{split}
\end{equation*}
provided $N_L,N_R\geq Q+m+1$, completing the proof of the first estimate.
The second estimate follows as a simple corollary.
\end{proof}

\begin{rmk}\label{RmkSmallN}
In our proof that $\sA'\subseteq \sA$, we will see that the only
reason we need to take $N_0$ large in Definition \ref{DefnsA} is
so that we may apply Lemma \ref{LemmaSumB}.  Because of this,
once we show that $\sA =\sA'$, we will see that we may
replace $N_0$ in Definition \ref{DefnsA} by $Q+m+1$.
\end{rmk}

\begin{lemma}\label{LemmaFirstElemOnBump}
Suppose $\phi_{2^{j_1},2^{k_1}}^{x}$ is a normalized $2^{j_1},2^{k_1}$ bump function of order $0$ centered at $x$, $\phi_{2^{j_2},2^{k_2}}^{z}$ is a normalized $2^{j_2},2^{k_2}$ bump
function of order $0$ centered at $z$, and $E_{2^{j_3},2^{k_3}}$ is a $2^{j_3},2^{k_3}$
elementary operator.  Then,
\begin{equation*}
\begin{split}
&\l|\Ltip{\phi_{2^{j_1},2^{k_1}}^x}{E_{2^{j_3},2^{k_3}} \phi_{2^{j_2},2^{k_2}}^z}\r|\\
&\quad\lesssim B\l( 2^{j_1\wedge j_2\wedge j_3}, 2^{k_1\wedge k_2\wedge k_3},0, 0,m,x,z \r)
\end{split}
\end{equation*}
with constants uniform in all the relevant parameters.
\end{lemma}
\begin{proof}
Let $j_0=j_1\wedge j_2$, and $k_0=k_1\wedge k_2$.  We prove the result in
the case when $k_3-k_0\geq j_3-j_0$, the other case being similar.  We
also proceed in the case when $k_3\geq j_3$, though the proof is 
essentially independent of this choice.
Consider, letting $\chi$ be as in Remark \ref{RmkConvoBoundBump},
\begin{equation}\label{EqnFirstElemBound}
\begin{split}
&\l|\Ltip{\phi_{2^{j_1},2^{k_1}}^x}{E_{2^{j_3},2^{k_3}} \phi_{2^{j_2},2^{k_2}}^z}\r|
\\&\lesssim \sum_{l\geq 0} \int_{\rlp{2^{j_3-k_3}}{y_1}{y_2}\approx 2^{l-j_3}} \phi_{2^{j_1},2^{k_1}}^x\l( y_1\r) 2^{-lN} \frac{1}{\Vl{2^{j_3-k_3}}{y_1}{2^{l-j_3}}} \phi_{2^{j_2},2^{k_3}}^z\l( y_2\r) dy_1 dy_2\\
&\lesssim \sum_{l\geq 0} 2^{-lN} {\rm Ker}\bigg( \Opl{\dil{\chi}{2^{j_1}}}\Opr{\dil{\chi}{2^{k_1}}}\Opl{\dil{\chi}{2^{j_3-l}}}\Opr{\dil{\chi}{2^{k_3-l}}} \\
&\quad\quad\quad\quad\quad\quad\quad\times\Opl{\dil{\chi}{2^{j_2}}}\Opr{\dil{\chi}{2^{k_2}}}\bigg)\l( x,z\r)
\end{split}
\end{equation}
at this point we may drop the assumption $k_3\geq j_3$, and note that
we could have just as easily shown (\ref{EqnFirstElemBound}) in the
case $k_3\leq j_3$.  In the above $N$ is any fixed integer we choose,
obtained from the rapid decrease of $\ker{E_{2^{j_3},2^{k_3}}}$.  Rearranging
terms, and using that, for instance,
$$\dil{\chi}{2^{j_1}}*\dil{\chi}{2^{j_3-l}}*\dil{\chi}{2^{j_2}} = \dil{\chit}{2^{j_0\wedge \l(j_3-l\r)}}$$
where $\chit$ is a bounded function of bounded support, with bounds independent
of all the relevant parameters above, we see that the left hand side of
(\ref{EqnFirstElemBound}) is
\begin{equation*}
\begin{split}
&\lesssim \sum_{l\geq 0} 2^{-lN}\ker{\Opl{\dil{\chit}{2^{j_0\wedge \l(j_3-l\r)}}}\Opr{\dil{\chitt}{2^{k_0\wedge \l(k_3-l\r)}}}} \l( x,z\r)
\\&\lesssim \sum_{l\geq 0} 2^{-lN} B\l( 2^{j_0\wedge \l(j_3-l\r)},2^{k_0\wedge\l(k_3-l\r)}, 0,0,m, x,z\r)
\end{split}
\end{equation*}
where we have applied (\ref{EqnSuppOfConvo}) and Theorem \ref{ThmBoundIntersection}.
We separate this sum into three sums.  The first:
\begin{equation*}
\sum_{0\leq l\leq k_3-k_0} 2^{-lN} B\l( 2^{j_0},2^{k_0}, 0,0,m, x,z\r) \lesssim B\l( 2^{j_0}, 2^{k_0},0,0,m,x,z\r)
\end{equation*}
with this sum $=0$ if $k_0> k_3$.  The second:
\begin{equation*}
\begin{split}
&\sum_{\l[\l(k_3-k_0\r) \vee 0\r]\leq l \leq j_3-j_0} 2^{-lN} B\l( 2^{j_0},2^{k_3-l},0,0,m,x,z\r) \\
&= \sum_{\l[\l(k_3-k_0\r)\vee 0\r]\leq l \leq j_3-j_0} 2^{-k_3N} B\l( 2^{j_0},2^{k_3-l},0,N,m,x,z\r)\\
&\approx 2^{-k_3N} B\l( 2^{j_0}, 2^{k_3\wedge k_0},0,N,m,x,z\r) \\
&\lesssim B\l( 2^{j_0},2^{k_3\wedge k_0},0,0,m,x,z\r)
\end{split}
\end{equation*}
with this sum $=0$ if $j_0>j_3$, and we have applied Lemma \ref{LemmaSumB} and
we have used that we may take $N$ large.  Finally,
\begin{equation*}
\begin{split}
&\sum_{\l[\l(j_3-j_0\r)\vee 0\r] \leq l} 2^{-lN} B\l( 2^{j_3-l},2^{k_3-l},m, 0,0,x,y\r) \\
&=\sum_{\l[\l(j_3-j_0\r)\vee 0\r] \leq l} 2^{-j_3N/2-k_3N/2} B\l( 2^{j_3-l},2^{k_3-l},m,\frac{N}{2},\frac{N}{2},m,x,y\r)\\
&\leq \sum_{\substack{\l[\l(j_3-j_0\r)\vee 0\r] \leq l_1\\ \l[\l(k_3-k_0\r) \vee 0\r]\leq l_2}} 2^{-j_3N/2-k_3N/2} B\l( 2^{j_3-l_1},2^{k_3-l_2},\frac{N}{2},\frac{N}{2},m,x,y\r)\\
&\lesssim 2^{-j_3N/2-k_3N/2} B\l( 2^{j_3\wedge j_0}, 2^{k_3\wedge k_0},  \frac{N}{2}, \frac{N}{2}, m,x,y\r)\\
&\lesssim B\l( 2^{j_3\wedge j_0}, 2^{k_3\wedge k_0},0,0,m,x,y\r) 
\end{split}
\end{equation*}
where again we have taken $N$ large and applied Lemma \ref{LemmaSumB}.
\end{proof}

\begin{cor}
Suppose $\phi_{2^{j_1},2^{k_1}}^{x}$ is a normalized $2^{j_1},2^{k_1}$ bump function centered
at $x$, $\phi_{2^{j_2},2^{k_2}}^{z}$ is a normalized $2^{j_2},2^{k_2}$ bump
function centered at $z$ (each of some large order, how large
will be implicit in the proof), and $E_{2^{j_3},2^{k_3}}$ is a $2^{j_3},2^{k_3}$
elementary operator.  Then,
\begin{equation*}
\begin{split}
&\l|\Ltip{\phi_{2^{j_1},2^{k_1}}^x}{\gl^{\alpha_1}\gr^{\beta_1}E_{2^{j_3},2^{k_3}}\gl^{\alpha_2}\gr^{\beta_2} \phi_{2^{j_2},2^{k_2}}^z
}\r|\\
&\quad\lesssim 
2^{\l( \l(j_1\wedge j_2\r)-j_3 \r)\wedge 0 +\l(\l(k_1\wedge k_2\r)-k_3\r)\wedge 0}\\&\quad\quad\quad\times B\l( 2^{j_1\wedge j_2\wedge j_3}, 2^{k_1\wedge k_2\wedge k_3},\l|\alpha_1\r|+\l|\alpha_2\r|,\l|\beta_1\r|+\l|\beta_2\r|, m, x,z\r)
\end{split}
\end{equation*}
with constants uniform in all the relevant parameters.
\end{cor}
\begin{proof}
Let $j_0=j_1\wedge j_2$ and $k_0=k_1\wedge k_2$.  We first prove the
result without the factor of:
$$2^{\l(j_0-j_3\r)\wedge 0+\l(k_0-k_3\r)\wedge 0}$$
Suppose that $j_2=j_0\wedge j_3$.  Then, we have:
\begin{equation*}
\begin{split}
&\Ltip{\phi_{2^{j_1},2^{k_1}}^x}{\gl^{\alpha_1}\gr^{\beta_1}E_{2^{j_3},2^{k_3}}\gl^{\alpha_2}\gr^{\beta_2} \phi_{2^{j_2},2^{k_2}}^z}
\\&= \sum_{\l|\alpha_3\r|=\l|\alpha_1\r|}\Ltip{\phi_{2^{j_1},2^{k_1}}^x}{\gr^{\beta_1}\Et_{2^{j_3},2^{k_3},\alpha_3}\gl^{\alpha_3}\gl^{\alpha_2}\gr^{\beta_2} \phi_{2^{j_2},2^{k_2}}^z}
\\&=2^{j_2\l( \l|\alpha_1\r|+\l|\alpha_2\r|\r)} \sum_{\l|\alpha_3\r|=\l|\alpha_1\r|}\Ltip{\phi_{2^{j_1},2^{k_1}}^x}{\gr^{\beta_1}\Et_{2^{j_3},2^{k_3},\alpha_3}\gr^{\beta_2} \phit_{2^{j_2},2^{k_2},\alpha_3}^z}
\end{split}
\end{equation*}
a finite sum of terms of the same form but with $\l|\alpha_1\r|=0=\l|\alpha_2\r|$, times 
$2^{j_2 \l( \l| \alpha_1 \r| + \l| \alpha_2 \r| \r)} $.  
We get a similar result when $j_1=j_0\wedge j_3$.  Finally, when
$j_3=j_3\wedge j_0$, we merely let all the $\gl$ derivatives land on the $E_{2^{j_3},2^{k_3}}$,
\begin{equation*}
\begin{split}
&\Ltip{\phi_{2^{j_1},2^{k_1}}^x}{\gl^{\alpha_1}\gr^{\beta_1}E_{2^{j_3},2^{k_3}}\gl^{\alpha_2}\gr^{\beta_2} \phi_{2^{j_2},2^{k_2}}^z}\\
&= 2^{j_3 \l( \l| \alpha_1 \r| + \l| \alpha_2 \r| \r) } \Ltip{\phi_{2^{j_1},2^{k_1}}^x}{\gr^{\beta_1}\Et_{2^{j_3},2^{k_3}}\gr^{\beta_2} \phi_{2^{j_2},2^{k_2}}^z}
\end{split}
\end{equation*}
Doing a similar proof with the $k$s, we see:
\begin{equation*}
\begin{split}
&\l|\Ltip{\phi_{2^{j_1},2^{k_1}}^x}{\gl^{\alpha_1}\gr^{\beta_1}E_{2^{j_3},2^{k_3}}\gl^{\alpha_2}\gr^{\beta_2} \phi_{2^{j_2},2^{k_2}}^z}\r|\\
&\leq 2^{j_0\wedge j_3\l( \l| \alpha_1 \r| + \l| \alpha_2 \r| \r) k_0\wedge k_3\l(\l|\beta_1\r|+\l|\beta_2\r|\r)} \sum \Ltip{\phit_{2^{j_1},2^{k_1}}^x}{\Et_{2^{j_3},2^{k_3}} \phit_{2^{j_2},2^{k_2}}^z}
\end{split}
\end{equation*}
where the sum denotes a finite sum of such terms.  Applying Lemma \ref{LemmaFirstElemOnBump}, we see:
\begin{equation*}
\begin{split}
&\lesssim  2^{j_0\wedge j_3\l( \l| \alpha_1 \r| + \l| \alpha_2 \r| \r) k_0\wedge k_3\l(\l|\beta_1\r|+\l|\beta_2\r|\r)} B\l( 2^{j_0\wedge j_3},2^{k_0\wedge k_3}, 0,0,m, x,z\r)
\\&= B\l( 2^{j_0\wedge j_3},2^{k_0\wedge k_3},\l| \alpha_1 \r| + \l| \alpha_2\r|, \l|\beta_1\r|+\l|\beta_2\r|, m, x,z \r)
\end{split}
\end{equation*}
Which completes the proof, without the factor of $2^{\l(j_0-j_3\r)\wedge 0+\l(k_0-k_3\r)\wedge 0}$.  To see how to obtain that factor, suppose we are
in the case when $j_0=j_1\leq j_3$.  Then we ``pull derivatives out'' of
$E_{2^{j_3},2^{k_3}}$ and let them land on $\phi_{2^{j_1},2^{k_1}}^x$; indeed,
\begin{equation*}
\begin{split}
&\Ltip{\phi_{2^{j_1},2^{k_1}}^x}{\gl^{\alpha_1}\gr^{\beta_1}E_{2^{j_3},2^{k_3
}}\gl^{\alpha_2}\gr^{\beta_2} \phi_{2^{j_2},2^{k_2}}^z}\\
&=\sum_{\l|\alpha\r|=1} 2^{-j_3}\Ltip{\phi_{2^{j_1},2^{k_1}}^x}{\gl^{\alpha_1}\gl^{\alpha}\gr^{\beta_1}E_{2^{j_3},2^{k_3}}\gl^{\alpha_2}\gr^{\beta_2} \phi_{2^{j_2},2^{k_2}}^z}\\
&= \sum_{\l|\alpha_3\r|=\l|\alpha_1\r|} 2^{j_1-j_3} \Ltip{\phit_{2^{j_1},2^{k_1},\alpha_3}^x}{\gl^{\alpha_3}\gr^{\beta_1}E_{2^{j_3},2^{k_3}}\gl^{\alpha_2}\gr^{\beta_2} \phi_{2^{j_2},2^{k_2}}^z}
\end{split}
\end{equation*}
which is $2^{j_1-j_3}=2^{j_0-j_3}$ times a finite sum of terms of the
original form.  A similar proof works for when $j_2=j_0$ and for the $k$s.
\end{proof}

\begin{thm}
Suppose $T\in \sA'$, then $T\in \sA$.
\end{thm}
\begin{proof}
We apply Theorem \ref{ThmDecompOfsAp} to decompose $T$:
$$T=\sum_{j,k\in \Z} E_{2^j,2^k}$$
where $E_{2^j,2^k}$ are uniformly $2^j,2^k$ elementary operators, and
this sum converges in the strong operator topology as operators $L^2\rightarrow L^2$.
Fix $m$ and fix $N_L,N_R\geq Q+m+1$.
Suppose $\phi_{2^{j_1},2^{k_1}}^{x}$ is a normalized $2^{j_1},2^{k_1}$ bump function centered
at $x$, $\phi_{2^{j_2},2^{k_2}}^{z}$ is a normalized $2^{j_2},2^{k_2}$ bump
function centered at $z$ (each of some large order), and suppose
that $\l|\alpha_1\r|+\l|\alpha_2\r|=N_L$, $\l|\beta_1\r|+\l|\beta_2\r|=N_R$.  Then, letting $j_0=j_1\wedge j_2$, $k_0=k_1\wedge k_2$, we see:
\begin{equation*}
\begin{split}
&\l|\Ltip{\phi_{2^{j_1},2^{k_1}}^x}{\gl^{\alpha_1}\gr^{\beta_1}T\gl^{\alpha_2}\gr^{\beta_2} \phi_{2^{j_2},2^{k_2}}^z}\r|\\
&\leq \sum_{j,k} \l|\Ltip{\phi_{2^{j_1},2^{k_1}}^x}{\gl^{\alpha_1}\gr^{\beta_1}E_{2^j,2^k}\gl^{\alpha_2}\gr^{\beta_2} \phi_{2^{j_2},2^{k_2}}^z}\r|\\
&\lesssim \sum_{j,k}2^{\l(j_0-j\r)\wedge 0 +\l( k_0-k\r)\wedge 0} B\l( 2^{j_0\wedge j}, 2^{k_0\wedge k}, N_L,N_R, m, x,z \r)\\
&\lesssim B\l( 2^{j_0}, 2^{k_0}, N_L,N_R, m, x, z\r)
\end{split}
\end{equation*}
where we have applied Lemma \ref{LemmaSumB} to get the last line,
completing the proof.
\end{proof}

We now turn to showing that $\sA\subset \sA'$.  Fix $T\in \sA$.  Then,
we wish to show that $TE_{r_L,r_R} = \Et_{r_L,r_R}$.  
As we have seen before, it will suffice to prove this result for
$r_L=2^j$, $r_R=2^k$.  This follows by choosing $j,k$ to minimize
$\l|r_L-2^j\r|+\l|r_R-2^k\r|$.

\begin{lemma}\label{LemmasAPointwise}
Given $r_L,r_R>0$, $y\in G$, $m\geq 0$, there exists a $N_0$ such that if 
$N_L,N_R\geq N_0$,
and $\l|\alpha_1\r|+\l|\alpha_2\r|=N_L$, $\l|\beta_1\r|+\l|\beta_2\r|=N_R$,
and $\phi_{r_L,r_R}^y$ is a $r_L,r_R$ bump function centered at $y$,
\begin{equation*}
\l|\gl^{\alpha_1}\gr^{\beta_1}T\gl^{\alpha_2}\gr^{\beta_2} \phi_{r_L,r_R}^y \l(x\r)\r|\lesssim B\l( r_L,r_R, N_L,N_R, m, x,y\r)
\end{equation*}
\end{lemma}
\begin{proof}
This follows directly from Definition \ref{DefnsA}, by taking $\phi_{\dil{r_L}{1},\dil{r_R}{1}}\rightarrow \delta_x$, by taking $\dil{r_L}{1},\dil{r_R}{1}\rightarrow \infty$.

To see that we can do this, merely take $\phi$ supported in the unit ball
$B$ such that $\int \phi =1$.  Then, 
$$\Opl{\dil{\phi}{2^j}}\rightarrow I$$
as $j\rightarrow \infty$; similarly for $\Opr{\dil{\phi}{2^k}}$.  Thus
if we set:
$$\phi_{2^j,2^k}^x\l(z\r) = \ker{\Opl{\dil{\phi}{2^j}}\Opr{\dil{\phi}{2^k}}}\l(x,z\r)$$
we see that $\phi_{2^j,2^k}^x\rightarrow \delta_x$.  Since we saw in 
Section \ref{SectionBumpFuncsAndElem} $\phi_{2^j,2^k}^x$ is essentially
an normalized bump function (it may really have support in a 
ball with radius a constant factor times the ball it is supposed to 
be supported in,
and need to be multiplied by a constant, but these only affect the
answer by a constant), we are done.
\end{proof}

\begin{prop}\label{PropsApSatisGrowth}
$\ker{TE_{2^j,2^k}}$ satisfies the estimates (\ref{EqnEleZero})
with $r_L=2^j$, $r_R=2^k$,
uniformly in the relevant parameters.
\end{prop}
\begin{proof}
We proceed in the case when $j\leq k$, the other case being similar.
Consider,
\begin{equation*}
\begin{split}
&\glx^{\alpha_1}\glx^{\beta_1}\glz^{\alpha_2}\glz^{\beta_2}\ker{TE_{2^j,2^k}}\l( x,z\r)\\
&= \glx^{\alpha_1}\glx^{\beta_1}\l( -1\r)^{\l| \alpha_2\r|+\l|\beta_2\r|}\ker{T E_{2^j,2^k}\gl^{\alpha_2}\gr^{\beta_2}}\l( x,z\r)\\
&= 2^{j\l|\alpha_2\r|+k\l|\beta_2\r|} \glx^{\alpha_1}\glx^{\beta_1} \ker{T \Et_{2^j,2^k}} \l( x,z\r)
\end{split}
\end{equation*}
which is $2^{j\l|\alpha_2\r|+k\l|\beta_2\r|}$ times a term of the original form.
Thus, it will suffice to prove the result when $\l|\alpha_2\r|=0=\l|\beta_2\r|$.

Fix $z\in G$.  Let us consider the function of $x$ given by:
\begin{equation*}
\begin{split}
&\glx^{\alpha_1}\grx^{\beta_1} \ker{T E_{2^j,2^k}}\l(x,z\r) \\&= 2^{-j\l|\alpha_2\r|-k\l|\beta_2\r|}\ker{\gl^{\alpha_1}\gr^{\beta_1}T\gl^{\alpha_2}\gr^{\beta_2}\Et_{2^j,2^k}}\l( x,z\r)
\end{split}
\end{equation*}
Here $\alpha_2$ and $\beta_2$ are not the same ordered multi-incides as before;
rather, we have applied Definition \ref{DefnElemKernels}, and the term
on the right hand side of the above equation really denotes a finite sum
of such terms.  Letting $\phi_{2^j,2^k} \l( x\r) = \ker{\Et_{2^j,2^k}}\l( x,z\r)$, we are considering
the function given by:
\begin{equation*}
2^{-j\l|\alpha_2\r|-k\l|\beta_2\r|} \gl^{\alpha_1}\gr^{\beta_1}T\gl^{\alpha_2}\gr^{\beta_2} \phi_{2^j,2^k}
\end{equation*}
where $\phi_{2^j,2^k}$ is a $2^j,2^k$ elementary kernel, and $\l|\alpha_2\r|,\l|\beta_2\r|$
can be as large as we like.

Theorem \ref{ThmBumpFuncs} allows us to create a partition of unity $\psi_l\l( x\r)$ ($l\geq 0$), such that:
\begin{enumerate}
\item $0\leq \psi_l \leq 1$, for every $l$
\item $\psi_0$ is supported where $\rlp{2^{j-k}}{x}{z}\lesssim 2^{-j}$
\item $\psi_l$ is supported where $\rlp{2^{j-k}}{x}{z}\approx 2^{-j+l}$, $l\ne 0$
\item $\l| \gl^\alpha \gr^\beta \phi\r|\lesssim 2^{\l| \alpha\r|\l(j-l\r)+\l|\beta\r|\l(k-l\r)}$ and, if $\l|\alpha\r|+\l|\beta\r|>0$, is supported where
$\rlp{2^{j-k}}{x}{z}\approx 2^{l-j}$, even if $l=0$.
\end{enumerate}
This follows from Theorem \ref{ThmBumpFuncs} directly for $z$ in the
closed unit ball, and for $2^{-j+l}$ small.  However Theorem \ref{ThmBumpFuncs}
really holds for all points in $G$ and all distances.  Creating a bump
function of radius $r$ centered at $z$ is equivalent to creating a bump
function of radius $\delta$ centered at $\frac{\delta}{r}z$.  Thus
Theorem \ref{ThmBumpFuncs} extends to all points and all radii, by
homogeneity (just take $\delta$ small enough), giving us the above partition
of unity.

Define $\phi_{l}\l(x\r) = \psi_l\l(x\r) \phi_{2^j,2^k}\l(x,z\r)$ (thinking of $z$ as fixed), so that $\sum_{l\geq 0} \phi_l = \phi_{2^j,2^k}$.  We wish to show that $\phi_l$ is $2^{-lN}$ times a $2^{j-l},2^{k-l}$ normalized
bump function, where $N$ is any integer we choose, and we really mean
a constant times a normalized bump function, with support in, perhaps,
a constant times the radius of the support it's supposed to have.

We already know that the support of $\phi_l$ is correct, by the properties
of $\psi_l$, so we turn to estimating derivatives of $\phi_l$.  When
$l>0$, we have:
\begin{equation*}
\begin{split}
\l|\gl^\alpha \gr^\beta \phi_l\r| &\leq \sum_{\substack{a_1+a_2=\l|\alpha\r|\\b_1+b_2=\l|\beta\r|}} 2^{\l(a_1\l(j-l\r) +b_1\l(k-l\r) \r)} \frac{2^{a_2j+b_2k}}{2^{lN} \Vl{2^{j-k}}{z}{2^{l-j}}}\\
&\lesssim \frac{2^{\l|\alpha\r|\l( j-l \r) + \l| \beta\r|\l(k-l\r)}}{2^{lN'} \Vl{2^{j-k}}{z}{2^{l-j}}}
\end{split}
\end{equation*}
here, $a_1$ represents the number of $\gl$ derivatives that land on $\psi_l$,
and $a_2$ represents the number that land on $\phi_{2^j,2^k}$.  We have
used in the last line that we make take $N$ large.  This establishes
that $\phi_l$ ($l>0$) is $2^{-lN}$ times a $2^j,2^k$ bump function.
When $l=0$, a nearly identical proof establishes the result.

Thus, by Lemma \ref{LemmasAPointwise}, and using the fact that we may
take $\l|\alpha_2\r|,\l|\beta_2\r|$ as large as we like, we have:
\begin{equation*}
\begin{split}
&2^{-j\l|\alpha_2\r|-k\l|\beta_2\r|} \l| \gl^{\alpha_1}\gr^{\beta_1} T \gl^{\alpha_2} \gr^{\beta_2} \phi_l  \r| \\
&\quad\quad\lesssim 2^{-j\l|\alpha_2\r|-k\l|\beta_2\r|-lN} B\l( 2^{j-l}, 2^{k-l}, \l|\alpha_1\r| + \l|\alpha_2\r|, \l|\beta_1\r|+\l|\beta_2\r|, m, x, z\r)
\\&\quad\quad \leq 2^{-lN} B\l( 2^{j-l}, 2^{k-l}, \l|\alpha_1\r|, \l|\beta_1\r|, m, x, z\r)
\\&\quad\quad \leq 2^{-lN'} B\l( 2^{j}, 2^k, \l|\alpha_1\r|, \l|\beta_1\r|, m, x, z\r)
\end{split}
\end{equation*}
Where we have used in the last step that the function:
$$2^{-l(m+1)} B\l(2^{j-l},2^{k-l}, \l|\alpha_1\r|, \l|\beta_1\r|, m, x, z \r)$$
decreases as $l$ increases.
Hence,
\begin{equation*}
2^{-j\l|\alpha_2\r|-k\l|\beta_2\r|} \l| \gl^{\alpha_1}\gr^{\beta_1}T\gl^{\alpha_2}\gr^{\beta_2} \phi_{2^j,2^k}\r| \lesssim  B\l( 2^{j}, 2^k, \l|\alpha_1\r|, \l|\beta_1\r|, m, x, z\r)
\end{equation*}
completing the proof.
\end{proof}

Proposition \ref{PropsApSatisGrowth} shows that $\ker{TE_{2^j,2^k}}$ satisfies the
growth estimates of a $2^j,2^k$ elementary kernel.  Thus, to show that
$TE_{2^j,2^k}$ is a $2^j,2^k$ elementary operator, it now remains
to show that we may ``pull out'' derivatives, as in Definition \ref{DefnElemKernels}.
To do this, it will suffice to show that the class of operators $\sA$
commutes with $\gl$ and $\gr$ derivatives.  By this, we mean:

\begin{thm}\label{ThmsACommutesWithDeriv}
Suppose $T\in \sA$.  Then,
\begin{equation*}
T\gl^{\alpha}\gr^{\beta} = \sum_{\substack{\l|\alpha_1\r|=\l|\alpha\r|\\\l|\beta_1\r|=\l|\beta\r|}} \gl^{\alpha_1}\gr^{\beta_1} T_{\alpha_1,\beta_1}
\end{equation*}
where $T_{\alpha_1,\beta_1}\in \sA$.
\end{thm}

To see why Theorem \ref{ThmsACommutesWithDeriv} completes the proof
that $\sA\subseteq \sA'$, consider:
\begin{equation*}
\begin{split}
TE_{2^j,2^k} &= \sum_{\substack{\l|\alpha_1\r|=N_1\\\l|\alpha_2\r|=N_2\\\l|\beta_1\r|=N_3\\\l|\beta_2\r|=N_4}} 2^{-j\l(N_1+N_2\r) -k\l(N_3+N_4\r)}T\gl^{\alpha_1}\gr^{\beta_1} \Et_{2^j,2^k} \gl^{\alpha_2} \gr^{\beta_2}\\
&= \sum 2^{-j\l(N_1+N_2\r) -k\l(N_3+N_4\r)} \gl^{\alpha_1}\gr^{\beta_1} T_{\alpha_1,\beta_1} \Et_{2^j,2^k} \gl^{\alpha_2} \gr^{\beta_2}
\end{split}
\end{equation*}
a finite sum of terms satisfying the proper bounds associated to elementary kernels.  Hence,
we conclude this section by proving Theorem \ref{ThmsACommutesWithDeriv}.

\begin{proof}[Proof of Theorem \ref{ThmsACommutesWithDeriv}]
We will show that if $X_L=\gl^{\alpha}$, $\l|\alpha\r|=1$,  then
$TX_L = \sum_{\l|\alpha_1\r|=1} \gl^{\alpha_1} T_{\alpha_1}$,
and the whole result will follow by symmetry and induction.
Let $\sJ$ be the homogeneous fundamental solution to the sublaplacian:
$$\sum_{\l|\alpha\r|=1} \l( \gl^{\alpha}\r)^2$$
See \cite{FollandSubellipticEstimatesAndFunctionsSpacesOnNilpotentLieGroups}
for background on $\sJ$.  Note that:
\begin{equation*}
TX_L= \sum_{\l|\alpha\r|=1} \l( \gl^{\alpha}\r)^2 \sJ TX_L= \sum_{\l|\alpha\r|=1} \gl^{\alpha} \l(\l(\gl^{\alpha}\sJ\r) TX_L\r)
\end{equation*}
Thus, it suffices to show that $STX_L\in \sA$, where $S$ is a left invariant
convolution operator, with kernel of type $1$, in the sense of
\cite{FollandSubellipticEstimatesAndFunctionsSpacesOnNilpotentLieGroups}.
Hence, we wish to estimate terms like:
\begin{equation*}
\begin{split}
\Ltip{\phi_{\dil{r_L}{1},\dil{r_R}{1}}^x}{\gl^{\alpha_1}\gr^{\beta_1} STX_L \gl^{\alpha_2}\gr^{\beta_2} \phi_{\dil{r_L}{2},\dil{r_R}{2}}^y}
\end{split}
\end{equation*}
where everything above is as in Definition \ref{DefnsA}.
However, $\gr^{\beta_1} S =S \gr^{\beta_1}$ and $\gl^{\alpha_1}S = \sum_{\l|\alpha_1'\r|=\l|\alpha_1\r|} S_{\alpha_1'} \gl^{\alpha_1'}$,
where $S_{\alpha_1'}$ is a left invariant operator with convolution kernel of type $1$.
Thus, it suffices to bound terms of the form:
\begin{equation*}
\begin{split}
\Ltip{\phi_{\dil{r_L}{1},\dil{r_R}{1}}^x}{S\gl^{\alpha_1}\gr^{\beta_1} T X_L\gl^{\alpha_2}\gr^{\beta_2} \phi_{\dil{r_L}{2},\dil{r_R}{2}}^y}
\end{split}
\end{equation*}
where $S$ is an operator with convolution kernel of type $1$.  (It is easy to see
that all the integrals involved converge absolutely, by Lemma \ref{LemmasAPointwise}.)
Let $\psi$ be a $C_0^\infty$ bump function supported on the unit ball,
which is $1$ on the ball of radius $3/4$, and $0\leq\phi \leq 1$.
Let $K\l( x\r)$ be the convolution kernel of $S$, and define:
$\phi \l( x\r) = \l( \phi\l( x\r) - \phi\l( 2x\r) \r) K\l( x\r)$.
Then,
$$S=\sum_{j\in \Z} 2^{-j}\Opl{\dil{\phi}{2^j}}$$
Applying Lemma \ref{LemmaConvoOnBump}, we see:
\begin{equation}\label{EqnsACommute1}
\begin{split}
&\l| \Ltip{\phi_{2^{j_1},2^{k_1}}^x}{S\gl^{\alpha_1}\gr^{\beta_1} T X_L\gl^{\alpha_2}\gr^{\beta_2} \phi_{2^{j_2},2^{k_2}}^y} \r|\\
&\leq \sum_{j\in\Z} 2^{-j}\l| \Ltip{\Opl{\dil{\phi}{2^j}}\phi_{2^{j_1},2^{k_1}}^x}{\gl^{\alpha_1}\gr^{\beta_1} T X_L\gl^{\alpha_2}\gr^{\beta_2} \phi_{2^{j_2},2^{k_2}}^y} \r|\\
&\lesssim \sum_{j\in\Z} 2^{-j}\l| \Ltip{\phit_{2^{j_1\wedge j},2^{k_1}}^x}{\gl^{\alpha_1}\gr^{\beta_1} T X_L\gl^{\alpha_2}\gr^{\beta_2} \phi_{2^{j_2},2^{k_2}}^y} \r|\\
&\lesssim 2^{-j}\sum_{j\in \Z} B\l( 2^{j_1\wedge j\wedge j_2}, 2^{k_1\wedge k_2}, \l|\alpha_1\r|+\l|\alpha_2\r|+1,\l|\beta_1\r|+\l|\beta_2\r|, x,y \r)
\end{split}
\end{equation}
Let $j_0=j_1\wedge j_2$, $k_0=k_1\wedge k_2$, $a=\l|\alpha_1\r|+\l|\alpha_2\r|$,
$b=\l|\beta_1\r|+\l|\beta_2\r|$.  Then, we separate the
sum on the left hand side of (\ref{EqnsACommute1}) into two sums:
\begin{equation*}
\begin{split}
\sum_{j\geq j_0} 2^{-j} B\l( 2^{j_0},2^{k_0},a+1,b,m,x,y  \r) 
&= 2^{-j_0} B\l( 2^{j_0},2^{k_0},a+1,b,m,x,y  \r)\\
&= B\l( 2^{j_0},2^{k_0},a,b,m,x,y  \r)
\end{split}
\end{equation*}
and,
\begin{equation*}
\begin{split}
\sum_{j\leq j_0} 2^{-j} B\l( 2^{j},2^{k_0},a+1,b,m,x,y  \r) 
&=\sum_{j\leq j_0} B\l( 2^{j},2^{k_0},a,b,m,x,y  \r) \\
&\approx B\l( 2^{j_0},2^{k_0},a,b,m,x,y \r)
\end{split}
\end{equation*}
where we have used that we may take $a$ and $b$ large, and we have
applied Lemma \ref{LemmaSumB}.  This completes the proof.
\end{proof}

\begin{rmk}
Theorem \ref{ThmsACommutesWithDeriv} is the only place where we use
the crucial hypothesis that we have a cancellation condition that
happens on \it both \rm sides of $T$ at once.
\end{rmk}

\section{$L^p$ Boundedness}\label{SectionLpBounded}
	In this section, we show that operators in $\sA'$ extend to bounded
operators on $L^p$, $1<p<\infty$.  To do this, we will need a relevant
Littlewood-Paley square function, and a relevant maximal function.  Fortunately,
we will be able construct both out of the building blocks of the analogous operators
for left and right convolution operators.

We begin with the maximal functions.  Define:
\begin{equation*}
\begin{split}
\l(\ML \l(f\r)\r)\l(x\r) &= \sup_{R>0} \frac{1}{R^Q} \int_{\l\| y^{-1}x\r\|<R} \l|f\l( y\r)\r| dy\\
&= \sup_{R>0} \Opl{\dil{\chi_B}{\frac{1}{R}}} \l|f\r|
\end{split}
\end{equation*}
and similarly,
\begin{equation*}
 \MR \l( f\r)  = \sup_{R>0} \Opr{\dil{\chi_B}{\frac{1}{R}}} \l| f\r|
\end{equation*}
It follows from the results in \cite{SteinHarmonicAnalysis} that
$$\LpNormG{\ML f}{p}\lesssim \LpNormG{f}{p}$$
for $1<p<\infty$, and similarly for $\MR$.
For us, the relevant maximal function will be:
\begin{equation*}
 \M f = \sup_{R_1>0, R_2>0} \Opl{\dil{\chi_B}{\frac{1}{R_1}}} \Opr{\dil{\chi_B}{\frac{1}{R_2}}}\l| f\r|
\end{equation*}
It is easy to see that 
\begin{equation}\label{EqnBoundMaximal}
\M f \leq \ML\MR f
\end{equation}
and therefore,
\begin{equation*}
\LpNormG{\M f}{p}\lesssim \LpNormG{f}{p}
\end{equation*}
Corresponding to each $1\geq\epsilon>0$, we get a maximal function for $\rl{\epsilon}$
(and one for $\rr{\epsilon}$, but let's focus on $\rl{\epsilon}$), defined
by:
\begin{equation*}
\sM_{\epsilon} f \l(x\r) = \sup_{R>0} \frac{1}{\Vl{\epsilon}{x}{R}}\int_{\rlp{\epsilon}{x}{y}\leq R} \l| f\l( y\r)\r| dy
\end{equation*}
But, then, taking $\chi$ as in Remark \ref{RmkConvoBoundBump}, we see:
\begin{equation*}
\sM_{\epsilon} f \lesssim \sup_{R>0} \Opl{\dil{\chi}{\frac{1}{R}}}\Opr{\dil{\chi}{\frac{1}{\epsilon R}}} \l| f\r| \lesssim \M f
\end{equation*}
so we see that $\M$ uniformly bounds the maximal functions corresponding to all
of the geometries we are considering.

\begin{lemma}\label{LemmaElemBoundedByMax}
Suppose $\ker{E_{r_L,r_R}}$ satisfies the bounds (\ref{EqnEleZero}) of 
Definition \ref{DefnElemKernels} without any derivatives (ie it is like an elementary kernel,
but we need not be able to ``pull out'' derivatives or take derivatives).  Then, for
$f\in \sS$, $\l| E_{r_L,r_R}f\r|\lesssim \M f$, with constants uniform
in all the relevant parameters.
\end{lemma}
\begin{proof}
We prove this in the case $r_L=2^j$, $r_R=2^k$, the more general case
following from this one.  We assume $k\geq j$, the other case following
in the same manner.
Consider,
\begin{equation*}
\begin{split}
\l| E_{2^j,2^k}f\l(x\r)\r| &\lesssim \int_{\rlp{2^{j-k}}{x}{y}\leq 2^{-j}} \frac{\l| f\l(y\r)\r|}{ \Vl{2^{k-j}}{x}{2^{-j}}} dy \\&\quad\quad+ \sum_{l\geq 1} \int_{\rlp{2^{j-k}}{x}{y}\approx 2^{l-j}} \frac{\l| f\l(y\r)\r|}{ 2^{lN}\Vl{2^{k-j}}{x}{2^{l-j}}} dy\\
&\leq \sum_{l\geq 0} 2^{-lN} \sM_{2^{k-j}} f\\
&\lesssim \M f
\end{split}
\end{equation*}
completing the proof.
\end{proof}

Recall the definition $\sS_{(m)}$.  $\phi\in \sS_{(m)}$ if and only if
$\phi = \sum_{\l| \alpha \r|=m} \gl^{\alpha} \psi_{\alpha}$, where
$\psi\in \sS$.  We get essentially the same space if we replace
$\gl$ by $\gr$ (by that we mean $\sS_{(m_1)}^R \subset \sS_{(m)}^L \subset \sS_{(m_2)}^R$ where $m_1,m_2\rightarrow \infty$ as $m\rightarrow \infty$).
In short, being a high order of $\gl$ derivatives is the same as being
a high order of $\gr$ derivatives, which is the same as moments
up to a high order vanishing.
We have:

\begin{lemma}\label{LemmaReproducing}
For any $N\in \N$ there exist function $\phi_1,\ldots, \phi_M, \psi_1,\ldots,\psi_M\in \sS_{(N)}$ (here $M$ depends on $N$) such that:
\begin{equation*}
\sum_{l=1}^M \sum_{j\in \Z} \dil{\phi_l}{2^j} *\dil{\psi_l}{2^j} = \delta_0
\end{equation*}
\end{lemma}
\begin{proof}
This follows directly from Theorem 1.61 in \cite{FollandSteinHardySpacesOnHomogeneousGroups}.
\end{proof}

We will be able to use the $\psi_l$ and $\phi_l$ from Lemma \ref{LemmaReproducing}
to construct a relevant Littlewood-Paley square function.  Henceforth, we fix such
$\psi_l$ and $\phi_l$, thinking of $N$ as large (how large $N$ will
have to be will be implicit in our proof).

We define
\begin{equation*}
\begin{split}
\Lambda_{j,k}^{l_1,l_2} &= \Opl{\dil{\phi_{l_1}}{2^j}}\Opr{\dil{\psi_{l_2}}{2^k}}\\
P_{j,k}^{l_1,l_2} &= \Opl{\dil{\psi_{l_1}}{2^j}}\Opr{\dil{\phi_{l_2}}{2^k}}
\end{split}
\end{equation*}
so that,
\begin{equation*}
\sum_{l_1,l_2} \sum_{j,k\in \Z} P_{j,k}^{l_1,l_2} \Lambda_{j,k}^{l_1,l_2} = I
\end{equation*}
and we define our square function:
\begin{equation*}
\Lambda\l( f\r) = \l(\sum_{l_1,l_2}\sum_{j,k\in \Z} \l| \Lambda_{j,k}^{l_1,l_2} f\r|^2\r)^{\frac{1}{2}}
\end{equation*}

\begin{thm}
For $1<p<\infty$,
\begin{equation*}
\LpNormG{f}{p}\approx \LpNormG{\Lambda\l( f\r)}{p}
\end{equation*}
\end{thm}
\begin{proof}
Fix $l_1,l_2$ (recall, $l_1$ and $l_2$ just range over a finite set).  The theorem will
follow if we can show that for any sequence of $1$s and $-1$s, $\epsilon_{j,k}$,
we have that 
$$\sum_{j,k} \epsilon_{j,k} \Lambda_{j,k}^{l_1,l_2}$$
is bounded on $L^p$, uniformly in the choice of the sequence $\epsilon_{j,k}$
(and with a similar result for the $P_{j,k}^{l_1,l_2}$, which will follow in the same
way).
To see why this is enough, see p. 267 of \cite{SteinHarmonicAnalysis} and
Chapter 4, Section 5 of \cite{SteinSingularIntegralsAndDifferentiablilityPropertiesOfFunctions}.

However,
\begin{equation*}
\begin{split}
\sum_{j,k} \epsilon_{j,k} \Lambda_{j,k}^{l_1,l_2} &= \sum_{j,k} \Opt{\epsilon_{j,k} \dil{\phi_{l_1}}{2^j}\l(x\r) \dil{\psi_{l_2}}{2^k}\l(y\r)}\\
&= \Opt{\sum_{j,k} \epsilon_{j,k} \dil{\phi_{l_2}}{2^j}\l(x\r) \dil{\psi_{l_2}}{2^k}\l(y\r) }
\end{split}
\end{equation*}
and $$\sum_{j,k} \epsilon_{j,k} \dil{\phi_{l_1}}{2^j}\l(x\r) \dil{\psi_{l_2}}{2^k}\l(y\r)$$
converges to a product kernel, uniformly in the choice of $\epsilon_{j,k}$ (see \cite{NagelRicciSteinSingularIntegralsWithFlagKernels}, Theorem 2.2.1).
Hence, Corollary \ref{CorTwoSidedBounded} shows us that
$$\sum_{j,k} \epsilon_{j,k} \Lambda_{j,k}^{l_1,l_2}$$
is uniformly bounded on $L^p$.
\end{proof}

\begin{thm}
Suppose $T\in \sA'$.  Then, $T$ extends to a bounded operator $L^p\rightarrow L^p$, $1<p<\infty$.
\end{thm}
\begin{proof}
We first prove the result imagining that $\phi_l,\psi_l\in \sSz$.  At
the end, we will explain why it is enough to have them in $\sS_{(N)}$
for some large $N$.  This proof is more or less standard, however
we include it to help make clear where we are using $\phi_l,\psi_l\in \sS_{(N)}$.

Since $\phi_l,\psi_l\in \sSz$, we have (by Theorem \ref{ThmGensisElemKer})
that $\Lambda_{j,k}^{l_1,l_2},P_{j,k}^{l_1,l_2}$ are $2^j,2^k$ elementary kernels (uniformly
in $j,k$).  Hence, we have (for $f\in \sS$):
\begin{equation*}
\begin{split}
\Lambda_{j_1,k_1}^{l_1,l_2} T P_{j_2,k_2}^{l_1',l_2'} f &= \Lambda_{j_1,k_1}^{l_1,l_2} E_{2^{j_2},2^{k_2}}\\
&= 2^{-\l|j_1-k_1\r| - \l| j_2-k_2\r|} E_{2^{j_1},2^{k_1}} f\\
&\lesssim 2^{-\l|j_1-k_1\r| - \l| j_2-k_2\r|} \M f
\end{split}
\end{equation*}
where we have used the definition of $\sA'$, Lemma \ref{LemmaCompElemKer},
and Lemma \ref{LemmaElemBoundedByMax}, and $E_{2^j,2^k}$ just represents
some $2^j,2^k$ elementary kernel that may change from line to line.

Define $F_{j,k}^{l_1,l_2} = \Lambda_{j,k}^{l_1,l_2} T f$.  Then,
\begin{equation*}
\begin{split}
\l| F_{j,k}^{l_1,l_2} \r| &= \l| \sum_{l_1',l_2'}\sum_{j_1,k_1} \Lambda_{j,k}^{l_1,l_2} T P_{j_1,k_1}^{l_1',l_2'}\Lambda_{j_1,k_1}^{l_1',l_2'}f\r|\\
&\lesssim \sum_{l_1',l_2'} 2^{-\l|j_1-j\r| - \l|k_1-k\r|} \M\l( \Lambda_{j_1,k_1}^{l_1',l_2'} f\r)
\end{split}
\end{equation*}
and hence,
\begin{equation*}
\begin{split}
\l| F_{j,k}^{l_1,l_2}\r|^2 &\lesssim \l[\sum_{l_1',l_2'}\sum_{k_1,j_1} 2^{-\l|j_1-j\r| - \l|k_1-k\r|} \l( \M\l(\Lambda_{j_1,k_1}^{l_1',l_2'}f\r) \r)^2\r] \l[ \sum_{l_1',l_2'}\sum_{j_1,k_1} 2^{-\l|j_1-j\r| - \l|k_1-k\r|} \r]\\
&\approx \sum_{l_1',l_2'}\sum_{k_1,j_1} 2^{-\l|j_1-j\r| - \l|k_1-k\r|} \l( \M\l(\Lambda_{j_1,k_1}^{l_1',l_2'}f\r) \r)^2
\end{split}
\end{equation*}
and so, we have
\begin{equation*}
\begin{split}
\sum_{l_1,l_2}\sum_{j,k} \l| F_{j,k}^{l_1,l_2}\r|^2 \lesssim \sum_{l_1,l_2}\sum_{j,k} \l(\M\l( \Lambda_{j,k}^{l_1,l_2} f \r)\r)^2
\end{split}
\end{equation*}
Putting all of this together, we see:
\begin{equation*}
\begin{split}
\LpNormG{Tf}{p} &\approx \LpNormG{\l(\sum_{l_1,l_2}\sum_{j,k} \l|F_{j,k}^{l_1,l_2}\r|^2\r)^{\frac{1}{2}}}{p}\\
&\lesssim \LpNormG{\l( \sum_{l_1,l_2}\sum_{j,k} \l(\M\l(\Lambda_{j,k}^{l_1,l_2} f\r)\r)^2 \r)^{\frac{1}{2}}}{p}\\
&\lesssim \LpNormG{\Lambda\l(f\r)}{p}\\
&\approx \LpNormG{f}{p}
\end{split}
\end{equation*}
where we have used the vector valued maximal function, see \cite{SteinHarmonicAnalysis}, Chapter 2, Section 1.  The vector valued inequality comes from
(\ref{EqnBoundMaximal}), and the corresponding inequalities for $\sM_L,\sM_R$
as shown in \cite{SteinHarmonicAnalysis}.

Now we turn to explaining why we only need $\phi_l,\psi_l\in \sS_{(N)}$ for
some fixed large $N$.  Indeed, this proof used only a finite number of
the semi-norms that define the elementary operators.  This follows from
the fact that every proof we have done about elementary operators
was continuous.  For example Lemma \ref{LemmaCompElemKer} showed:
$$E_{2^{j_1},2^{k_1}}E_{2^{j_2},2^{k_2}}=2^{-\l|j_1-j_2\r|-\l|k_1-k_2\r|}\Et_{2^{j_1},2^{k_1}}$$
Where each semi-norm of $\Et$ was bounded in terms of a finite number
of semi-norms of the terms on the left hand side.  Thus, this proof
only required a finite number of semi-norms, which we may control
by taking $N$ large.  The only potential worry is the line where
we used the definition of $\sA'$ (ie, $T$ takes elementary
operators to elementary operators); since this was a definition,
and we do not have a priori continuity in the above sense.
However, this continuity follows from a combination of Theorem
\ref{ThmDecompOfsAp} and Lemma \ref{LemmaCompElemKer}.
\end{proof}

\section{Pseudolocality}\label{SectionPseudolocal}
	In this section, we show that the operators in $\sA'$ are pseudolocal,
and calculate bounds for derivatives of the kernel away from the diagonal;
although we will not put these bounds in a closed form.
In Section \ref{SectionHeisenbergGroup}, however, we will derive a closed
form for the growth of the kernel off the diagonal in the case that
$G$ is the three dimensional Heisenberg group.

Fix $T\in \sA'$.  Decompose $T$ as in Theorem \ref{ThmDecompOfsAp}:
\begin{equation*}
T=\sum_{j,k\in \Z} E_{2^j,2^k}
\end{equation*}
We will imagine this is a finite sum, and show that the result is
uniformly $C^\infty$ off of the diagonal.  It will then follow
that $T$ is pseudolocal.
In fact, we will prove the bounds separately for
\begin{equation*}
\begin{split}
T_L &= \sum_{k\geq j} E_{2^j, 2^k}\\
T_R &= \sum_{j\geq k} E_{2^j,2^k}
\end{split}
\end{equation*}
We focus on $T_L$, the bounds for $T_R$ being the same, with the roles
of right and left reversed.  Set $K\l(x ,z\r) = \ker{T_L}$.

Let us return to the notation from Section \ref{SectionStatmentOfResults}:
$$\lieg = V_1 \bigoplus \cdots \bigoplus V_n$$
It is easy to see that each vector field $\gr^{\alpha}$ with $\l|\alpha\r|=1$ can be written
in the form:
\begin{equation}\label{EqnDefnqj}
\sum_{j=1}^n q_{j,\alpha}\l( x\r)X_{L,\alpha,j}
\end{equation}
where $X_{L,\alpha,j}\in V_j$ (when thought of as left invariant vector fields),
and $q$ is a homogeneous polynomial of degree $j-1$.  For example, on the
Heisenberg group, $X_R=X_L-4y\partial_t$ (see Section \ref{SectionHeisenbergGroup} for this
notation).
\begin{lemma}\label{LemmaByPartsElem}
Let $\phi\l( x, z\r)$ be a $2^j,2^k$ elementary kernel (we are still
assuming $k\geq j$).  Then,
if $\l|\alpha_1\r|+\l|\alpha_2\r|=a$, $\l|\beta_1\r|+\l|\beta_2\r|=b$,
\begin{equation*}
\begin{split}
&\l|\glx^{\alpha_1} \grx^{\beta_1}\glz^{\alpha_2}\grz^{\beta_2} \phi\l(x,z\r) \r|
\\& \lesssim  \frac{2^{ja+kb}\l(1\wedge \l( \sum_{\l|\alpha\r|=1}\sum_{s=1}^n \l| q_{s,\alpha}\l(x\r) \r| 2^{sj-k} \r)\r)}{\l( 1+2^j\rlp{2^{j-k}}{x}{z} \r)^N \Vl{2^{j-k}}{x}{2^{-j}+\rlp{2^{j-k}}{x}{z}}}
\end{split}
\end{equation*}
where $N\geq 0$ is fixed and as large as we like.
\end{lemma}
\begin{proof}
It is easy reduce to the case when $a=0=b$, just from the definition
of an elementary kernel.  The case when $1=1\wedge \l( \sum_{\l|\alpha\r|=1}\sum_{s=1}^n \l| q_{s,\alpha}\l(x\r) \r| 2^{sj-k} \r)$ follows directly
from the definition of an elementary kernel.  For the other case,
consider (in what follows, $\psi$ with any subscript will denote a $2^j,2^k$
elementary operator):
\begin{equation*}
\begin{split}
\phi\l( x, z\r) &= \sum_{\l|\beta\r|=1} 2^{-k}\grx^\beta \psi_{\beta}\l(x,z\r)\\
&= \sum_{\l|\beta\r|=1} \sum_{s=1}^n 2^{-k} q_{s,\beta}\l( x\r) X_{L,\beta,s}\psi_{\beta}\l( x,z\r)\\
&= \sum_{\l|\beta\r|=1} \sum_{s=1}^n 2^{sj-k} q_{s,\beta}\l( x\r) \psi_{s,\beta}\l(x,z\r)
\end{split} 
\end{equation*}
now the claim follows by taking absolute values and applying the definition
of elementary kernels.
\end{proof}

\begin{thm}\label{ThmPseudolocality}
Let $a=\l|\alpha_1\r|+\l|\alpha_2\r|$.  Then,
for $x\ne z$, we have
\begin{equation*}
\begin{split}
&\l| \glx^{\alpha_1}\glz^{\alpha_2} K\l( x,z\r) \r|\\
&\lesssim \sum_{l=0}^\infty \rlp{2^{-l}}{x}{z}^{-a}\l( \frac{1}{\Vl{2^{-l}}{x}{\rlp{2^{-l}}{x}{z}}}\wedge \frac{\sum_{\l|\alpha\r|=1} \sum_{s=1}^n \l| q_{s,\alpha}\l(x\r) \r|  \rlp{2^{-l}}{x}{z}^{-s+1}  }{2^{l} \Vl{2^{-l}}{x}{\rlp{2^{-l}}{x}{z}}}\r)
\end{split}
\end{equation*}
and therefore $T_L$ is pseudolocal.  A similar result holds for $T_R$,
thereby showing that $T$ is pseudolocal.
\end{thm}
\begin{proof}
First let's see why this shows that $T_L$ is pseudolocal.  We claim
that the above sum converges absolutely.  This can be seen by using the facts
that $\rl{\epsilon}\geq \rl{1}$ and $V^L_{\epsilon} \geq V^L_{0}$ (for all
$\epsilon\in [0,1]$) and
thus if we only take the right part of the $\wedge$ we get a geometric
series.  Hence, the whole series converges absolutely, showing that
$T_L$ is pseudolocal.

Let $\phi_{2^j,2^k} = \ker{E_{2^j,2^k}}$.  Let $0\leq l=k-j$, and to
save space, define $\delta_{l} = \rlp{2^{-l}}{x}{z}$.  We think
of $l\geq 0$ as fixed, and sum over all those $k$, $j$ such that $k-j=l$.
Using Lemma \ref{LemmaByPartsElem}, we see:
\begin{equation*}
\begin{split}
&\l| \glx^{\alpha_1} \glz^{\alpha_2} \sum_{k-j=l} \phi_{2^j,2^k}\l(x,z\r)\r|\\
&\lesssim \sum_{k-j=l} 2^{ja}\frac{1\wedge \l(2^{-l}\sum \l| q_{s,\alpha}  \r| 2^{j\l( s-1 \r)} \r)}{\l(1+2^{j}\delta_l\r)^N \Vl{2^{-l}}{x}{2^{-j}+\delta_l}} \\
&= \sum_{j\in \Z} 2^{ja}\frac{1\wedge \l(2^{-l}\sum \l| q_{s,\alpha}  \r| 2^{j\l( s-1 \r)} \r)}{\l(1+2^{j}\delta_l\r)^N \Vl{2^{-l}}{x}{2^{-j}+\delta_l}} \\
\end{split}
\end{equation*}
where, in the numerator, $\sum = \sum_{\l|\alpha\r|=1} \sum_{s=1}^n$, and
we have suppressed the $x$ in $q_{s,\alpha}\l( x\r)$.
We separate the above sum into two sums: when $2^j\geq \frac{1}{\delta_l}$ and
when $2^j\leq \frac{1}{\delta_l}$.  Now,
\begin{equation*}
\begin{split}
\sum_{2^j\geq \frac{1}{\delta_l}} &\approx \sum_{2^j\geq \frac{1}{\delta_l}} 2^{ja} \frac{1\wedge \l(2^{-l}\sum \l| q_{s,\alpha}  \r| 2^{j\l( s-1 \r)} \r)}{\l(2^j \delta_l\r)^N \Vl{2^{-l}}{x}{\delta_l}}
\\& \approx \delta_l^{-a} \frac{1\wedge \l(2^{-l}\sum \l| q_{s,\alpha}  \r| \delta_l^{\l( 1-s \r)} \r)}{\Vl{2^{-l}}{x}{\delta_l}}
\end{split}
\end{equation*}
since the second term is a geometric sum (when $N$ is sufficiently large), and therefore bounded by its first term.
This is precisely the bound we were striving for.  We now turn to the sum
when $2^j\leq \frac{1}{\delta_j}$:
\begin{equation*}
\begin{split}
\sum_{2^j\leq \frac{1}{\delta_j}} &\approx \sum_{2^j\leq \frac{1}{\delta_j}} 2^{ja} \frac{1\wedge \l(2^{-l}\sum \l| q_{s,\alpha}  \r| 2^{j\l( s-1 \r)} \r)}{\Vl{2^{-l}}{x}{2^{-j}}} \\
&\approx \delta_l^{-a} \frac{1\wedge \l(2^{-l}\sum \l| q_{s,\alpha}  \r| \delta_l^{\l( 1-s \r)} \r)}{\Vl{2^{-l}}{x}{\delta_l}}
\end{split}
\end{equation*}
where this follows since the above sum is geometric, and therefore bounded
by its first term, completing the proof.
\end{proof}

\begin{rmk}
It seems likely that Theorem \ref{ThmPseudolocality} is the best we
can do (at least when $a=0$).  Indeed, if $K_1$ and $K_2$ are Calder\'on-Zygmund
kernels, we can decompose:
\begin{equation*}
\begin{split}
K_1 &= \sum_j \dil{\gl\cdot\phi_j}{2^j}\\
K_2 &= \sum_k \dil{\gr\cdot \psi_k}{2^k}
\end{split}
\end{equation*}
where $\phi_j,\psi_k$ are d-tuples of $\C_0^\infty$ functions supported on the unit ball $B$ (see \cite{NagelRicciSteinSingularIntegralsWithFlagKernels}).
Then we consider:
$$\Opl{K_1}\Opr{K_2} = \sum_{j,k} \Opl{\dil{\gl\cdot\phi_j}{2^j}}\Opr{\dil{\gr\cdot \psi_k}{2^k}}$$
One wishes to use the fact that $\gl\cdot \phi_j$ and $\gr \cdot \phi_k$
are derivatives of functions to yield a gain over the estimate
given in Theorem \ref{ThmBoundIntersection}.  The standard way
of doing this, when $k\geq j$, is to integrate the $\gr$ by parts over to
$\phi_j$.  However, this process is exactly the one we used in Lemma
\ref{LemmaByPartsElem}.  

The observant reader will note, however, that the bound in Theorem \ref{ThmPseudolocality}
is not actually symmetric in $x$ and $z$, as the optimal bound should be.  And
that, moreover, we could use the same proof to prove a seemingly better
symmetric bound.
This turns out to not be an essential point, and indeed the bound is essentially
symmetric in $x$ and $z$.  This is exemplified in Section \ref{SectionHeisenbergGroup} in the case of the Heisenberg group.
Thus, without some new idea, one is
unable to do better than Theorem \ref{ThmPseudolocality}.
\end{rmk}

	\subsection{The Heisenberg Group}\label{SectionHeisenbergGroup}
		In this section, we derive a closed form for the bound in Theorem
\ref{ThmPseudolocality}, in the case of the three dimensional
Heisenberg group, $\Ho$.  As a manifold $\Ho=\C\times \R$, and
we give it coordinates $\l( z,t\r)= \l( x,y,t\r)$.
The multiplication is given by $\l( z,t\r) \l( w,s\r) = \l( z+w, t+s+2{\rm Im}\l(z\overline{w}\r)\r)$.
The dilation is given by $r\l( z,t\r) = \l( rz, r^2t\r)$.  The
left invariant vector fields of order $1$ are spanned by $X_L= \partial_x+2y\partial_t$, $Y_L = \partial_y-2x\partial_t$, while the right invariant
vector fields of order $1$ are spanned by $X_R=\partial_x-2y\partial_t$, $Y_R=\partial_y+2x\partial_t$.
We also have:
\begin{equation*}
\l[ X_L,Y_L\r] = -4\partial_t = -\l[X_R,Y_R\r] 
\end{equation*}
and so $\partial_t$ spans the left (and right) invariant vector fields of
order $2$.

We fix $\l(z,t\r),\l(w,s\r)\in \Ho$, $\l(z,t\r)\ne \l(w,s\r)$, and we again define (for $l\geq 0$),
$\delta_l = \rlp{2^{-l}}{\l(z,t\r)}{\l(w,s\r)}$.  Note that
$\delta_\infty = \rl{0}$, $\delta_0=\rl{1}$.  Fix $\alpha\in \N$, $\alpha\geq 0$
($\alpha$ will play the role of $a$ in Theorem \ref{ThmPseudolocality}),
and let $\zeta=\l( z,t\r)$.
We will show:
\begin{thm}\label{ThmHeisBound}
\begin{equation}\label{EqnHeisBound}
\sum_{l=0}^\infty \delta_l^{-\alpha}\l( \frac{1}{\Vl{2^{-l}}{\zeta}{\delta_l}} \wedge \frac{1+\l|z\r| \delta_l^{-1}}{2^l\Vl{2^{-l}}{\zeta}{\delta_l}}\r) \approx \frac{1}{\delta_\infty^2 \delta_0^{2+\alpha}} 
\end{equation}
Note that this sum is exactly the one that appears in Theorem \ref{ThmPseudolocality}, in the case of $\Ho$.
\end{thm}
The first question to address is:  When do we use each side of the $\wedge$?
Namely, we are interested in the question, when is $\rl{2^{-l}} =\delta_l \gtrsim 2^{-l}\l|z\r|$?  So let us investigate the question:  When is $\epsilon \l|z\r| \lesssim \rl{\epsilon}$?  The answer is that it is true precisely when
$\rl{\epsilon} \approx \rl{0}$.  Indeed, suppose we are on the scale
$\rl{\epsilon} \approx \delta$, here $\delta$ is just some number $>0$, not
to be confused with $\delta_l$. (We assume $\delta<1$, and that
we are working very close to $0$ and then extend the results by
homogeneity.)  Then, by Theorem \ref{ThmCanUseMaximalDet}, we wish
to find the maximal determinant among $3\times 3$ submatricies of:
\begin{equation*}
\begin{bmatrix}
\delta & 0 & -2y\delta\\
0 & \delta & 2x\delta\\
0 & 0 & \delta^2 \\
\delta\epsilon & 0 & 2y\delta\epsilon \\
0 & \delta \epsilon & -2x\delta\epsilon\\
0 & 0 & \delta^2 \epsilon^2
\end{bmatrix}
\end{equation*}
Now the largest three determinants are given by:  Rows $(1,2,3) = \delta^4$,
$(1,2,4)= 4 \delta^3 y \epsilon$, and $(1,2,5)= -4\delta^3x\epsilon$.
Thus, when $\delta\geq \l|z\r|\epsilon$, $(1,2,3)$ is the largest determinant (up to a constant),
and therefore on this scale $\rl{\epsilon}\approx \rl{0}$ (since the
first $3$ rows corresponded to the left invariant vector fields).
\begin{rmk}
Actually, this proof extends to an arbitrary stratified group.  That is,
the right part of the $\wedge$ in Theorem \ref{ThmPseudolocality}
is less than the left part precisely when $\rl{2^{-l}}\approx \rl{0}$.
We leave the details to the interested reader.
\end{rmk}

\begin{prop}\label{PropEquivHeisDist}
Suppose $\rlp{\epsilon}{\l( z,t\r)}{\l( w,s\r)}\lesssim \l|z\r| \epsilon$.
Then,
$$\rlp{\epsilon}{\l(z,t \r)}{\l( w,s\r)} \approx \l| z-w \r|+ \frac{1}{\epsilon\l|z\r|}\l|t-s+2{\rm Im}\l(z\overline{w} \r) \r|$$
\end{prop}
\begin{proof}
When $\epsilon=1$, the result follows easily.
Moreover, if $1\geq \epsilon\geq\frac{1}{2}$, we have
$\rl{\epsilon} \approx \rl{1}$, and the result follows from the
case when $\epsilon =1$.  Henceforth,
we restrict our attention to the case $\epsilon < \frac{1}{2}$.
Fix $\delta\lesssim \l|z\r| \epsilon$. 
We will show that the following conditions are equivalent:
\begin{enumerate}
\item $\rlp{\epsilon}{\l(z,t\r)}{\l(w,s\r)}\lesssim \delta$.
\item $\exists w_0,s_0$ such that $\l| w-w_0\r| \lesssim \delta$, $\l| s_0-s+2{\rm Im}\l( w_0\overline{w} \r) \r|\lesssim \delta^2$, $\l| z-w_0\r| \lesssim \epsilon \delta$, $\l| t-s_0\r|\lesssim \epsilon \l| z\r|\delta$.
\item $\exists s_0$ such that $\l| z-w\r| \lesssim \delta$, $\l| s_0-s +2{\rm Im}\l( z\overline{w} \r)\r|\lesssim \delta^2$, $\l| t-s_0\r|\lesssim \delta \l|z\r| \epsilon$.
\item $\exists s_0$ such that $\l| z-w\r| \lesssim \delta$, $\l| s_0-s +2{\rm Im}\l( z\overline{w} \r)\r|\lesssim \delta\l|z\r| \epsilon$, $\l| t-s_0\r|\lesssim \delta \l|z\r| \epsilon$.
\item $\l| z-w\r| \lesssim \delta$, $\l| t-s+2{\rm Im} \l(z\overline{w}\r)\r|\lesssim \delta \l| z\r| \epsilon$.
\item $\l| z-w \r| + \frac{1}{\epsilon \l| z\r|} \l| t-s+2{\rm Im}\l( z\overline{w} \r)\r|\lesssim \delta$.
\end{enumerate}
this will complete the proof, since the statement of the proposition is
$1\Leftrightarrow 6$.  We will show the equivalences in the following
order:
\begin{equation*}
1\Leftrightarrow 2, \quad 3\Rightarrow 2\Rightarrow 4\Rightarrow 3, \quad 4\Leftrightarrow 5, \quad 5 \Leftrightarrow 6
\end{equation*}
To see that $1\Leftrightarrow 2$, we apply Theorem \ref{ThmNewDist} to see:
$$\rl{\epsilon} = \l(\epsilon \rr{0} \r) \circ \l(\epsilon \rl{0} \r)\circ \l( \l(1-\epsilon\r) \rl{0} \r) =  \l( \epsilon \rl{1}\r)\circ \l( \l( 1-\epsilon\r) \rl{0} \r)  \approx  \l( \epsilon \rl{1}\r)\circ \rl{0}$$
where we have used that $\epsilon<\frac{1}{2}$, and the obvious
fact that $\l( \epsilon \rl{0}\r) \circ \l(\l(1-\epsilon\r)\rl{0}\r)=\rl{0}$.  The statement
$1\Leftrightarrow 2$ now follows directly from the definition of
$\l( \epsilon \rl{1}\r)\circ \rl{0}$.
For $3\Rightarrow 2$ merely take $w_0=z$.  Suppose we have $2$.  Consider,
$$\l| z-w\r|\leq \l| z-w_0\r| + \l| w_0-w\r| \lesssim \delta+ \epsilon \delta \lesssim \delta$$
\begin{equation*}
\begin{split}
\l| s_0-s +2{\rm Im}\l( z\overline{w} \r) \r| &\leq \l|2{\rm Im}\l( \l( z-w_0 \r) \overline{w}\r) \r| + \l| s_0-s + 2{\rm Im}\l( w_0\overline{w}\r) \r|
\\&\lesssim \delta^2 + \l| {\rm Im}\l( \l( z-w_0 \r)\overline{w} \r) \r|
\\& = \delta^2 +\l| {\rm Im}\l( \l( z-w_0 \r) \l(\overline{ w-w_0+z }\r) \r) \r|
\\& \leq \delta^2 + \l| z-w_0\r| \l| w-w_0\r| + \l| z-w_0\r| \l|z\r|
\\& \lesssim \delta^2 + \epsilon \delta^2 + \l|z\r| \epsilon \delta
\\& \lesssim \l|z\r| \epsilon \delta
\end{split}
\end{equation*}
and we therefore have $2\Rightarrow 4$.  Suppose we have $4$.  Define
$s_1  = s-2{\rm Im}\l( z\overline{w} \r)$.  We will show that $3$
holds with $s_1$ in place of $s_0$ (the notation $s_0$ has already
been used in the definition of $4$).  Indeed,
$$\l| s_1-s+2{\rm Im}\l(z\overline{w}\r) \r| =0 \lesssim \delta^2$$
$$\l| t-s_1 \r|=\l| t-s+2{\rm Im}\l( z\overline{w} \r) \r|\leq \l| s_0-s+2{\rm Im}\l( z\overline{w} \r)\r| + \l| s_0-t \r|\lesssim \l|z\r| \delta \epsilon$$
and so $3$ holds.

$5\Rightarrow 4$ follows just by taking $s_0=t$.  Suppose we have $4$, consider:
$$\l| t-s+2{\rm Im} \l(z\overline{w} \r) \r| \leq \l| s_0-s+2{\rm Im}\l(z\overline{w}\r) \r|+ \l| t-s_0\r|\lesssim \l|z\r| \delta \epsilon$$
and so $5$ holds.  Finally, $6\Leftrightarrow 5$ is obvious.
\end{proof}

We now claim that the condition $\epsilon \l|z\r|\lesssim \rl{\epsilon}$
is the same as the condition $\epsilon \l| z\r| \lesssim \rl{0}$.  
Indeed, $\rl{\epsilon}\leq \rl{0}$, and so one direction is clear.
For the other direction, fix $\eta>0$ and consider those $\epsilon$
such that $\epsilon \l| z\r| \geq \eta \rl{\epsilon}$.  The left
hand side decreases to $0$ as $\epsilon$ decreases, while the right
hand side increases.  Thus there is a least $\epsilon$ (call it
$\epsilon_0$) for which it holds.  For this $\epsilon_0$, we have
$\epsilon_0 \l| z\r| = \eta \rl{\epsilon_0}$, so by the above
remarks, we have $\rl{\epsilon_0}\approx \rl{0}$, and thus,
$\epsilon_0\l| z\r| \geq C \rl{0}$.  Hence for all $\epsilon$
such that $\epsilon \l| z\r| \geq \eta \rl{\epsilon}$ (namely
all $\epsilon\geq \epsilon_0$), we have $\epsilon\l| z\r| \geq C \rl{0}$.

We separate the sum (\ref{EqnHeisBound}) into two parts:  when $2^l \delta_\infty\leq \l| z\r|$ and when $2^l\delta_\infty\geq \l|z\r|$.  We first look at the
case when $2^l \delta_\infty\leq \l| z\r|$, so that by Proposition \ref{PropEquivHeisDist}, we have
\begin{equation*}
\delta_l \approx \l|z-w\r| + 2^l\frac{\l| t-s+2{\rm Im}\l(z\overline{w}\r) \r|}{\l| z\r|}
\end{equation*}
and we see by Theorem \ref{ThmEstBalls} and the remarks above that
\begin{equation*}
\Vl{2^{-l}}{\zeta}{\delta_l}\approx 2^{-l}\l| z \r|\delta_l^3
\end{equation*}
to save space, we denote $b=\l| z-w\r|$ and $a= \frac{\l| t-s+2{\rm Im}\l(z\overline{w}\r) \r|}{\l|z\r|}$, so that
\begin{equation*}
\delta_l \approx b+2^la
\end{equation*}
Also, let $c=\frac{\l| z\r|}{\delta_\infty}$.  Thus, we consider:
\begin{equation*}
\begin{split}
&\sum_{2^l\delta_\infty\leq \l|z\r|} \delta_l^{-\alpha}\l( \frac{1}{\Vl{2^{-l}}{\zeta}{\delta_l}} \wedge \frac{1+\l|z\r| \delta_l^{-1}}{2^l\Vl{2^{-l}}{\zeta}{\delta_l}}\r) \\
&\approx\sum_{2^l\delta_\infty\leq \l|z\r|} \delta_l^{-\alpha} \frac{1}{\Vl{2^{-l}}{\zeta}{\delta_l}}  \\
\end{split}
\end{equation*}
\begin{equation*}
\begin{split}
&\approx \sum_{2^l\leq c} \frac{2^l}{\l|z\r| \l( b+2^l a \r)^{3+\alpha}}\\
&\approx  \frac{1}{\l| z\r| a^{3+\alpha}} \sum_{2^l\leq c} \frac{2^l}{\l( \frac{b}{a} +2^l \r)^{3+\alpha}}\\
&\approx \frac{1}{\l| z\r| a^{3+\alpha}} \int_0^{\log_2\l( c\r)}\frac{2^t}{\l(\frac{b}{a}+2^t\r)^{3+\alpha}} dt\\
&\approx \frac{1}{\l| z\r| a^{3+\alpha}} \int_1^c \frac{1}{\l(\frac{b}{a}+u\r)^{3+\alpha}} du
\end{split}
\end{equation*}
\begin{equation*}
\begin{split}
\approx \frac{1}{\l| z\r| a^{3+\alpha}} \frac{\l( c+\frac{b}{a} \r)^{2+\alpha} - \l( 1+\frac{b}{a} \r)^{2+\alpha}}{\l(c+\frac{b}{a} \r)^{2+\alpha}\l(1+\frac{b}{a}\r)^{2+\alpha}}
\end{split}
\end{equation*}
Recall, $l\geq 0$ and this sum is nonzero only when $c\geq 2^l$.  In
fact, let us ignore the possibility that this sum is nonzero when
$c\leq 2$ (this case can be taken care of in a similar manner
to our sum when $c\leq 2^l$).
Thus, we have:
\begin{equation*}
\begin{split}
&\approx \frac{1}{\l|z\r| a^{3+\alpha}} \frac{c \l( \sum_{j=1}^{2+\alpha} \l(\frac{b}{a}\r)^{2+\alpha-j}c^{j-1} \r)}{\l(c+\frac{b}{a} \r)^{2+\alpha}\l(1+\frac{b}{a}\r)^{2+\alpha}} \\
&\approx \frac{1}{\l|z\r| a^{3+\alpha}} \frac{c\l( c+\frac{b}{a} \r)^{2+\alpha-1}}{\l(c+\frac{b}{a} \r)^{2+\alpha}\l(1+\frac{b}{a}\r)^{2+\alpha}}
\end{split}
\end{equation*}
\begin{equation*}
\begin{split}
&\approx \frac{1}{\l|z\r| a^{3+\alpha}} \frac{c}{\l(c+\frac{b}{a} \r)\l(1+\frac{b}{a}\r)^{2+\alpha}}\\
&= \frac{1}{\delta_\infty \l( \frac{\l|z\r|}{\delta_\infty}a+b \r)\l( a+b \r)^{2+\alpha}}\\
&\approx \frac{1}{\delta_\infty \l( \frac{\l|z\r|}{\delta_\infty}a+b \r)\delta_0^{2+\alpha}}
\end{split}
\end{equation*}
where, in the last line, we have used that $\l( a+b\r) \approx \delta_0$,
by Proposition \ref{PropEquivHeisDist}.  Finally, we will be done with this sum provided we can
show
$$\frac{\l|z\r|}{\delta_\infty}a+b\approx \delta_\infty$$
To see this, consider $\epsilon_0$ such that $\epsilon_0 \l|z\r| = \rl{\epsilon_0}$ (as before).  Then, for this $\epsilon_0$, we have $\epsilon_0 \l|z\r|=\rl{\epsilon_0} \approx \rl{0}$.  However, we also have $\epsilon_0 \l|z\r| \gtrsim \rl{\epsilon_0}$,
and therefore,
\begin{equation*}
\delta_\infty = \rl{0} \approx \rl{\epsilon_0} \approx \l( b+\frac{a}{\epsilon_0}\r) \approx \l( b+ \frac{\l| z\r|}{\rl{0}} a\r) =\l( b+ \frac{\l| z\r|}{\delta_\infty} a\r)
\end{equation*}
completing the proof for this sum.

We now turn to the sum when $2^l\delta_\infty \geq \l|z\r|$.  In this case,
$\delta_l\approx \delta_\infty$, and $\Vl{2^{-l}}{\zeta}{\delta_l}\approx \Vl{0}{\zeta}{\delta_l}$, by the remarks at the beginning of the proof.
And thus, we are considering:
\begin{equation*}
\begin{split}
&\sum_{2^l\delta_\infty\geq \l|z\r|} \delta_l^{-\alpha}\l( \frac{1}{\Vl{2^{-l}}{\zeta}{\delta_l}} \wedge \frac{1+\l|z\r| \delta_l^{-1}}{2^l\Vl{2^{-l}}{\zeta}{\delta_l}}\r) \\
&\approx \sum_{2^l \delta_\infty\geq\l|z\r|} \delta_\infty^{-\alpha} \l[ \frac{\l|z\r|}{2^l\delta_\infty\Vl{0}{\zeta}{\delta_\infty}} + \frac{1}{2^l \Vl{0}{\zeta}{\delta_\infty}} \r] \\
&\approx \sum_{2^l \delta_\infty \geq \l|z\r|} \frac{\l|z\r|}{2^l \delta_\infty^{5+\alpha}} + \frac{1}{2^l \delta_\infty^{4+\alpha}}
\end{split}
\end{equation*}
the first sum is geometric, and therefore bounded by its first term, and we have:
\begin{equation*}
\begin{split}
&\approx \frac{1}{\delta_\infty^{4+\alpha}} + \sum_{2^l\delta_\infty\geq \l|z\r|} \frac{1}{2^l \delta_\infty^{4+\alpha}}
\end{split}
\end{equation*}
The second term above is bounded by:
\begin{equation*}
\sum_{l\geq 0} \frac{1}{2^l\delta_\infty^{4+\alpha}} \approx \frac{1}{\delta_\infty^{4+\alpha}}
\end{equation*}
Hence, the whole sum is:
\begin{equation*}
\approx \frac{1}{\delta_\infty^{4+\alpha}}
\end{equation*}
Since $\delta_\infty\geq \delta_0$, this completes the proof.

\section{Closing Remarks}\label{SectionClosingRemarks}
	Definition \ref{DefnsA} only tested high derivatives of the operator $T$.
One can replace Definition \ref{DefnsA} with an equivalent definition
that works for derivatives of all orders, but with the price that
the $B$ must be quite a bit more complicated.  Let $q_{j,\alpha}^R$
be the functions from (\ref{EqnDefnqj}) and $q_{j,\alpha}^L$ be
the corresponding ones with the roles of left and right reverse.
We then define:
\begin{equation*}
\begin{split}
&\Bt\l( 2^{j_0}, 2^{k_0}, x,y, N_L,N_R\r)\\
&= \sum_{j\leq j_0, k\leq k_0} 2^{jN_L+kN_R} \l( 1\wedge\l( \sum_{\l|\alpha\r|=1}\sum_{s=1}^N \l|q_{s,\alpha}^R\l(x\r)\r|2^{sj-k} \r)\wedge \l(\sum_{\l|\alpha\r|=1}\sum_{s=1}^N \l|q_{s,\alpha}^L\l(x\r) \r|2^{sk-j}\r) \r)\\
&\quad\quad\quad\times K_{2^j,2^k}\l( x,y\r)
\end{split}
\end{equation*}
where $K_{2^j,2^k}$ is the function from Section \ref{SectionRelToConvo}.
Note that $\Bt$ does not involve $m$.  Then Definition \ref{DefnsA} with
$B$ replaced by $\Bt$ now works for all $N_L,N_R\geq 0$, and defines
the same algebra.  That these algebras are the same follow in a manner
similar to the bounds in the rest of this paper.

One may think of the operators in $\sA$ and $\sA'$ as ``smoothing of order $0$.''
To come up with an analogous definition for operators which
are ``smoothing of order $s_L$'' in the left invariant vector fields and 
``smoothing of order $s_R$'' in the right invariant vector fields, it
suffices to modify Definitions \ref{DefnsA} and \ref{DefnsAp} only
slightly.

Indeed, we say an operator $T:\sS_{(N)}\rightarrow \sS'$ (for some large $N$)
is in $\sA_{s_L,s_R}$ if it satisfies the conditions of Definition \ref{DefnsA}
with $B\l( \cdot,\cdot, N_L, N_R,\cdot,\cdot,\cdot\r)$ replaced by
$B\l( \cdot, \cdot, N_L-s_L,N_R-s_R,\cdot,\cdot,\cdot\r)$.

We say an operator $T:\sSz\rightarrow \sSz$ is in $\sA_{s_L,s_R}'$ if
$r_L^{s_L}r_R^{s_R} T E_{r_L,r_R}$ is an $r_L,r_R$ elementary operator
for every $r_L,r_R$ elementary operator $E_{r_L,r_R}$, uniformly
in the relevant parameters.
Then, $\sA_{s_L,s_R} = \sA_{s_L,s_R}'$ (under the obvious identification).
It is clear that if $T_1\in \sA_{s_L^1,s_R^1}', T_2\in \sA_{s_L^2,s_R^2}'$,
then we have $T_1T_2 \in \sA_{s_L^1+s_L^2, s_R^1+s_R^2}'$, and therefore
we have a similar result for $\sA_{s_L,s_R}$ (remember, we are just thinking
of these operators on $\sS_{(N)}$ for $N$ large).  Many of the results
of this paper extend to these operators in the obvious way.

Finally, let us consider the question of whether or not it is really
necessary to have a cancellation condition on both sides simultaneously
as in Definition \ref{DefnsA}, as opposed to something more along
the lines of the standard definitions of Calder\'on-Zygmund operators.  
One could think about this in two ways.  One could try to use a one
sided cancellation condition along the lines of Lemma \ref{LemmasAPointwise}
(or something slightly stronger, in terms of the $\Bt$ above), along with
a growth condition of the kernel of $T$ off of the diagonal.  However, in light of
Theorem \ref{ThmHeisBound}, any condition along these lines seems likely to be
necessarily weaker than our Definition \ref{DefnsA}.

Alternatively, let us go back to considering the composition of
$$\Opl{K_1}\Opr{K_2}$$
where $K_1,K_2$ are Calder\'on-Zygmund kernels.  We decompose
\begin{equation*}
\begin{split}
K_1 &= \sum_{j\in \Z} \dil{\phi_j}{2^j}\\
K_2 &= \sum_{k\in \Z} \dil{\psi_k}{2^k}
\end{split}
\end{equation*}
where $\phi_j,\psi_k$ form a bounded subset of $C_0^\infty$, are supported
in the unit ball, and have mean $0$ (see \cite{NagelRicciSteinSingularIntegralsWithFlagKernels}, Theorem 2.2.1).  This cancellation condition on the $\phi_j$
essentially tells us that if $\eta$ is another $C_0^\infty$ function supported
on the unit ball, we have:
\begin{equation*}
\begin{split}
\Opl{K_1} \Opl{\dil{\eta}{2^{j_0}}} &= \sum_{j\leq j_0} \Opl{\dil{\phit_j}{2^j}}
\end{split}
\end{equation*}
where the $\phit_j$ are essentially of the same form as the $\phi_j$, and
with a similar result for $\Opl{\dil{\eta}{2^{j_0}}}\Opl{K}$.  Thus,
if $\eta_1,\eta_2$ are of the same form as $\eta$, we have:
\begin{equation*}
\Opl{\dil{\eta_1}{2^{j_1}}}\Opl{K_1} \Opr{\dil{\eta_2}{2^{j_2}}} = \sum_{j\leq j_1\wedge j_2} \Opl{\dil{\phit_j}{2^j}}
\end{equation*}
Hence, for composition $\Opl{K_1}\Opr{K_2}$ we want:
\begin{equation*}
\begin{split}
&\Opl{\dil{\eta_1}{2^{j_1}}} \Opr{\dil{\eta_2}{2^{k_1}}} \Opl{K_1}\Opr{K_2} \Opl{\dil{\eta_3}{2^{j_2}}} \Opr{\dil{\eta_4}{2^{k_2}}}
\\&\quad = \sum_{\substack{j\leq j_1\wedge j_2\\k\leq k_1\wedge k_2}} \Opl{\dil{\phit_j}{2^j}}\Opr{\dil{\psit_k}{2^k}}
\end{split}
\end{equation*}
and so a cancellation condition on one side alone will be fine in the case
when $k_1\leq k_2$ and $j_1\leq j_2$ (or the reverse situation), but seems like
it will not be able to yield the desired estimate when $k_1<k_2$ and $j_2<j_1$
(or the reverse situation).

\bibliographystyle{amsalpha}

\bibliography{convo}

\end{document}